\documentclass[10pt]{article}

\usepackage{amsmath,amsfonts,amsthm}
\usepackage[text={6.75in,10in},centering,a4paper]{geometry}
\usepackage[colorlinks=true,linkcolor=blue,citecolor=blue]{hyperref}
\usepackage{color}

\newtheorem{Thm}{Theorem}

\newtheorem{Lem}{Lemma}

\newtheorem{Def}{Definition}
\newtheorem{Rem}{Remark}

\def\Ran{\operatorname{\rm Ran}}

\def\cl{\operatorname{\rm cl}}
\def\span{\operatorname{\rm span}}

\def\wU{\widehat{U}}
\def\wW{\widehat{W}}
\def\wA{\widehat{A}}
\def\wM{\widehat{M}}

\begin{document}

\title {Perturbations of embedded eigenvalues for the planar bilaplacian}

\author{%
Gianne Derks\\
Department of Mathematics\\
University of Surrey\\
Guildford, GU2 7XH, UK
\and
Sara Maad Sasane\footnote{Sara Maad Sasane was partially supported by
  the European Union under the Marie Curie Fellowship
  MEIF-CT-2005-024191 and by the Swedish Research Council.}\\
Department of Mathematics\\
Stockholm University\\
106 91 Stockholm, Sweden
\and
Bj\"orn Sandstede\footnote{Bj\"orn Sandstede was partially supported by a Royal Society--Wolfson Research Merit Award and by the NSF under grant DMS-0907904.}\\
Division of Applied Mathematics\\
Brown University\\
Providence, RI~02912, USA
}

\date{\today}
\maketitle

\begin{abstract}
Operators on unbounded domains may acquire eigenvalues that are embedded in the essential spectrum. Determining the fate of these embedded eigenvalues under small perturbations of the underlying operator is a challenging task, and the persistence properties of such eigenvalues is linked intimately to the multiplicity of the essential spectrum. In this paper, we consider the planar bilaplacian with potential and show that the set of potentials for which an embedded eigenvalue persists is locally an infinite-dimensional manifold with infinite codimension in an appropriate space of potentials.
\end{abstract}


\begin{section}{Introduction}

Determining the dependence of the spectrum of operators on perturbations is an important issue that is of relevance in many applications. Of course, much is known in this direction: the persistence of point eigenvalues and the behavior of the essential spectrum under small bounded perturbations, for instance, have been analyzed comprehensively, and we refer to \cite{tK76} for many results along these lines. Here, we consider differential operators that are posed on unbounded domains and are interested in the interaction between eigenvalues, with proper eigenfunctions in the underlying domain of the operator, and the essential spectrum. More precisely, we study the fate of eigenvalues that are embedded in the essential spectrum under small perturbations of the operator. Typically, such eigenvalues will disappear under generic perturbations of the potential, and it is therefore of interest to determine the class of perturbations for which an embedded eigenvalue persists. For the bilaplacian on cylindrical domains, we showed in our previous work \cite{DMS08} that the set of perturbations for which an embedded eigenvalue persists is an infinite-dimensional manifold of finite codimension. Furthermore, we showed that the codimension of this set is given by the multiplicity of the essential spectrum, defined as the number of independent continuum eigenfunctions or, more rigorously, via the spectral resolution of the Fourier transform of the bilaplacian (see e.g.\ \cite[Definition~2 in \S85]{AG}). In this paper, we continue the investigation that we began in \cite{DMS08} and consider the bilaplacian posed on the plane: the challenge is that the essential spectrum of the planar bilaplacian has infinite multiplicity. Thus, we may expect that the set of potentials for which an embedded eigenvalue persists is an infinite-dimensional manifold of infinite codimension, and this is indeed what we shall prove for an appropriate class of potentials.

Before stating our results, we briefly outline why embedded
eigenvalues are of interest. Our first motivation comes from quantum
mechanics: the eigenfunctions associated with eigenvalues of an energy
operator correspond to bound states that can be attained by the
physical system modelled by the energy operator. If such an eigenvalue
is embedded in the essential spectrum, then its fate under
perturbations of the potential determines whether the associated bound
states persists or not (see \cite{HS96,RS78} for examples). A second
example comes from inverse scattering theory, where eigenvalues
correspond to coherent soliton structures of the underlying integrable
system, while the essential spectrum describes radiative scattering
behavior. Thus, bifurcations of solitons are reflected by the
disappearance or persistence of embedded eigenvalues
\cite{PS98,PS00}. Finally, embedded eigenvalues provide a common
mechanism for the destabilization of travelling waves in
near-integrable Hamiltonian partial differential equations, and we
refer to \cite{S02} for further background information and pointers to
the literature.

As mentioned above, we focus in this paper on the persistence of embedded eigenvalues for the planar bilaplacian. Our primary reason for considering the bilaplacian is that this operator is complex enough to exhibit the underlying difficulties, while not adding technical complications that have nothing to do with the issue we are interested in. In other words, the planar bilaplacian provides a useful paradigm for the issues that we expect to encounter for other more complicated differential operators. Note also that the applications we mentioned above all involve selfadjoint operators, a feature shared by the bilaplacian.

We now describe the precise setting that we consider. Let $r_0>0$, and assume that $\theta\in C^\infty_0(B_{r_0}(0);\mathbb{R})$ is a radially symmetric potential. Hence, we use polar coordinates $(r,\varphi)$, write $\theta=\theta(r)$, and consider the multiplication operator on $L^2(\mathbb R^2)$ (also denoted by $\theta$) defined by
\begin{equation*}
   [\theta u](r,\varphi) := \theta(r)u(r,\varphi).
\end{equation*}
We define $\mathcal L:= \Delta^2 + \theta$ on $L^2(\mathbb R^2)$,
where $\Delta^2$ is the bilaplace operator which is densely defined
on $L^2(\mathbb R^2)$ with domain $H^4(\mathbb R^2)$. It is known
that the spectrum of $\Delta^2$ is $\sigma(\Delta^2) = [0,\infty)$.
Since $\theta$ has compact support, the essential spectra of
$\mathcal L$ and $\Delta^2$ coincide, and so $\sigma_c(\mathcal L)=[0,\infty)$. We assume that $\theta$ is chosen so that $\mathcal
L$ has a simple positive eigenvalue $\lambda_0$:
\begin{enumerate}
\item[(A1)] $\mathcal L$ has an eigenvalue $\lambda_0>0$ of multiplicity $1$.
\end{enumerate}
We are mainly interested in the case where $\lambda_0$ is an embedded
eigenvalue, i.e.\ when $\lambda_0{\ge} 0$, since when $\lambda_0$ is
isolated from the rest of the spectrum, the persistence of
eigenvalues is well known,~\cite[pp. 213--215]{tK76}. We also exclude
the case $\lambda_0=0$ which lies on the boundary between spectrum and resolvent set.
We denote by $u_*(r,\varphi)$ the eigenfunction associated with the embedded eigenvalue $\lambda_0$. Since $\theta$ is radially symmetric and the laplacian $\Delta$ is invariant under rotations of the underlying cartesian coordinates, we see that the functions $u_*(r,\varphi+\varphi_0)$ are, for each fixed $\varphi_0$, also eigenfunctions of $\mathcal{L}$ belonging to the eigenvalue $\lambda_0$. The simplicity of $\lambda_0$ required in assumption~(A1) therefore implies that $u_*$ is a radial function, and we henceforth write $u_*=u_*(r)$. It is clear by existence and uniqueness of solutions of ODEs that $u_*(r)$ cannot vanish for all $r\ge r_1$, and so we assume that 
\begin{enumerate}
 \item [(A2)] $r_1>r_0$ is such that $u_*(r_1)\ne 0$.
\end{enumerate}
Lemma~\ref{ex1} below shows that our hypotheses can be satisfied. We now perturb the potential $\theta$ by potentials $\rho$ in the weighted $L^2$-space $\mathcal R:= L^2([0,r_1],H^{1/2}(S^1),r\,dr)$ of functions that map the interval $[0,r_1]$ into $H^{1/2}(S^1)$, where the interval $[0,r_1]$ is the domain of the radial variable $r$, while $H^{1/2}(S^1)$ describes the dependence on the angular variable $\varphi$. Our main result is as follows.

\begin{Thm}\label{T:main}
  Let $0<r_0\leq r_1$, $\theta\in C_0^\infty(B_{r_0}(0);\, \mathbb R)$
  be radially symmetric, and assume that~(A1) and~(A2) hold.  Then
  there exists $\delta>0$ and a neighbourhood $\mathcal O$ of $0$ in
  $\mathcal R= L^2([0,r_1], H^{1/2}(S^1),r\, dr)$ such that the set
   \begin{equation*}
     \mathcal R_{emb} := \{\rho\in \mathcal O;\;
     \mathcal L + \rho \text{ has an embedded eigenvalue in }
     (\lambda_0 - \delta,\lambda_0 + \delta)\}
   \end{equation*}
   is a smooth manifold in $\mathcal R$ of infinite dimension and codimension.
\end{Thm}

Before commenting on the ideas behind the proof of Theorem~\ref{T:main}, we illustrate that our hypotheses can be met.

\begin{Lem}\label{ex1}
There exists a smooth radial potential $\theta(r)$ with compact support such that $\mathcal{L}=\Delta^2+\theta$ satisfies Hypothesis~(A1).
\end{Lem}

\begin{proof}
Let $K_0(r)$ denote the modified Bessel function of the second kind and define a smooth, strictly positive function $u_0(r)$ via
\[
u_0(r) =
\left\{\begin{array}{lcl}
1      &\quad& 0\leq r\leq1 \\
K_0(r) &\quad& 2\leq r
\end{array}\right.
\]
together with a smooth interpolation in the intermediate region $r\in[1,2]$. Note that $u_0$ decays to exponentially as $r\to\infty$ and can be chosen so that $u_0(r)>0$ for all $r$. Thus, the radial potential
\[
\theta := \frac{1}{u_0} (-\Delta^2+1) u_0 = \frac{1}{u_0} (\Delta+1)(-\Delta+1) u_0 =
\left\{\begin{array}{lcl}
1 && 0\leq r\leq1 \\
0 && 2\leq r
\end{array}\right.
\]
is well-defined, smooth, and has support in $[0,2]$ since $K_0(r)$ satisfies $(-\Delta+1)K_0=0$. Furthermore, we have
\[
\mathcal{L} u_0 = (\Delta^2+\theta) u_0 = \Delta^2 u_0 + \frac{1}{u_0} [(-\Delta^2+1) u_0] u_0 = u_0,
\]
and $u_0$ is a positive radial eigenfunction belonging to the embedded eigenvalue $\lambda_0=1$ of $\mathcal{L}$.

It remains to show that $\lambda_0=1$ is simple. Using the radial symmetry of $\theta$, the results presented in the rest of this paper imply that it suffices to show that the equation
\begin{eqnarray}\label{e:ex1}
\left[\partial_r^2 + \frac{1}{r}\partial_r - \frac{k^2}{r^2}\right]^2 u = (1-\theta) u
\quad\mbox{ for } r\in(0,2) \\ \nonumber
u_r(0) = u_{rrr}(0)=0, \quad u(r)=K_k(r) \mbox{ for } r\geq2
\end{eqnarray}
does not have a solution $u(r)$ for each integer $k\geq1$. We will now outline why (\ref{e:ex1}) will not have solutions for $k\geq1$ provided $\theta$ is modified appropriately but omit the straightforward details. Using variations of parameters, it can be shown that (\ref{e:ex1}) cannot have solutions for $k\gg1$. If it does have solutions for some or all of the remaining finitely many integers $k\geq1$, then we can modify the potential $\theta$ to remove these solutions while retaining the eigenfunction for $k=0$. Indeed, any solution of (\ref{e:ex1}) for $k\geq1$ is of the explicit form $u(r)=r^\ell$ for some integer $\ell=\ell(k)\geq1$ since $\theta(r)-1=0$ for $0\leq r\leq1$. Replacing $u_0$ in the above construction of $\theta$ by $u_0+\epsilon v_0$ for bounded functions $v_0$ with support in $(\frac12,1)$ and using the necessary expressions (\ref{E:F_0j-eqn}) derived in section~\ref{S:Lyapunov-Schmidt} for the persistence of eigenvalues, it is then not difficult to see that any nonzero choice of $v_0\geq0$ removes the solutions of (\ref{e:ex1}) for $k\geq1$ for $0\ll\epsilon\ll1$.
\end{proof}

The idea for proving Theorem~\ref{T:main} is to characterize embedded
eigenvalues as roots of a regular function, since such a
characterization would allow us to use the implicit function theorem.
As it appears difficult to find a functional-analytic characterization
of embedded eigenvalues, we pursue here a dynamical-systems
formulation similar to that used in our precursor work \cite{DMS08}
for the bilaplacian on cylinders. The eigenvalue problem can be
written as a system of differential equations in the radial evolution
variable $r$ posed on an appropriate function space $X$ of functions
that are defined in the angular variable $\varphi$. The issue is that
this system is ill-posed in the sense that, for given initial data,
solutions may not exist. Using a similar approach as in Scheel
\cite{aS98}, we will show, however, that this dynamical system has an
exponential dichotomy: there are two infinite-dimensional subspaces
$X^{cu}$ and $X^s$ of $X$ at $r=r_1$ so that solutions with initial
data in $X^{cu}$ exist and stay bounded for $r\leq r_1$, while
solutions with data in $X^s$ exist and decay as $r\to\infty$. The
intersection of these spaces corresponds to eigenfunctions of the
underlying operator, and our goal is therefore to characterize those
perturbations for which this intersection is nontrivial. We show that
there are infinitely many conditions that characterize such
intersections and prove that we can solve them using an implicit
function theorem.  {A key issue is the space for the perturbation
  $\rho$. For the conditions of the implicit function theorem to be
  satisfied, the space for $\rho$ needs to be
  $L^2([0,r_1];H^{1/2}(S^1), r\, dr)$, a space with very low
  regularity.  This low regularity forces us to work with different
function spaces for $r\le r_1$ (where $\rho$ has its support)
and for $r\ge r_1$ (where we have an explicit formulation of the
solutions of the system in terms of Bessel functions), and so we need
to take extra care when matching the solutions at $r=r_1$.}

The rest of this paper is organized as follows. In section~\ref{s:sd},
we introduce the spatial-dynamics formulations of the eigenvalue
problem. In sections~\ref{S:-infty} and \ref{S:r_small}, we prove the
existence of exponential dichotomies for the bilaplacian and for the
operator $\mathcal{L}$, respectively, near the core $r=0$. We then
construct dichotomies for $\mathcal{L}$ in the far field for $r\gg1$
in section~\ref{S:s large} and discuss similar properties for the
adjoint spatial dynamical system in section~\ref{S:adjoint}.  These
results are then used in section~\ref{S:Lyapunov-Schmidt} where we
match the solutions from the core and the far field 
{by using Lyapunov--Schmidt reduction and prove Theorem~\ref{T:main}.
  The paper is concluded with suggestions for extensions and some open
  problems.}
\end{section}


\begin{section}{Spatial-dynamics formulation}\label{s:sd}
If $\lambda$ is an eigenvalue of $\mathcal L + \rho$, then there
exists $u\in H^4(\mathbb R^2)$ such that
\begin{equation}\label{E:eigenvalue_equation}
   \Delta^2 u + (\theta + \rho) u = \lambda u.
\end{equation}
Let $r_3>r_2>\max(1,r_1)$. We introduce a new radial variable
\begin{equation}\label{E:s-var}
   s(r) =  \begin{cases}
      \log r & \text{if }r\le r_2, \\
      r & \text{if }r\ge r_3,
   \end{cases}
\end{equation}
and for $r\in (r_2,r_3)$, we define $s$ such that $s\in
C^\infty(\mathbb R^+;\mathbb R)$ is strictly increasing. Note that this implies
that there exist constants $c$ and $C$ such that $0<c<C$ and $c\le s'(r) \le C$
for every $r_2\le r \le r_3$.
We define $\tilde \theta$ and $\tilde \rho$ by $\tilde
\theta(s(r)) = \theta(r)$, etc.  Since $s$ is an increasing function,
it is invertible, and we denote the inverse function by $r(s)$.
Let $s_j := s(r_j)$, $j = 1,\dots,3$. Under the coordinate
transformation \eqref{E:s-var}, the space $\mathcal R$ transforms into
the space $\tilde {\mathcal R}$ given by
\begin{equation*}
   \tilde {\mathcal R} :=  L^2((-\infty,s_1];\; H^{1/2}(S^1), e^{2s} ds),
\end{equation*}
that is, the weighted $L^2$ space with values in $H^{1/2}(S^1)$ and weight
$e^{2s}$.

Setting $v=\Delta u$, equation \eqref{E:eigenvalue_equation}
is equivalent to the system
\begin{equation}\label{E:intermediate}
   \begin{aligned}
      \Delta u &= v, \\
      \Delta v &= (\lambda - \tilde \theta - \tilde\rho) u,
   \end{aligned}
\end{equation}
where in the variables $s$ and $\varphi$, the Laplacian is given by
\[
\Delta = \frac{1}{r'(s)^2}\,
\left[\,\frac{\partial^2}{\partial s^2} +
  \left(\frac{r'(s)}{r(s)}-\frac{r''(s)}{r'(s)}\right)
  \frac\partial{\partial s} + \left(\frac{r'(s)}{r(s)}\right)^2
  \frac{\partial^2}{\partial \varphi^2}\,\right].
\]
Rewriting this intermediate system as a first order system, with
$u_1=u$, $u_2=u'$, $u_3=v$ and $u_4=v'$, where $'$ denotes
differentiation with respect to $s$, we obtain a system of the form
\begin{equation}\label{E:perturbed system}
   U'(s) = A(s;\lambda,\tilde \rho) U,
\end{equation}
where $A(s;\lambda,\tilde \rho)$ is given by
\begin{equation}\label{E:A-matrix}
   A(s;\lambda,\tilde \rho) := \left(
   \begin{matrix}
       0 & 1 & 0 & 0 \\
       -\frac{r'(s)^2}{r(s)^2}\partial^2 &
       \frac{r''(s)}{r'(s)} - \frac{r'(s)}{r(s)} & r'(s)^2 & 0 \\
       0 & 0 & 0 & 1 \\
       (\lambda - \tilde \theta - \tilde \rho) r'(s)^2 & 0
       & - \frac{r'(s)^2}{r(s)^2}\partial^2 & \frac{r''(s)}{r'(s)} -
       \frac{r'(s)}{r(s)}
   \end{matrix}
   \right)
\end{equation}
where $\partial$ denotes differentiation with respect to $\varphi$,
i.e., $\partial=\frac{\partial}{\partial \varphi}$.  The expression
for $A(s;\lambda,\tilde \rho)$ simplifies significantly for $s<s_2$ or
$s>s_3$; see sections~\ref{S:-infty} and~\ref{S:s large}.

The perturbation $\tilde\rho\in \tilde{\mathcal R}$ is, in general, not
continuous or even bounded, so we need to study weak solutions of
\eqref{E:perturbed system}. Let
\begin{eqnarray*}
X &=& H^2(S^1)\times H^1(S^1)\times H^1(S^1)\times L^2(S^1);\\
Y&=& H^3(S^1)\times H^2(S^1)\times H^2(S^1)\times H^1(S^1).
\end{eqnarray*}
\begin{Def}\label{D:weak sol}
Let $J$ be an interval of $\mathbb R$.  A function $U:J\to X$ is a
weak solution of \eqref{E:perturbed system} in $J$ if
\begin{enumerate}
 \item $U\in L^2_{loc}(J;Y)\cap H^1_{loc}(J;X)$,
 \item for every $V\in C_0^\infty(J;X)$ we have
 \begin{equation*}
  -\int_J U(s) V'(s)\, ds = \int_J A(s;\lambda,\tilde \rho) U(s) V(s) \, ds.
 \end{equation*}
 \end{enumerate}
\end{Def}

\begin{Lem}\label{L:equation-system}
  Let $\lambda\in \mathbb R$. The eigenvalue
  equation~\eqref{E:eigenvalue_equation} has a solution $u\in
  H^4_{loc}(\mathbb R^2)$ if and only if \eqref{E:perturbed system}
  has a weak solution $U\in H^1_{loc}(\mathbb R;X)\cap
  L^2_{loc}(\mathbb R;Y)\cap L^\infty(\mathbb R_-;X)$.
\end{Lem}
\begin{proof}
  Suppose that $u\in H^4_{loc}(\mathbb R^2)$ is a solution of
  \eqref{E:eigenvalue_equation}, and let $U:= (u,u',\Delta u,(\Delta
  u)')^T$, where $'$ denotes differentiation with respect to $s$.

  We first consider $u$ as a function of $r$, and let $B_R(0)$ be a
  ball centered at $0$, with $R$ any positive radius. Then
   \begin{equation*}
      \begin{aligned}
         u&\in H^1((0,R);H^3(S^1), r\, dr)\cap L^2((0,R);H^4(S^1), r\, dr) \\
         &\qquad\qquad\subset
         H^1((0,R);H^2(S^1),r\, dr)\cap L^2((0,R);H^3(S^1), r\, dr), \\
         u'&= r' \frac{du}{dr}\in H^1((0,R);H^2(S^1),r\, dr)\cap
         L^2((0,R);H^3(S^1), r\, dr) \\
         &\qquad \qquad\subset H^1((0,R);H^1(S^1), r\, dr)\cap
         L^2((0,R);H^2(S^1),r\, dr), \\
         \Delta u &\in H^1((0,R);H^1(S^1),r\, dr)\cap
         L^2((0,R);H^2(S^1),r\, dr), \\
         (\Delta u)' &= r' \frac{d(\Delta u)}{dr} \in
         H^1((0,R);L^2(S^1),r\, dr)\cap L^2((0,R);H^1(S^1),r\, dr),
      \end{aligned}
   \end{equation*}
   where $r'=dr/ds = r$ for $r<r_2$.  By the Sobolev embedding
   theorem, $U\in C([0,R];X)$, and so $U(s(r))$ has a limit as $r\to
   0+$, or equivalently, as $s\to -\infty$. Hence, viewing $U$ as a
   function of $s$, $U\in L^\infty(\mathbb R_-;X)$. We also see that
   $U\in H^1_{loc}(\mathbb R;X)\cap L^2_{loc}(\mathbb R;Y)$.  It is
   clear from the construction that $U$ is a weak solution of
   \eqref{E:perturbed system}.

   Conversely, let $U\in H^1_{loc}(\mathbb R;X)\cap L^2_{loc}(\mathbb
   R;Y)\cap L^\infty(\mathbb R_-;X)$ be a weak solution of
   \eqref{E:perturbed system}, and let $u=U_1$. Viewing $u$ as a
   function of $r$ rather than of $s$, it is clear that $u\in
   H^1_{loc}(\mathbb R_+;H^2(S^1))\subset C((0,\infty);C(S^1))$, and
   so by \eqref{E:intermediate},
   \begin{equation*}
      \Delta^2 u = (\lambda -  \theta - \rho) u\in
      L^2_{loc}(\mathbb R_+;H^{1/2}(S^1),r\, dr)\subset
      L^2_{loc}(\mathbb R_+;L^2(S^1), r\, dr) =  L^2_{loc}(\mathbb R^2\setminus\{0\}).
   \end{equation*}
   Since we also have
   \begin{equation*}
      \begin{aligned}
         u&\in L^\infty([0,r_2];H^2(S^1))\subset L^\infty([0,r_2];C(S^1)), \\
         \lambda - \theta - \rho&\in L^2([0,r_2];H^{1/2}(S^1),r\,dr)
         \subset L^2([0,r_2];L^2(S^1),r\, dr),
      \end{aligned}
   \end{equation*}
   we see that
   $\Delta^2 u = (\lambda - \theta - \rho) u \in
   L^2([0,r_2];L^2(S^1),r\, dr) = L^2(B_{r_2}(0))$. We have proved that
   $u\in H^4(B_{r_2}(0))\cap H^4_{loc}(\mathbb R^2\setminus \{0\})
   = H^4_{loc}(\mathbb R^2)$.

   Since it is also clear that $u$ solves
   \eqref{E:eigenvalue_equation}, the proof is complete.
\end{proof}

Note that a weak solution satisfies $U\in C(\mathbb R;X)$ (see
e.g. \cite[p. 286]{lcE98}), and so the following definition for an
exponential dichotomy makes sense (see also~\cite{wC78} for the
standard definition for ODEs and~\cite{PSS97} for an extension to PDEs):
\begin{Def}\label{D:dichotomies}
  Let $J$ be an unbounded subinterval of $\mathbb R$.  We say that
  equation \eqref{E:perturbed system} has an exponential dichotomy in
  $X$ on $J$ if there exists a family of projections $P(s)$ for $s\in
  J$ such that for any $s\in J$, $P(s)\in \mathcal L(X)$, $P(s)^2 = P(s)$
  and $P(\cdot) U \in C(J;X)$ for every $U\in X$, and there
  exist constants $K>0$ and $\kappa^s<\kappa^u$ with the following
  properties:
   \begin{enumerate}
   \item[(i)] For each $t\in J$ and $U\in X$ there exists a unique weak
     solution $\Phi^s(s,t)U$ of \eqref{E:perturbed system} defined for
     $s\ge t$, $s, t\in J$ such that $\Phi^s(t,t) U = P(t) U$ and
      \begin{equation*}
         \|\Phi^s(s,t) U\|_X \le K e^{\kappa^s (s - t)}\|U\|_{X}
      \end{equation*}
      for all $s\ge t$, $s, t\in J$.
    \item[(ii)] For each $t\in J$ and $U\in X$ there exists a unique
      weak solution $\Phi^{u}(s,t) U$ of \eqref{E:perturbed system}
      defined for $s\le t$, $s, t\in J$ such that $\Phi^{u}(t,t) U = (I - P(t))U$ and
      \begin{equation*}
         \|\Phi^{u}(s,t)U\|_X \le K e^{\kappa^u (s - t)}\|U\|_X
      \end{equation*}
      for all $s\le t$, $s, t\in J$.
    \item[(iii)] The solutions $\Phi^s(s,t)U$ and $\Phi^u(s,t)U$
      satisfy
      \begin{equation*}
         \begin{array}{rll}
            \Phi^s(s,t) U &\in \Ran P(s) \qquad &
            \text{for every }s\ge t,\;  s, t\in J, \\
            \Phi^u(s,t) U &\in \ker P(s) \qquad &
            \text{for every }s\le t,\; s, t\in J.
         \end{array}
      \end{equation*}
   \end{enumerate}
\end{Def}

{We also need the definition of time-dependent exponential dichotomy, which will be used for $J=[s_1,\infty)$ and
with $\mathcal X^s:= H^1\times L^2\times H^1\times L^2$ with the $s$-dependent norm
\begin{equation*}
   \|U\|_{\mathcal X^s}^2:= \frac{1}{s^2}\|u_1\|_{H^1(S^1)}^2 + \|u_1\|_{L^2(S^1)}^2 + \|u_2\|_{L^2(S^1)} + \frac{1}{s^2}\|u_3\|_{H^1(S^1)}^2 + \|u_3\|_{L^2(S^1)}^2 + \|u_4\|_{L^2(S^1)},
\end{equation*}
where $u_j$ are the components of $U$, $j=1,\dots, 4$.}
{\begin{Def}
   Let $J$ be an unbounded subinterval of $\mathbb R$.  We say that
  equation \eqref{E:perturbed system} has a time-dependent exponential dichotomy in
  $\mathcal X^s$ on $J$ if there exists a family of projections $P(s)$ for $s\in
  J$ such that for any $s\in J$, $P(s)\in \mathcal L(\mathcal X^s)$, $P(s)^2 = P(s)$
  and $P(\cdot) U \in C(J;\mathcal{X})$ for every $U\in \mathcal{X}$, and there
  exist constants $K>0$ and $\kappa^s<\kappa^u$ with the following
  properties:
   \begin{enumerate}
   \item[(i)] For each $t\in J$ and $U\in \mathcal X^t$ there exists a unique
     solution $\Phi^s(s,t)U$ of \eqref{E:perturbed system} defined for
     $s\ge t$, $s, t\in J$ such that $\Phi^s(t,t) U = P(t) U$ and
      \begin{equation*}
         \|\Phi^s(s,t) U\|_{\mathcal X^s} \le K e^{\kappa^s (s - t)}\|U\|_{\mathcal X^t}
      \end{equation*}
      for all $s\ge t$, $s, t\in J$.
    \item[(ii)] For each $t\in J$ and $U\in \mathcal X^t$ there exists a unique
      solution $\Phi^{u}(s,t) U$ of \eqref{E:perturbed system}
      defined for $s\le t$, $s, t\in J$ such that $\Phi^{u}(t,t) U = (I - P(t))U$ and
      \begin{equation*}
         \|\Phi^{u}(s,t)U\|_{\mathcal X^s} \le K e^{\kappa^u (s - t)}\|U\|_{\mathcal X^t}
      \end{equation*}
      for all $s\le t$, $s, t\in J$.
    \item[(iii)] The solutions $\Phi^s(s,t)U$ and $\Phi^u(s,t)U$
      satisfy
      \begin{equation*}
         \begin{array}{rll}
            \Phi^s(s,t) U &\in \Ran P(s) \qquad &
            \text{for every }s\ge t,\;  s, t\in J, \\
            \Phi^u(s,t) U &\in \ker P(s) \qquad &
            \text{for every }s\le t,\; s, t\in J.
         \end{array}
      \end{equation*}
   \end{enumerate}
\end{Def}}

In the following sections, we will consider the intervals
\begin{equation*}
J_-:=(-\infty,s_1] \mbox{~~and~~ } J_+:= [s_1,\infty),
\end{equation*}
and show that the system~\eqref{E:perturbed system} has {an exponential
dichotomy on $J_-$ and a time-dependent exponential dichotomy on $J_+$.}
\end{section}


\begin{section}{Dichotomies for the system at $-\infty$}\label{S:-infty}
For $s\le s_1$, $r(s)=e^s$, and hence the system \eqref{E:perturbed
    system} is given by
\begin{equation}\label{E:core equation}
   \left\{
   \begin{aligned}
      u_1' &= u_2, \\
      u_2' &= -\partial^2 u_1 +e^{2s} u_3, \\
      u_3' &= u_4, \\
      u_4' &= \left(\lambda - \tilde \theta(s) -
        \tilde \rho(s,\cdot)\right)e^{2s}u_1 - \partial^2 u_3.
   \end{aligned}
   \right.
\end{equation}

In the limit as $s\to -\infty$ we have the system
\begin{equation*}
   \left(
   \begin{matrix}
      u'_1 \\
      u'_2 \\
      u'_3 \\
      u'_4 \\
   \end{matrix}
   \right) = \left(
   \begin{matrix}
      0 & 1 & 0 & 0 \\
      -\partial^2 & 0 & 0 & 0 \\
      0 & 0 & 0 & 1 \\
      0 & 0 & -\partial^2 & 0
   \end{matrix}
   \right)\left(
   \begin{matrix}
      u_1 \\
      u_2 \\
      u_3 \\
      u_4
   \end{matrix}
   \right),
\end{equation*}
or
\begin{equation}\label{E:-infty-system}
   U' = A_- U.
\end{equation}

We expand $U = (u_1,u_2,u_3,u_4)^T$ as a Fourier series in the
$\varphi$ variable, and denote the $k$th Fourier coefficient by
$\widehat U_k(s)$.  For $j\in \mathbb R$, we define the weighted $l^2$
spaces $l^2_j$ with norm defined by
   \begin{equation*}
      \|\{a_k\}_{k\in \mathbb Z}\|_{l^2_j}^2 :=
      \sum_{k\in \mathbb Z} (1 + k^2)^j |a_k|^2
   \end{equation*}
The function space induced by~$X$ is
\begin{equation}\label{E:X hat}
\widehat X := l^2_2\times l^2_1\times l^2_1\times l^2.
\end{equation}
The system~\eqref{E:-infty-system} decouples in the Fourier space and
for $k\in \mathbb Z$ we have
\begin{equation}\label{E:kth Fourier system}
   \widehat U_k'(s)= \widehat A_{-}(k)\widehat U_k(s),
\end{equation}
where
\begin{equation*}
   \widehat A_{-}(k) := \left(
   \begin{matrix}
      0 & 1 & 0 & 0 \\
      k^2 & 0 & 0 & 0 \\
      0 & 0 & 0 & 1 \\
      0 & 0 & k^2 & 0
   \end{matrix}
   \right).
\end{equation*}
The eigenvalues of $\widehat A_{-}(k)$ are $\pm |k|$, and for $k\ne 0$
both eigenvalues have geometric multiplicity $2$.  The eigenvectors
for $k\ne 0$ are $(\pm 1/k^2,1/|k|,0,0)^T$ and
$(0,0,\pm 1/|k|,1)$ (we normalize the eigenvectors so that their
$\widehat X$ norm is approximately constant and bounded away from $0$
as $k\to\infty$) . Let
\begin{equation*}
   M_k := \left(
   \begin{matrix}
      -1/k^2 & 0 & 1/k^2 & 0 \\
      1/|k| & 0 & 1/|k| & 0 \\
      0 & -1/|k| & 0 & 1/|k| \\
      0 & 1 & 0 & 1
   \end{matrix}
   \right)
\end{equation*}
and
\begin{equation*}
  D_k := \left(
  \begin{matrix}
      -|k| & 0 & 0 & 0 \\
      0 & -|k| & 0 & 0 \\
      0 & 0 & |k| & 0 \\
      0 & 0 & 0 & |k|
   \end{matrix}
   \right),
\end{equation*}
so that $\widehat A_{-}(k) = M_k D_k M_k^{-1}$ for $k\ne 0$. Note that
\begin{equation*}
   M_k^{-1} = \frac{1}{2}\left(
   \begin{matrix}
      -k^2 & |k| & 0 & 0 \\
      0 & 0 & -|k| & 1 \\
      k^2 & |k| & 0 & 0 \\
      0 & 0 & |k| & 1
   \end{matrix}
   \right),
\end{equation*}

For $k=0$, the eigenvalue $0$ has algebraic multiplicity $4$ and
geometric multiplicity $2$, so $\widehat A_{-}(0)$ is not
diagonalizable. Note however that $\widehat A_{-}(0)$ is already in
Jordan normal form. We define $M_0 = I$ and $D_0 = \widehat A_{-}(0)$.

\begin{Lem}\label{L:closed}
  The operator $A_{-}:X\to X$ is a closed densely defined operator
  with spectrum $\sigma(A_{-}) = \mathbb Z$.
\end{Lem}
\begin{proof}
Recall that $X = H^2(S^1)\times H^1(S^1)\times H^1(S^1)\times
L^2(S^1)$, and $Y=H^3(S^1)\times H^2(S^1)\times H^2(S^1)\times
H^1(S^1)$.  It is easy to check that the domain of $A_{-}$ is $Y$,
which is dense in $X$. To see that $A_{-}$ is closed, let $U_j \in
Y$ be such that $U_j\to U$ in $X$ and $A_{-} U_j\to f$ in
$X$. We write $U_j = (u_{1,j},u_{2,j},u_{3,j},u_{4,j})^T$ etc. By
the definition of $A_{-}$ we have
   \begin{equation*}
      \begin{array}{r l l l}
         u_{1,j}&\to & u_1 \quad &\text{in }H^2, \\
         u_{2,j}&\to & u_2 \quad &\text{in }H^1, \\
         u_{3,j}&\to & u_3 \quad &\text{in }H^1, \\
         u_{4,j}&\to & u_4 \quad &\text{in }L^2,
      \end{array}
   \end{equation*}
   while
   \begin{equation*}
      \begin{array}{r l l l}
         u_{2,j}&\to & f_1 \quad &\text{in }H^2, \\
         -\partial^2 u_{1,j}&\to & f_2 \quad &\text{in }H^1, \\
         u_{4,j}&\to & f_3 \quad &\text{in }H^1, \\
         -\partial^2 u_{3,j}&\to & f_4 \quad &\text{in }L^2.
      \end{array}
   \end{equation*}
   It follows that $u_2 = f_1\in H^2$, and that $u_{1,j}$ converges in
   $H^3$. Since $u_{1,j}\to u_1$ in $H^2\supset H^3$, and since limits
   (in $H^2$) are unique if they exist, we also have $u_{1,j}\to u_1$
   in $H^3$, and so $-\partial^2 u_1=f_2$.  It follows in
   exactly the same way that $u_4 = f_3\in H^1$, that $u_3\in H^2$ and
   that $-\partial^2 u_3 = f_4$.  This shows that $U\in Y$
   and $A_{-} U = F$, and so $A_{-}:X\to X$ is closed.

   The operator $A_{-}:X\to X$ induces an operator $\widehat
   A_{-}:\widehat X\to \widehat X$ defined by
   \begin{equation*}
      (\widehat A_{-} \widehat U)_k := \widehat A_{-}(k) \widehat U_k.
   \end{equation*}
   Then $\widehat A_{-}$ is a densely defined operator on $\widehat X$
   with domain $\widehat Y := l^2_3\times l^2_2\times l^2_2\times
   l^2_1$.

   It is clear that $(A_{-} - \mu I):X\to X$ has a bounded inverse if
   and only if $(\widehat A_{-} - \mu I):\widehat X\to \widehat X$ has
   a bounded inverse. It is also clear that $k\in \sigma(\widehat
   A_{-})$ for $k\in \mathbb Z$.  To prove that there are no other
   points in the spectrum of $A_{-}$, let $\mu\in \mathbb C\setminus
   \mathbb Z$.

   Define $\widehat M:l^2\times l^2\times l^2\times l^2\to l^2_2\times
   l^2_1\times l^2_1\times l^2$ by
   \begin{equation*}
      (\widehat M \widehat U)_k = M_k \widehat U_k,
   \end{equation*}
   and note that $\widehat M$ is a linear homeomorphism between these
   spaces. Define also the unbounded operator $\widehat D$ on
   $l^2\times l^2\times l^2\times l^2$ by
   \begin{equation*}
      (\widehat D \widehat U)_k = D_k \widehat U_k.
   \end{equation*}
   Note that $\widehat D$ is a closed densely defined operator with
   domain $l^2_1\times l^2_1\times l^2_1\times l^2_1$, and that
   $\sigma(\widehat D) = \mathbb Z$.

   If $\mu\in \mathbb C\setminus\mathbb Z$, then
   \begin{equation*}
      (\widehat A_{-} - \mu I)^{-1} = \widehat M (\widehat D - \mu I)^{-1} \widehat M^{-1}
   \end{equation*}
   It is now easy to see that $(\widehat A_{-}-\mu I)^{-1}:\widehat
   X\to \widehat X$ is bounded, and consequently also $(A_{-}-\mu
   I)^{-1}:X\to X$.
\end{proof}

Having established that the spectrum of $A_{-}$ consists exactly of
its eigenvalues, we define the (generalized) spectral projections
$P^s$, $P^c$, $P^u$ in $X$, corresponding to the negative, the zero
and the positive eigenvalues of $A_{-}$, respectively. Let $X^s = P^s
X$, etc. so that $X = X^s\oplus X^c\oplus X^u$, where $X^s$ and $X^u$
are infinite-dimensional whereas $X^c$ is four-dimensional. {We also
define corresponding spectral projections $P_k^s$, $P_k^u$ of $\widehat A_-(k)$,
in the spaces $X_k$, $k\in \mathbb Z\setminus\{0\}$ and note that if
$U=\sum_{k\in \mathbb Z}\widehat U_k e^{ik\cdot}\in X$, then
\begin{equation*}
   \begin{aligned}
      P^s U &= \sum_{k\in\mathbb Z\setminus\{0\}} P_k^s \widehat U_k e^{ik\cdot},\\
      P^u U &= \sum_{k\in\mathbb Z\setminus\{0\}} P_k^u \widehat U_k e^{ik\cdot},\\
      P^c U &= \widehat U_0.
   \end{aligned}
\end{equation*}}

{\begin{Lem}\label{L:resolvent-estimate}
  The operator $A_{-}$ possesses an exponential dichotomy
  in $X$ on $J_-=(-\infty,s_1]$ with constant $K$ and rates $\kappa^s=0$ and $\kappa^u=1$, and another exponential dichotomy in $X$ on $J_-$
  with constant $K$ and rates
  $\kappa^s=-1$ and $\kappa^u=0$.
\end{Lem}}
\begin{proof}
  {Let $\eta\in (0,1)$ be arbitrary. We apply Lemma 2.1 of \cite{PSS97} for the operators $A_- -\eta I$ and $A_- + \eta I$, and obtain exponential dichotomies with constant $K$ and rates $\kappa^s=-\eta$ and $\kappa^u=1-\eta$, and $\kappa^s=-1+\eta$ and $\kappa^u=\eta$, respectively. The existence of exponential dichotomies for $A_-$ with rates $\kappa^s=-1$ and $\kappa^u=0$, and $\kappa^s=0$ and $\kappa^u=1$, respectively then follows by using the transformation $V=e^{\pm \eta \cdot} U$.}

  {We only consider the operator  $A_{-} -\eta I$, since the proof for $A_- + \eta I$ is similar.}
  The result follows from Lemma 2.1 of \cite{PSS97} if we can verify
  condition (H1) of \cite{PSS97} for the operator $A_{-} -\eta I$,
  namely
   \begin{enumerate}
      \item[(H1)] Suppose that there exists a constant $C>0$ such that
      \begin{equation*}
         \bigl\|\left(A_{-} - \eta I -i\mu I\right)^{-1}\bigr\|
         _{\mathcal L(X)}\le \frac{C}{1 + |\mu|}
      \end{equation*}
      for every $\mu\in \mathbb R$.
   \end{enumerate}

   As in the proof of Lemma~\ref{L:closed}, it suffices to prove that
   there exists a constant $\tilde C$ such that
   \begin{equation*}
      \left\|\left(\widehat D - \eta I - i\mu I\right)^{-1}
        \widehat U\right\|_{l^2\times l^2\times l^2\times l^2}
      \le \frac{\tilde C}{1 + |\mu|}
      \|\widehat U\|_{l^2\times l^2\times l^2\times l^2}
   \end{equation*}
   for every $\widehat U\in l^2\times l^2 \times l^2\times l^2$.
   Note that for $k\ne 0$ we have
   \begin{equation*}
      \left|\left(D_k  - \eta I - i\mu I\right)^{-1} \widehat
        U_k\right|^2
      = \frac{|\widehat U_k|^2}{(k - \eta)^2 + \mu^2}\le
      \frac{2}{\min(\eta^2,(1+\eta)^2)(1 + |\mu|)^2}|\widehat U_k|^2,
   \end{equation*}
   and it is not difficult to see that a similar estimate holds for
   $|(D_0 - \eta I - i\mu I) \widehat U_0|^2$.  Hence
   \begin{equation*}
      \begin{aligned}
        \left\|\left(\widehat D - \eta I - i \mu I\right)^{-1}
          \widehat U\right\|_{l^2\times l^2\times l^2\times l^2}^2 &        \sum_{k\in \mathbb Z}
        \left| \left(D_k - \eta I - i\mu I\right)^{-1}\widehat U_k\right|^2  \\
        &\le \frac{\tilde C^2}{(1 + |\mu|)^2} \sum_{k\in\mathbb
          Z}|\widehat U_k|^2.
      \end{aligned}
   \end{equation*}
\end{proof}


\end{section}


\begin{section}{{Dichotomies near the core}}\label{S:r_small}
The system~\eqref{E:core equation} can be abbreviated and written as
\begin{equation}\label{E:core equation'}
   U' = (A_- + B(s;\lambda,\tilde\rho)) U,
\end{equation}
where
\begin{equation}
  \label{E:B in J-}
  B(s;\lambda,\tilde \rho) := e^{2s}\,
  \begin{pmatrix}
    0&0&0&0\\
    0&0&1&0\\
    0&0&0&0\\
    \lambda - \tilde\theta(s) - \tilde\rho(s,\cdot)&0&0&0
  \end{pmatrix}
\end{equation}
We will show that the system~\eqref{E:core
  equation'} has an exponential dichotomy on the interval
$J_-=(-\infty,s_1]$.

To show this, we would like to apply Theorem~1
of~\cite{PSS97}. This is not possible, however, since $\tilde \rho\in
\tilde{\mathcal R}$ is not smooth enough in $s$. We are interested in
$\tilde \rho$ small and consider therefore first
$\tilde\rho=0$, and show that the $\lambda$-perturbed system
\begin{equation}\label{E:approximated core equation'}
   U' = \bigl(A_{-} + B(s;\lambda, 0)\bigr) U
\end{equation}
possesses an exponential dichotomy in $X$ on $J_-$. Then we will use the
implicit function theorem to show that also the system \eqref{E:core
  equation'} possesses an exponential dichotomy. Note that from its
definition, it follows immediately that $B(s;\lambda,0)\in
\mathcal{L}(X)$.

As $\tilde\theta$ does not depend on $\varphi$, the
system~\eqref{E:approximated core equation'} decouples in Fourier
space, just as the limiting system \eqref{E:-infty-system}. Using the same notation in the
Fourier spaces as before, we get for $k\in\mathbb{Z}$
\begin{equation}
  \label{E:Fourier_unperturbed}
  \wU_k'(s) = \left[\wA_-(k) + B(s;\lambda,0)\right]\,\wU_k(s).
\end{equation}
As $\wA_-(k)=M_kD_kM_k^{-1}$, we rescale both
$\wU_k$ and $s$ to get estimates which are uniform in~$k$. For $k\neq
0$, define
$\tau = |k|(s-s_1)$ and $V_k(\tau) = M^{-1}_k \wU_k(\tau/|k|+s_1)$.
Then~\eqref{E:Fourier_unperturbed} becomes
\begin{equation}\label{E:Fourier_unperturbed_scaled}
  \frac{d}{d\tau} V_k = \left[D_1 + \frac{1}{|k|}
    M_k^{-1}B(\tau/|k|;\lambda,0)M_k  \right] \, V_k  .
\end{equation}
A short calculation shows that
\begin{equation*}
   |M_k^{-1}B(\tau/|k|;\lambda,0)M_k|\leq \frac {e^{2\tau/|k|}}2 \sup_{s\leq s_1} \{1,
|\lambda-\tilde\theta(s)|/k^2\}.
\end{equation*}
Hence there exists a constant $C$
such that for all $k\neq 0$ , $|M_k^{-1}B(\tau/|k|;\lambda,0)M_k|\leq
2Ce^{2\tau/|k|}$ and
\begin{equation}\label{E:estimate_B}
\int^0_{-\infty}\frac{1}{|k|}
    |M_k^{-1}B(\tau/|k|;\lambda,0)M_k|\, d\tau \leq
\int^0_{-\infty}\frac{2Ce^{2\tau/|k|}}{|k|}\, d\tau = C.
\end{equation}
{By the proof of the
roughness theorem for ordinary dichotomies (see \cite{wC78} for details),
the system~\eqref{E:Fourier_unperturbed_scaled} has an
exponential dichotomy which we denote by $\Psi_k^{u/s}(\tau,\sigma)$,
with constants $K$, $\kappa^u=1$, $\kappa^s=-1$.
We choose the dichotomy in such a way that
$\Ran\Psi_-^{s}(s_1,s_1;\lambda,0)\subset \span\{e_1,e_2\}$, where $e_j$,
$j=1,\dots, 4$ are the standard basis vectors of $\mathbb C^4$ (again see \cite{wC78}).}
This implies that the stable and unstable solutions satisfy
\begin{equation*}
\begin{aligned}
  \left|\Psi^s_k(\sigma,\tau) V_k \right| & \leq K\, e^{-(\sigma-\tau)}\,
  |V_k|, &\quad \tau\leq\sigma\leq0;\\
  \left|\Psi^u_k(\sigma,\tau) V_k \right| & \leq K\, e^{(\sigma-\tau)}\,
  |V_k|, &\quad \sigma\leq\tau\leq0 .
\end{aligned}
\end{equation*}
The norm in $X$ induces a norm on the Fourier space~$X_k$ with
\begin{equation*}
   \|\wU_k\|_{X_k}^2: = \|\wU_k {e^{ik\cdot}}\|_{X}^2 = (k^2+1)^2([\wU_k]_1)^2+
(k^2+1)([\wU_k]_2)^2+(k^2+1)([\wU_k]_3)^2+([\wU_k]_4)^2.
\end{equation*}
As seen in the
proof of Lemma~\ref{L:closed}, $M_k$ is a linear homeomorphism between
$\mathbb{C}^4$ and $X_k$. {Thus if we denote the exponential dichotomy of the unscaled
system~\eqref{E:Fourier_unperturbed} by $\Phi_k^{u/s}$, then
$\Phi_k^{u/s}(s,t)=M_k\Psi_k^{u/s}(|k|(s-s_1),|k|(t-s_1)) M_k^{-1}$, and they
satisfy for $\wU_k\in X_k$}
\begin{equation}\label{E:unperturbed_dichotomy_k_-}
\begin{aligned}
  \|\Phi^s_k(s,t) \wU_k \|_{X_k} &\leq K\, e^{-|k|(s-t)}\,
  \|\wU_k\|_{X_k}{ \leq K\, e^{-(s-t)}\|\wU_k \|_{X_k}}, &\quad t\leq s\leq s_1;\\
  \|\Phi^u_k(s,t) \wU_k \|_{X_k} &\leq K\, e^{|k|(s-t)}\,
  \|\wU_k\|_{X_k}{\leq K\, e^{(s-t)}\|\wU_k\|_{X_k}}, &\quad s\leq t\leq s_1 .
\end{aligned}
\end{equation}
for some constant $K$, which is independent of $k$.

For the central space, corresponding to $k=0$, the scaling $V_0(s)=e^{\pm\epsilon s} \wU_0(s)$
and the integrability of $B(s;\lambda,0)$ shows that, for any $\epsilon>0$
\begin{equation}\label{E:unperturbed_dichotomy_0_-}
\begin{aligned}
  |{\Phi}_0(s,t) \wU_0| &\leq K e^{\epsilon(s-t)} |\wU_0|,
  &\quad t\leq s \leq s_1;\\
  |{\Phi}_0(s,t) \wU_0| &\leq K e^{-\epsilon(s-t)} |\wU_0|,
  &\quad s\leq t \leq s_1.
\end{aligned}
\end{equation}

Thus for the full solutions, we can define the stable and
center--unstable solutions
\begin{equation*}
\begin{aligned}
  \Phi_-^s(s,t;\lambda,0) U &=  \sum_{k\in\mathbb{Z}\backslash\{0\}}\Phi^s_k(s,t)\wU_ke^{ik\cdot},
  & s\leq t\leq s_1,\\
  \Phi_-^{cu}(s,t;\lambda,0) U &= \Phi_0(s,t) \wU_0 +
  \sum_{k\in\mathbb{Z}\backslash\{0\}}\Phi^u_k(s,t)\wU_ke^{ik\cdot},
  &t\leq s\leq s_1,
\end{aligned}
\end{equation*}
and the unstable and center--stable solutions
\begin{equation*}
\begin{aligned}
  \Phi_-^{cs}(s,t;\lambda,0) U &= \Phi_0(s,t) \wU_0 +
  \sum_{k\in\mathbb{Z}\backslash\{0\}}\Phi^s_k(s,t)\wU_ke^{ik\cdot},
  & s\leq t\leq s_1,\\
  \Phi_-^{u}(s,t;\lambda,0) U &=
  \sum_{k\in\mathbb{Z}\backslash\{0\}}\Phi^u_k(s,t)\wU_ke^{ik\cdot},
  &t\leq s\leq s_1.
\end{aligned}
\end{equation*}

These solutions are related to dichotomies for~\eqref{E:approximated
  core equation'} in $X$ on $J_-$.

\begin{Lem}\label{L:dich-smooth}
   {Let $-1=\kappa^s<\kappa^{cu}<0<\kappa^{cs}<\kappa^u=1$ and $\lambda\in\mathbb{R}$.
   Then the system \eqref{E:approximated core equation'} has an exponential
   dichotomy in $X$ on $J_-$ with
   constant $K$ and rates $\kappa^{cu}$ and $\kappa^s$, and another with
   constant $K$ and rates $\kappa^u$ and $\kappa^{cs}$. The dichotomies can
   be chosen such that $\Ran \Phi_-^s(s_1,s_1;\lambda,0) = P^s$ and
   $\Ran \Phi_-^{cs}(s_1,s_1;\lambda,0)= P^{cs}$.
   Moreover, for any $t\in(-\infty,s_1]$ and $U_0\in X$, the solutions
   $\Phi_-^{cu}(\cdot,t;\lambda,0)U_0$ and
   $\Phi_-^{s}(\cdot,t;\lambda,0)U_0$ belong to $C^\infty((-\infty,t);X)$
   and $C^\infty((t,s_1);X)$, respectively. Similarly, the solutions
   $\Phi_-^{u}(\cdot,t;\lambda,0)U_0$ and
   $\Phi_-^{cs}(\cdot,t;\lambda,0)U_0$ belong to
   $C^\infty((-\infty,t);X)$ and $C^\infty((t,s_1);X)$, respectively.
   All solutions also depend smoothly on the parameter $\lambda$.}
\end{Lem}
\begin{proof}
  The scaling $e^{\pm \eta s}U$ for
  $0<\eta<1$ and the dichotomy estimates
  in~\eqref{E:unperturbed_dichotomy_k_-}
  and~\eqref{E:unperturbed_dichotomy_0_-} immediately prove the first
  part of the Lemma. {The dichotomies satisfy $\Ran \Phi_-^s(s_1,s_1;\lambda),0) = P^s$ and
  $\Ran \Phi_-^{cs}(s_1,s_1;\lambda,0)= P^{c} + P^{s}$ since we have chosen the
  $\Psi_k^{u/s}$ above to satisfy
  $\Ran\Psi_-^{s}(s_1,s_1;\lambda,0)\subset \span\{e_1,e_2\}$ (cf. the definition of $D_k$).}

  The smoothness with respect to $s$ follows since $\theta$ is smooth in $s$ and
  smoothness in $\lambda$ can be proved using an implicit function theorem argument.
  First observe that for any $\lambda$, $\tilde\lambda$ close to each
  other, the solutions $\Phi^{cu}$ and $\Phi^{s}$ satisfy the integral
  equations
\begin{equation*}
  \begin{aligned}
    0 = &-\Phi_-^{cu}(s,t;\tilde\lambda,0) +\Phi_-^{cu}(s,t;\lambda,0))
    +(\tilde\lambda-\lambda)\,\left[ \int_{-\infty}^s
      \Phi_-^{s}(s,\tau;\lambda,0) e^{2\tau} B_0
      \Phi_-^{cu}(\tau,t;\tilde\lambda,0)\, d\tau \right.\\
    &{} - \left.\int_{s}^t
      \Phi_-^{cu}(s,\tau;\lambda,0) e^{2\tau} B_0
      \Phi_-^{cu}(\tau,t;\tilde\lambda,0)\, d\tau
      + \int^{s_1}_t
      \Phi_-^{cu}(s,\tau;\lambda,0) e^{2\tau} B_0
      \Phi_-^{s}(\tau,t;\tilde\lambda,0)\, d\tau\right],\\
   0 = &-\Phi_-^{s}(s,t;\tilde\lambda,0) + \Phi_-^{s}(s,t;\lambda,0))
    -(\tilde\lambda-\lambda)\,\left[ \int_{-\infty}^t
      \Phi_-^{s}(s,\tau;\lambda,0) e^{2\tau} B_0
      \Phi_-^{cu}(\tau,t;\tilde\lambda,0)\, d\tau \right.\\
    &{} - \left.\int_{t}^s
      \Phi_-^{s}(s,\tau;\lambda,0) e^{2\tau} B_0
      \Phi_-^{s}(\tau,t;\tilde\lambda,0)\, d\tau
      + \int_{s}^{s_1}
      \Phi_-^{cu}(s,\tau;\lambda,0) e^{2\tau} B_0
      \Phi_-^{s}(\tau,t;\tilde\lambda,0)\, d\tau\right],\\
  \end{aligned}
\end{equation*}
{for $s\le t\le s_1$ and $t\le s\le s_1$, respectively,}
where $B_0$ is the matrix
\begin{equation}\label{E:def B_0}
B_0:=\left(
      \begin{matrix}
        0 & 0 & 0 & 0 \\
        0 & 0 & 0 & 0 \\
        0 & 0 & 0 & 0 \\
        1 & 0 & 0 & 0
      \end{matrix}
    \right).
\end{equation}
Define the function spaces
\begin{equation*}
  \begin{aligned}
    X^s&:=\{\Phi^s\,;\,\Phi^s(s,t)\in \mathcal{L}(X) \mbox{ is defined
      and continuous for } t\leq s\leq s_1\\
    &\qquad\mbox{ with }
    \|\Phi^s\|_s := \sup_{t\leq s\leq s_1}e^{-\kappa^s(s-t)}\|\Phi^s(s,t)\|_{\mathcal L(X)} \};\\
    X^{cu} &:= \{\Phi^{cu}\,;\,\Phi^{cu}(s,t)\in \mathcal{L}(X)
    \mbox{ is defined and continuous for } s\leq t\leq s_1 \\
    &\qquad\mbox{ with } \|\Phi^{cu}\|_{cu} := \sup_{s\leq t\leq
      s_1}e^{-\kappa^{cu}(s-t)}\|\Phi^{cu}(s,t)\|_{\mathcal L(X)} \} .
  \end{aligned}
\end{equation*}
For $\lambda$ fixed, the integral equations can be written as
$F(\Phi^{cu},\Phi^s;\tilde\lambda)=0$, where $F:X^{cu}\times X^s
\times \mathbb{R} \to X^{cu}\times X^s$.  The estimates of the exponential dichotomies
immediately give that $F$ is indeed a mapping between those spaces,
for example,
\[
\left\|\int_{-\infty}^s \Phi_-^{s}(s,\tau;\lambda,0) e^{2\tau} B_0
      \Phi_-^{cu}(\tau,t;\tilde\lambda,0)\,
      d\tau\right\|_{\mathcal{L}(X)} \leq
\int_{-\infty}^s
K^2e^{\kappa^s(s-\tau)}e^{2\tau}e^{\kappa^{cu}(\tau-t)}d\tau =
\frac{K^2e^{2s_1}\,e^{\kappa^{cu}(s-t)}}{2+\kappa^{cu}-\kappa^s}.
\]
The other integrals can be estimated in a similar way. Since
$D_{(\Phi^{cu},\Phi^s)}F
(\Phi^{cu}(s,t;\lambda,0),\Phi^s(s,t;\lambda,0);\lambda)=I$,
the implicit function theorem can be applied and the smoothness
with respect to $\lambda$ follows immediately.
\end{proof}

\begin{Rem}
  The $\varphi$-independence of $\tilde\theta$ is not essential in
  Lemma~\ref{L:dich-smooth}. The lemma can be proved for smooth
  $\varphi$-dependent functions $\theta$ by using Theorem~1
  of~\cite{PSS97} and verifying the conditions (H1), (H2), (H3) and
  (H5) of that paper.
\end{Rem}

Next, we prove four technical lemmas needed in the proof of the
existence of exponential dichotomies for the full system~\eqref{E:core
  equation'} with $\tilde \rho\in\tilde{\mathcal R}$. We work in exponentially
weighted spaces, and for an unbounded interval $J\subset J_-$ and
$\eta\in\mathbb{R}$ , we let $C_\eta(J;X)$ be the space defined by
\begin{equation*}
   C_{\eta}(J;X) := \{U\in C(J;X);\; \|U\|_{C_\eta} := \sup_{s\in J} e^{\eta s} \|U(s)\|_X <\infty\}.
\end{equation*}
Hence $C_0(J,X)$ is the space of continuous functions with an $X$-norm
that is uniformly bounded in $J$.

\begin{Lem}\label{L:product}
  Let $J\subset J_-$ and pick $u\in C_0(J;H^2(S^1))$ and $\rho\in
  L^2(J;H^{1/2}(S^1))$, then $\rho u\in L^2(J;H^{1/2}(S^1))$.
\end{Lem}
\begin{proof}
   We need to prove that  for $s$ fixed,
   \begin{equation}\label{E:spaces}
     \|\rho(s) u(s)\|_{H^{1/2}(S^1)} \le
     C \|u(s)\|_{H^2(S^1)} \|\rho(s)\|_{H^{1/2}(S^1)}.
   \end{equation}
   Indeed, if this is proved, the claim follows, since
   \begin{equation*}
      \begin{aligned}
        \|\rho u\|_{L^2(J;H^{1/2}(S^1))}^2 &=
        \int_{J} \|\rho(s) u(s)\|_{H^{1/2}(S^1)}^2\, ds \\
        &\le C^2 \int_{J} \|\rho(s)\|_{H^{1/2}(S^1)}^2\|u\|_{H^2(S^1)}^2\, ds \\
        &\le C^2 \sup_{s\in J} \|u(s)\|_{H^2(S^1)}^2
        \int_{J} \|\rho(s)\|_{H^{1/2}(S^1)}^2\, ds \\
        &=C^2 \|u\|_{C_0(J;H^2(S^1))}^2 \|\rho\|_{L^2(J;H^{1/2}(S^1))}^2.
      \end{aligned}
   \end{equation*}
   To prove \eqref{E:spaces}, let $u\in H^2(S^1)$ and $\rho\in
   H^{1/2}(S^1)$ (we suppress the variable $s$ for simplicity of
   notation).  Let $\widehat \rho_k$ and $\widehat u_k$ be the Fourier
   coefficients of $\rho$ and $u$, respectively. We have
   \begin{equation*}
      \begin{aligned}
         \|u\|_{H^2}^2 &= \sum_{k\in \mathbb Z} \widehat u_k^2 (1 + k^2)^2, \\
         \|\rho\|_{H^{1/2}}^2 &= \sum_{k\in \mathbb Z} \widehat \rho_k^2 (1 + k^2)^{1/2}.
      \end{aligned}
   \end{equation*}
   Then $(\widehat{u\rho})_k = \sum_{j\in \mathbb Z} \widehat u_j
   \widehat \rho_{k-j}$, and so $\|u \rho\|_{H^{1/2}}^2 = \sum_{k\in
     \mathbb Z}\left(\sum_{j\in \mathbb Z} \widehat u_j \widehat
     \rho_{k-j}\right)^2(1 + k^2)^{1/2}$.  Let $v$ and $\sigma$ be the
   functions with Fourier coefficients $\widehat u_k(1+k^2)^{1/4}$ and
   $\widehat \rho_k(1+k^2)^{1/4}$, respectively. Note that $v\in
   H^{3/2}(S^1)$ and $\sigma\in L^2(S^1)$.

   Now observe that
   \begin{equation*}
      1+k^2=1+((k-j)+j)^2\leq 2(1+j^2)+2(1+(k-j)^2),
   \end{equation*}
   and hence
   \begin{equation*}
      (1+k^2)^{1/4} \le 2^{1/4} ((1+j^2)^{1/4} + (1 + (k-j)^2)^{1/4})
   \end{equation*}
   for any $j\in\mathbb{Z}$. Thus
   \begin{equation*}
      \begin{aligned}
         \|u\rho\|_{H^{1/2}}^2 &\le {\sqrt{2}} \sum_{k\in \mathbb Z} \biggl(\sum_{j\in \mathbb Z} \widehat u_j \widehat
         \rho_{k-j} (1 + j^2)^{1/4} + \sum_{j\in \mathbb Z} \widehat u_j \widehat \rho_{k-j}(1 +
         (k-j)^2)^{1/4}\biggr)^2 \\
         &\le {2 \sqrt{2}} \sum_{k\in \mathbb Z} \biggl(\biggl( \sum_{j\in \mathbb Z} \widehat u_j \widehat \rho_{k-j} (1 +
         j^2)^{1/4}\biggr)^2 + \biggl(\sum_{j\in \mathbb Z} \widehat u_j \widehat \rho_{k-j} (1 +
         (k-j)^2)^{1/4}\biggr)^2\biggr) \\
         &= {2\sqrt{2}} \bigl(\|v \rho\|_{L^2}^2 + \|u\sigma\|_{L^2}^2\bigr) \\
         &\le {2\sqrt{2}} \bigl(\sup_{\varphi\in S^1}|v(\varphi)|^2 \|\rho\|_{L^2}^2 + \sup_{\varphi\in S^1} |u(\varphi)|^2
         \|\sigma\|_{L^2}^2\bigr) \\
         &\le C^2\bigl(\|v\|_{H^{3/2}}^2 \|\rho\|_{L^2}^2 + \|u\|_{H^2}^2 \|\sigma\|_{L^2}^2\bigr) \\
         &\le C^2\bigl(\|u\|_{H^2}^2 \|\rho\|_{H^{1/2}}^2\bigr)
      \end{aligned}
   \end{equation*}
   for some constant $C>0$.
   This completes the proof.
\end{proof}
For each $\tilde \rho\in\tilde{\mathcal R}$ and  $s\in J_-$, let
\begin{equation*}
    \delta B(s;\tilde \rho):=
    B(s;\lambda,\tilde \rho) - B(s;\lambda,0)
    = -e^{2s}\tilde \rho(s)\, B_0,
\end{equation*}
where $B_0$ has been defined in \eqref{E:def B_0}.
Note that for any $s\in J_-$,
\begin{equation}\label{E:delta-B}
  \begin{aligned}
    \|\delta B(s;\tilde \rho)\|_{\mathcal L(X)} &\le \sup_{
      \begin{smallmatrix}
        u_1\in H^2(S^1) \\
        \|u_1\|_{H^2}=1
      \end{smallmatrix}} \|e^{2s}\tilde \rho(s)\,
    u_1\|_{L^2(S^1)}
    \le C \sup_{
      \begin{smallmatrix}
        u_1\in H^2(S^1) \\
        \|u_1\|_{H^2}=1
      \end{smallmatrix}} \sup_{\varphi\in S^1} |u_1(\varphi)|
    \|e^{2s}\tilde \rho(s)\|_{L^2(S^1)}\\
    &\le C e^{2s} \|\tilde \rho(s)\|_{L^2(S^1)},
  \end{aligned}
\end{equation}
where we use the notation $C$ for the different constants occurring. It follows that
\begin{equation}\label{E:delta-B int}
  \begin{aligned}
    \int_{-\infty}^{s_1} \|\delta B(s;\tilde \rho)\|_{\mathcal
      L(X)}^2\, ds &\le C^2 \int_{-\infty}^{s_1} e^{2s}
    \|\tilde \rho(s)\|_{L^2(S^1)}^2 e^{2s}\, ds
    \le C^2 e^{2 s_1} \int_{-\infty}^{s_1}\| \tilde
    \rho(s)\|_{L^2(S^1)}^2 e^{2s}\, ds \\
    &\le C^2 e^{2 s_1}\|\tilde \rho\|_{\widetilde {\mathcal R}}^2.
      \end{aligned}
   \end{equation}
\begin{Lem}\label{L:constant}
  {For $\eta\in(-1,\kappa^{cu})$, where $\kappa^{cu}$ is as in Lemma
  \ref{L:dich-smooth}}, pick $U_-^{cu}\in C_{\eta}(J_-;X)$,
  and $\tilde \rho\in \widetilde{\mathcal R}$. Let $s\in J_-$. Then
  the integral
  \begin{equation*}
    I:=  \int_{-\infty}^s A_- e^{A_- P^s(s-\tau)}
    P^s \delta B(\tau;\tilde \rho) U_-^{cu}(\tau)\, dt
   \end{equation*}
   belongs to $X$.
\end{Lem}
\begin{proof}
  Let $H(\tau):= e^{(\eta - 1) \tau} \delta B(\tau;\tilde
  \rho)U_-^{cu}(\tau)$. By the definition of $\delta B(\tau,\tilde
  \rho)$,
  \begin{equation*}
    H(\tau) = (0,0,0,-e^{(\eta+1) \tau}\tilde \rho(\tau)\,u(\tau))^T
    =: (0,0,0,h(\tau))^T,
   \end{equation*}
   where $u(\tau)$ is the first component of $U_-^{cu}(\tau)$. Then
   $e^{\eta \cdot} u\in C_0(J_-;H^2(S^1))$ and $e^{\cdot} \tilde
   \rho\in L^2(J_-;H^{1/2}(S^1))$, and so by Lemma~\ref{L:product},
   $h\in L^2(J_-;H^{1/2}(S^1))$. For $k\in \mathbb Z$, let $H_k(\tau)$
   and $h_k(\tau)$ be the Fourier coefficients of $H(\tau)$ and
   $h(\tau)$, respectively. Let $\widetilde P_k^s:= M_k^{-1} P_k^s
   M_k$. To show that $I$ exists in $X$, it suffices to show that
   $\{I_k\}_{k\in \mathbb Z}\in \widehat X$ (see~\eqref{E:X hat}), where
   \begin{equation*}
      \begin{aligned}
         I_k :&= \int_{-\infty}^s e^{(1-\eta)\tau} M_k D_k e^{D_k
           \tilde P_k^s (s-\tau)} M_k^{-1} P^s H_k(\tau)\, d\tau  \\
         &= \frac{1}{2} \int_{-\infty}^s e^{(1-\eta)\tau}
         e^{-|k|(s-\tau)} h_k(\tau)\, d\tau (0,0,1,-|k|)^{T}.
      \end{aligned}
   \end{equation*}
   We therefore need to prove that
   \begin{equation*}
      \left\{ |k| \int_{-\infty}^s e^{(1-\eta)\tau} e^{-|k|(s-\tau)}
        h_k(\tau)\, d\tau\right\}_{k\in \mathbb Z} \in  l^2.
   \end{equation*}
   Using that $\eta<0$ and $h_k\in L^2(J_-)$ (as $h\in
   L^2(J_-;H^{1/2}(S^1))$), we note that
   \begin{equation*}
      \begin{aligned}
         |k| \left|\int_{-\infty}^s e^{(1-\eta)\tau} e^{-|k|(s-\tau)} h_k(\tau)\, d\tau\right| &\le |k|
         e^{-|k|s}\left(\int_{-\infty}^s e^{2(1-\eta + |k|)\tau}d\tau\right)^{1/2} \left( \int_{-\infty}^s
         h_k(\tau)^2\, d\tau\right)^{1/2} \\
         &= |k| \frac{1}{\sqrt{2(1-\eta+|k|)}} e^{(1-\eta)s} \|h_k\|_{L^2((-\infty,s])} \\
         &\le e^{(1-\eta)s_1} (1 + |k|^2)^{1/4} \|h_k\|_{L^2(J_-)}.
      \end{aligned}
   \end{equation*}
   Since $\{(1 + |k|^2)^{1/4} \|h_k\|_{L^2(J_-)}\}_{k\in \mathbb
     Z}\in l^2$, the proof is complete.
\end{proof}
\begin{Lem}\label{L:B-integral}
  For $-1<\eta<\kappa^{cu}$, where $\kappa^{cu}<0$ is as in Lemma
  \ref{L:dich-smooth}, pick $U_-^{cu}\in C_\eta(J_-;X)$
  and $\tilde \rho\in \widetilde{\mathcal R}$. Then the integrals
   \begin{equation*}
     \int_{-\infty}^s \Phi_-^s(s,\tau;\lambda,0) \delta
     B(\tau;\tilde \rho) U_-^{cu}(\tau)\, d\tau
\mbox{~~and~~}
\int_{-\infty}^s B(s;\lambda,0)\Phi_-^s(s,\tau;\lambda,0) \delta
     B(\tau;\tilde \rho) U_-^{cu}(\tau)\, d\tau
   \end{equation*}
   exist in $X$ for each $s\in J_-$.
\end{Lem}
\begin{proof}
  We use \eqref{E:delta-B} and compute
  \begin{equation*}
    \begin{aligned}
      &\left\| \int_{-\infty}^s \Phi_-^s(s,\tau;\lambda,0)
        \delta B(\tau;\tilde \rho) U_-^{cu}(\tau)\, d\tau \right\|_X \\
      &\qquad \le C \int_{-\infty}^s
      \|\Phi_-^s(s,\tau;\lambda,0)\|_{\mathcal L(X)}
      e^{2\tau} \|\tilde \rho(\tau)\|_{L^2(S^1)} \|U_-^{cu}(\tau)\|_X\, d\tau \\
      &\qquad \le C \int_{-\infty}^s e^{2\tau} e^{\kappa^s(s-\tau)}
      \|\tilde \rho(\tau)\|_{L^2(S^1)}
      \|U_-^{cu}(\tau)\|_X\, d\tau \\
      &\qquad \le C \|U_-^{cu}\|_{C_\eta(J_-;X)} e^{\kappa^s
        s} \left(\int_{-\infty}^s e^{2(1 - \kappa^s -
          \eta)\tau}\,d\tau\right)^{1/2} \left(\int_{-\infty}^s
        e^{2\tau} \|\tilde \rho(\tau)\|_{L^2(S^1)}^2\, d\tau
      \right)^{1/2} \\
      &\qquad \le
      C \|U_-^{cu}\|_{C_\eta(J_-;X)}\frac{1}{\sqrt{2(1-\kappa^s-
          \eta)}} e^{(1 - \eta)s} \|\tilde \rho\|_{\tilde {\mathcal R}}.
      \end{aligned}
   \end{equation*}
   {Using that $B(s;\lambda,0)\in \mathcal L(X)$, it follows that both integrals converge in $X$.}
\end{proof}
\begin{Lem}\label{L:A-integral}
  Let $-1<\eta<\kappa^{cu}$, where $\kappa^{cu}<0$ is as in Lemma
  \ref{L:dich-smooth}. Let $U_-^{cu}\in C_\eta(J_-;X)$
  and  $\tilde \rho\in \widetilde{\mathcal R}$.  For every $s\in
  J_-$, the integral
  \begin{equation}\label{E:A_-estimate}
    \int_{-\infty}^s A_- \Phi_-^s(s,\tau;\lambda,0)
    \delta B(\tau;\tilde\rho) U_-^{cu}(\tau)\, d\tau
   \end{equation}
   exists in $X$.
\end{Lem}
\begin{proof}
  By (3.1) of \cite{PSS97}, for $\tau \le s \le s_1$,
  \begin{equation*}
    \begin{aligned}
      \Phi_-^s(s,\tau;\lambda,0) &= e^{A_- P^s(s-\tau)}
      P^s - \int_{-\infty}^\tau e^{A_- P^s(s-\xi)} P^s
      B(\xi;\lambda,0) \Phi_-^{cu}(\xi,\tau;\lambda,0)\, d\xi   \\
      &\qquad + \int_\tau^s e^{A_- P^s (s-\xi)} P^s
      B(\xi;\lambda,0) \Phi_-^s(\xi,\tau;\lambda,0)\, d\xi \\
      &\qquad- \int_s^{s_1} e^{A_- P^{cu}(s-\xi)} P^{cu}
      B(\xi;\lambda,0) \Phi_-^s(\xi,\tau;\lambda,0)\, d\xi.
      \end{aligned}
   \end{equation*}
   Note that we used here that $\Ran \Phi_-^s(s_1,s_1;\lambda,0)$ has been
   chosen so that it coincides with $\Ran P^s$.

   Substituting this into \eqref{E:A_-estimate}, we have four
   integrals to estimate, the first of which was dealt with in
   Lemma~\ref{L:constant}. The other three integrals are
   \begin{equation*}
     \begin{aligned}
       I_1 &:= \int_{-\infty}^s A_- \int_{-\infty}^\tau
       e^{A_-  P^s(s-\xi)} P^s B(\xi;\lambda,0)
       \Phi_-^{cu}(\xi,\tau;\lambda,0)\, d\xi \delta B(\tau;\tilde
       \rho) U_-^{cu}(\tau)\,
       d\tau, \\
       I_2 &:= \int_{-\infty}^s A_- \int_\tau^s e^{A_- P^s (s-\xi)}
       P^s B(\xi;\lambda,0) \Phi_-^s(\xi,\tau;\lambda,0)\, d\xi
       \delta B(\tau;\tilde \rho) U_-^{cu}(\tau)\, d\tau, \\
       I_3 &:= \int_{-\infty}^s A_- \int_s^{s_1} e^{A_- P^{cu}(s-\xi)}
       P^{cu} B(\xi;\lambda,0) \Phi_-^s(\xi,\tau;\lambda,0)\, d\xi
       \delta B(\tau;\tilde \rho) U_-^{cu}(\tau)\, d\tau.
      \end{aligned}
   \end{equation*}
   We carry out the calculations for $I_1$, since the others are
   similar.  Let $(\phi_{jl}(\xi,\tau))$, $j,l=1,\dots,4$, be the entries of the matrix
   corresponding to $\Phi_-^{cu}(\xi,\tau;\lambda,0)$, and as in the
   proof of Lemma~\ref{L:constant}, let $h(\tau) =- e^{(\eta+1)
     \tau}\tilde \rho(\tau)\, u(\tau)$. Recall that $h\in
   L^2(J_-;H^{1/2}(S^1))$.  A short calculation shows that
   \begin{equation*}
      B(\xi;\lambda,0) \Phi_-^{cu}(\xi,\tau;\lambda,0)
      \delta B(\tau;\tilde\rho) U_-^{cu}(\tau)      e^{2 \xi + (1-\eta)\tau}\left(
      \begin{matrix} 0 \\ \phi_{34}(\xi,\tau)h(\tau) \\ 0 \\
        (\lambda - \tilde \theta(\tau))\phi_{14}(\xi,\tau) h(\tau)
      \end{matrix}
      \right).
   \end{equation*}
   Note that $\phi_{34}(\xi,\tau)$ and $\phi_{14}(\xi,\tau)$ map
   $L^2(S^1)$ {boundedly} into $H^1(S^1)$ and $H^2(S^1)$, respectively, and that
   by Lemma~\ref{L:dich-smooth} for $\xi\le\tau\le s_1$
   \begin{equation*}
      \begin{aligned}
         \|\phi_{34}(\xi,\tau)\|_{\mathcal L(L^2;H^1)}&\le
         K e^{\kappa^{cu}(\xi-\tau)}, \\
         \|\phi_{14}(\xi,\tau)\|_{\mathcal L(L^2;H^2)}&\le
         K e^{\kappa^{cu}(\xi-\tau)}.
      \end{aligned}
   \end{equation*}

   Introducing the notation $f(\xi,\tau):= e^{-\kappa^{cu}(\xi-\tau)}
   \phi_{34}(\xi,\tau) h(\tau)$, and $g(\xi,\tau):   e^{-\kappa^{cu}(\xi-\tau)} (\lambda - \tilde
   \theta(\tau))\phi_{14}(\xi,\tau) h(\tau)$, we note that
   $\max(\|f(\xi,\tau)\|_{H^1},\|g(\xi,\tau)\|_{{H^2}})\le K
   \|h(\tau)\|_{L^2}$, for $\xi<\tau<s_1$. The Fourier coefficients of
   $f(\xi,\tau)$ and $g(\xi,\tau)$ are denoted by $\widehat
   f_k(\xi,\tau)$ and $\widehat g_k(\xi,\tau)$, respectively.

   To prove that $I_1\in X$, it suffices to prove that $\{J_k\}_{k\in
     \mathbb Z}\in \widehat X$, where
   \begin{equation}\label{E:J_k}
      \begin{aligned}
         J_k :&= \int_{-\infty}^s  M_k D_k \int_{-\infty}^\tau e^{2 \xi + (1-\eta)\tau} e^{\kappa^{cu}(\xi-\tau)} e^{D_k \widetilde P_k^s(s-\xi)} M_k^{-1} P^s \left(
         \begin{matrix}
            0 \\ \widehat f_k(\xi,\tau) \\ 0 \\ \widehat g_k(\xi,\tau)
         \end{matrix}
         \right)\,  d\xi\, d\tau \\
         &= \frac{1}{2} \int_{-\infty}^s \int_{-\infty}^\tau e^{2 \xi + (1-\eta)\tau} e^{\kappa^{cu}(\xi-\tau)} e^{-|k|(s-\xi)}
         \left(
         \begin{matrix}
            \widehat f_k(\xi,\tau) \\ -|k| \widehat f_k(\xi,\tau) \\ \widehat g_k(\xi,\tau) \\ -|k| \widehat
            g_k(\xi,\tau)
         \end{matrix}
         \right)\, d\xi\, d\tau,
      \end{aligned}
   \end{equation}
   {where $\tilde P_k^s= M_k^{-1} P_k^s M_k$ as before.}
   The first component of $J_k$ can be written
   \begin{equation*}
      \begin{aligned}
         &\frac{1}{2} e^{-|k|s} \int_{-\infty}^s e^{(1-\kappa^{cu}-\eta)\tau} \int_{-\infty}^\tau
         e^{(1+\kappa^{cu}+|k|)\xi} e^{\xi}\widehat f_k(\xi,\tau)\, d\xi\, d\tau \\
         &\qquad\qquad \le \frac{1}{2} e^{-|k|s} \int_{-\infty}^s
         \frac{1}{\sqrt{2(1+\kappa^{cu}+|k|)}} e^{(2+|k|-\eta)\tau} \left(\int_{-\infty}^\tau e^{2\xi} \widehat f_k(\xi,\tau)^2\,
         d\xi\right)^{1/2}\, d\tau \\
         &\qquad\qquad\le \frac{1}{4}
         \frac{e^{(2-\eta)s}}{\sqrt{(1+\kappa^{cu}+|k|)(2+|k|-\eta)}}\left(\int_{-\infty}^s\int_{-\infty}^\tau e^{2\xi} \widehat
         f_k(\xi,\tau)^2\, d\xi\, d\tau\right)^{1/2}.
      \end{aligned}
   \end{equation*}
   The square of the $l_2^2$ norm of the first component of
   $\{J_k\}_{k\in\mathbb Z}$ can then be estimated  by
   \begin{equation*}
      \begin{aligned}
         &\frac{e^{2(2-\eta)s}}{16} \sum_{k\in \mathbb Z} \frac{1+k^2}{(1+\kappa^{cu}+|k|)(2+|k|-\eta)}
         \int_{-\infty}^s\int_{-\infty}^\tau (1+k^2)e^{2\xi}\widehat f_k(\xi,\tau)^2\, d\xi\, d\tau \\
         &\qquad\qquad \le \frac{e^{2(2-\eta)s}}{16(1+\kappa^{cu})} \sum_{k\in\mathbb Z} \int_{-\infty}^s \int_{-\infty}^\tau
         e^{2\xi}\widehat h_k(\tau)^2\, d\xi\, d\tau \\
         &\qquad \qquad = \frac{e^{2(2-\eta)s}}{32(1+\kappa^{cu})}\sum_{k\in\mathbb Z} \int_{-\infty}^s
         e^{2\tau}\widehat h_k(\tau)^2\, d\tau \\
         &\qquad \qquad \le \frac{e^{2(3-\eta)s_1}}{32(1+\kappa^{cu})}
         \|h\|_{L^2(J_-;L^2(S^1))}^2.
      \end{aligned}
   \end{equation*}
   The other three components of $\{J_k\}_{k\in\mathbb Z}$ are
   estimated in a completely similar way, and by adding these
   estimates we see that $I_1\in X$.
\end{proof}

We are now ready to prove the existence of exponential dichotomies for
the full system.
\begin{Thm}\label{T:dich-nonsmooth}
  Let $-1<\kappa^{s}<\kappa^{cu}<0$ and $0<\kappa^{cs}<\kappa^u<1$.
  Then there exists a neighbourhood $\mathcal{U}$ of~$0$ in~$\tilde{\mathcal
    R}$ such that for any $\tilde\rho\in\mathcal{U}$ and any
  $\lambda\in\mathbb{R}$, the system~\eqref{E:perturbed
    system}
  has an exponential dichotomy on $J_-$ with constants $K$ and rates
  $\kappa^{cu}$, $\kappa^s$, and another with constants $K$ and rates
  $\kappa^u$, $\kappa^{cs}$. Moreover, the projections and
  evolution operators depend smoothly on $\lambda\in\mathbb{R}$ and
  $\tilde \rho\in \mathcal{U}$.
  The dichotomies are denoted by $\Phi_-^{s}(s,t;\lambda,\tilde
  \rho)$, $\Phi_-^{cu}(s,t;\lambda,\tilde \rho)$, and
  $\Phi_-^{cs}(s,t;\lambda,\tilde \rho)$, $\Phi_-^u(s,t;\lambda,\tilde
  \rho)$, respectively.  The associated projections will be denoted by
  $P^s_-(s;\lambda,\tilde\rho)(:=\Phi_-^{s}(s,s;\lambda,\tilde \rho))$,
  $P^{cu}_-(s;\lambda,\tilde\rho)$, $P^{cs}_-(s;\lambda,\tilde\rho)$,
  and $P^{u}_-(s;\lambda,\tilde\rho)$, respectively.
\end{Thm}
\begin{proof}
  We will show that there exists a neighbourhood of~0 in $\tilde {\mathcal
    R}$ such that if $\tilde \rho$ belongs to this neighbourhood then there
    exist exponential dichotomies for the system \eqref{E:core
    equation'} with this $\tilde \rho$.
  Let $U_0\in X$ and $t\in J_-$ be fixed but arbitrary.  We will use
  the implicit function theorem to solve the system of integral
  equations for the pair of functions $(U_-^{cu},U_-^s)$ as functions
  of the parameters $\lambda\in \mathbb{R}$ and $\tilde\rho
  \in\tilde{\mathcal{R}}$ near~$0$
  \begin{equation}\label{E:dichotomyequation}
    \begin{aligned}
      0 &= \Phi_-^{cu}(s,t;\lambda,0)U_0 - U_-^{cu}(s) +
      \int_{-\infty}^s \Phi_-^{s}(s,\tau;\lambda,0)
      \delta B(\tau;\tilde \rho) U_-^{cu}(\tau)\, d\tau \\
      & \qquad - \int_s^t \Phi_-^{cu}(s,\tau;\lambda,0)
      \delta B(\tau;\tilde \rho) U_-^{cu}(\tau)\, d\tau \\
      & \qquad + \int_t^{s_1} \Phi_-^{cu}(s,\tau;\lambda,0) \delta
      B(\tau;\tilde \rho) U_-^{s}(\tau)\, d\tau,
      \qquad \text{for }s \le t \le s_1, \\
      0 &= \Phi_-^{s}(s,t;\lambda,0)U_0 - U_-^{s}(s) -
      \int_{-\infty}^t \Phi_-^{s}(s,\tau;\lambda,0)
      \delta B(\tau;\tilde \rho) U_-^{cu}(\tau)\, d\tau \\
      & \qquad +\int_t^s \Phi_-^{s}(s,\tau;\lambda,0) \delta
      B(\tau;\tilde \rho) U_-^{s}(\tau)\, d\tau \\
      & \qquad - \int_s^{s_1} \Phi_-^{cu}(s,\tau;\lambda,0) \delta
      B(\tau;\tilde \rho,0) U_-^{s}(\tau)\, d\tau, \qquad
      \text{for }t\le s\le s_1.
      \end{aligned}
   \end{equation}
   By Lemma~\ref{L:dich-smooth}, the dichotomies
   $\Phi_-^{cu}(s,t;\lambda,0)U_0$ and
   $\Phi_-^{s}(s,t;\lambda,0)U_0$ exist and have constants
   $K$, $\tilde\kappa^{s}=-1$ and $\tilde \kappa^{cu}\in(-1,0)$.

   Let $\eta\in (\kappa^s,\kappa^{cu})$ and rewrite equation
   \eqref{E:dichotomyequation} as
   $F(U_-^{cu},U_-^{s};\lambda,\tilde\rho)=0$, where
   $F:C_\eta((-\infty,t];X)\times C_\eta([t,s_1];X)\times \mathbb
   R\times \tilde{\mathcal R} \to C_\eta((-\infty,t];X)\times
   C_\eta([t,s_1];X)$ is the right hand side of
   \eqref{E:dichotomyequation}.

   We first verify that $F$ is indeed a map between the above
   spaces. We do the estimates for the first integral in the first
   equation of \eqref{E:dichotomyequation}. The other estimates are
   similar.  Lemma~\ref{L:dich-smooth} gives that for any $s\in
   (-\infty,t]$ and $U^{cu}\in C_\eta((-\infty,t];X)$:
   \begin{equation}\label{E:dichotomy estimate}
     \begin{aligned}
       e^{\eta s} \bigg\| \int_{-\infty}^s
       \Phi_-^{s}&(s,\tau;\lambda,0) \delta B(\tau;\tilde\rho)
       U_-^{cu}(\tau)\, d\tau \bigg\|_X \\
       &\le K \sup_{\tau\in(-\infty,s]} \left(e^{\eta \tau}
         \|U_-^{cu}(\tau)\|_X\right) \int_{-\infty}^s e^{(\kappa^{s}-\eta) (s-\tau)}
       \|\delta B(\tau;\tilde \rho)\|_{\mathcal L(X)}\, d\tau \\
       &\le K \|U_-^{cu}\|_{C_\eta((-\infty,t],X)} \int_{-\infty}^s
       \left(e^{2(\kappa^s - \eta)(s-\tau)} +
         \|\delta B(\tau;\tilde\rho)\|_{\mathcal L(X)}^2\right)\, d\tau  \\
       &\le K \|U_-^{cu}\|_{C_\eta((-\infty,t],X)}\left(
         \frac{1}{2(\eta - \kappa^s)} + e^{2s_1}\|\tilde
         \rho\|_{\tilde {\mathcal R}}^2 \right),
      \end{aligned}
   \end{equation}
   where we have used \eqref{E:delta-B int}. After taking the supremum
   over all $s\in (-\infty,s_1]$ we see that the function defined by
   the first integral in~\eqref{E:dichotomyequation} belongs to
   $C_\eta(J_-,X)$. Using similar estimates for the other integrals,
   we can conclude that $F$ is indeed a map between the spaces as
   stated.

   That $F$ is smooth with respect to $\lambda$ and $\tilde \rho$
   follows since the evolution operators
   $\Phi_-^{cu}(\cdot,t;\lambda,0)U_0$ and
   $\Phi_-^{s}(\cdot,t;\lambda,0)U_0$ are smooth in $\lambda$ by
   Lemma~\ref{L:dich-smooth} (using that the $H^1$ norm is weaker than
   the $C^1$ norm on bounded intervals), and since $\delta B$ depends
   smoothly on $\tilde \rho$ (indeed, $\delta B$ is a bounded linear
   mapping with respect to~$\tilde \rho$). Note that
   \begin{equation*}
     F(\Phi_-^{cu}(\cdot,t;\lambda,0)U_0,\Phi_-^{s}(\cdot,t;\lambda,0)U_0;
     \lambda,0) = 0.
   \end{equation*}
   The Fr\'echet derivative of $F$ with respect to its two
   first variables evaluated at 
   $(\Phi_-^{cu}(\cdot,t;\lambda,0)U_0,
   \Phi_-^{s}(\cdot,t;\lambda,0) U_0;\lambda,0)$ is $-I$ on
   $C_\eta((-\infty,t];X)\times C_\eta([t,s_1];X)$. In particular,
   this derivative is a linear homeomorphism on this space, and so the
   implicit function theorem is applicable, and we obtain solutions
   $\Phi_-^{cu}(\cdot,t;\lambda,\tilde \rho)U_0:= U_-^{cu}$ and
   $\Phi_-^{s}(\cdot,t;\lambda,\tilde \rho) U_0:= U_-^{s}$ of the
   integral equation \eqref{E:dichotomyequation}, which exist in a
   neighbourhood of $(\lambda_0,0)$ in $\mathbb R\times
   \tilde {\mathcal R}$. Smoothness of these solutions with respect to
   parameters also follows from a corollary of the implicit function
   theorem (see e.g. \cite[p.115]{mB77}).

   Next, we need to verify that $\Phi_-^{cu}(\cdot,t;\lambda_0,\tilde
   \rho)U_0$ and $\Phi_-^{s}(\cdot,t;\lambda_0,\tilde \rho) U_0$ are
   weak solutions of \eqref{E:core equation'}, and that they satisfy
   the conditions of Definition~\ref{D:dichotomies}.
   We first check that $\Phi_-^{cu}(\cdot,t;\lambda,\tilde \rho)U_0$
   is a weak solution on the interval $(-\infty,t]$. By
   Lemma~\ref{L:dich-smooth}, $\Phi_-^{cu}(\cdot,t;\lambda,0)U_0$ is a
   $C^\infty$ solution of
   \begin{equation*}
      U' = (A_{-} + B(s;\lambda,0)) U
   \end{equation*}
   on $(-\infty,t]$, and hence it is also a weak solution of this
   equation. Next we deal with the integral terms. For the first
   integral we use the abbreviation
   \begin{equation*}
     g(s):= \int_{-\infty}^s f(s,\tau)\, d\tau, \mbox{~~with~~}
     f(s,\tau) = \Phi_-^{s}(s,\tau;\lambda,0)
     \delta B(\tau;\tilde \rho) \Phi_-^{cu}(\tau,t;\lambda_0,\tilde
     \rho)U_0.
   \end{equation*}
   Thus $f$ is $C^{\infty}$ in the first variable and $L^1$ in the second.
   From its definition, it follows immediately that $g$ is
   continuous. We will see that $g$ is weakly differentiable and that
   \begin{equation}\label{E:weak derivative}
     g'(s) = f(s,s) + \int_{-\infty}^s \frac{\partial f}{\partial
       s}(s,\tau)\, d\tau.
   \end{equation}
   In order to prove this, we need to check that the integral on the
   right hand side of \eqref{E:weak derivative} exists, and that the
   equality \eqref{E:weak derivative} holds.
   The integral in the right hand side of \eqref{E:weak derivative} is
   \begin{equation*}
     \int_{-\infty}^s \frac{\partial f}{\partial s}(s,\tau)\, d\tau
     = \int_{-\infty}^s (A_- + B(s;\lambda,0))
     \Phi_-^s(s,\tau;\lambda,0) \delta B(\tau;\tilde \rho)
     \Phi_-^{cu}(\tau,t;\lambda,\tilde\rho)U_0\, d\tau,
   \end{equation*}
   and it exists in $X$ by Lemma~\ref{L:B-integral} and
   Lemma~\ref{L:A-integral}.

   Next, we calculate the distributional derivative of $g$ and let
   $V\in C_0^{\infty}((-\infty,t];X)$ be a test function. Then by
   Fubini's Theorem and integration by parts
   \begin{equation*}
      \begin{aligned}
         \int_{-\infty}^t g'(s) V(s)\, ds &= -\int_{-\infty}^t g(s) V'(s)\, ds \\
         &=-\int_{-\infty}^t \int_{-\infty}^s f(s,\tau)\, d\tau V'(s)\, ds \\
         &=-\int_{-\infty}^t \int_\tau^t f(s,\tau) V'(s)\, ds\, d\tau \\
         &= \int_{-\infty}^t\left(f(\tau,\tau) V(\tau) + \int_\tau^t
           \frac{\partial f}{\partial s}(s,\tau) V(s)\, ds\right)\,
         d\tau \\
         &= \int_{-\infty}^t \left(f(s,s) + \int_{-\infty}^s
           \frac{\partial f}{\partial s}(s,\tau)\, d\tau\right) V(s)\,
         ds,
      \end{aligned}
   \end{equation*}
   and we see that the weak derivative of $g$ is indeed given by
   \eqref{E:weak derivative}.  Hence
   \begin{equation}
      \begin{aligned}\label{E:first integral}
         \frac{d}{ds}\int_{-\infty}^s
         &\Phi_-^{s}(s,\tau;\lambda,0)\delta B(\tau;\tilde \rho)
         \Phi_-^{cu}(\tau,t;\lambda,\tilde \rho)U_0\, d\tau \\
         &=(I - P_-^{cu}(s;\lambda,0))\delta B(s;\tilde \rho)
         \Phi_-^{cu}(s,t;\lambda,\tilde \rho)U_0 \\
         &+ \int_{-\infty}^s (A_{-} + B(s;\lambda,0))
         \Phi_-^{s}(s,\tau;\lambda,0)\delta B(\tau;\tilde \rho)
         \Phi_-^{cu}(\tau,t;\lambda,\tilde \rho)U_0\, d\tau.
      \end{aligned}
   \end{equation}
   We have already noticed that $g(s)=\int_{-\infty}^s
   \Phi_-^{s}(s,\tau;\lambda,0)\delta B(\tau;\rho)
   \Phi_-^{cu}(\tau,t;\lambda,\tilde \rho)U_0\, d\tau$ is continuous,
   and so it belongs to $L^2_{loc}((-\infty,t];X)$. The right hand
   side of \eqref{E:first integral} also belongs to
   $L^2_{loc}((-\infty,t];X)$ since the first term belongs to
   $L^2_{loc}((-\infty,t];X)$ and the second term is continuous on
   $(-\infty,t]$. This shows that
   \[
   \int_{-\infty}^s \Phi_-^{s}(s,\tau;\lambda,0)\delta B(\tau;\tilde
   \rho) \Phi_-^{cu}(\tau,t;\lambda,\tilde \rho)U_0\, d\tau
   \]
   belongs to $H^1_{loc}((-\infty,t];X)$.

   Similar calculations for the other integral terms of the first
   equation of \eqref{E:dichotomyequation} show that these are also
   weakly differentiable on $(-\infty,t]$ and belong to
   $H^1_{loc}((-\infty,t],X)$. After adding the terms up, we conclude that
   \begin{equation*}
      \frac{d}{ds}\Phi_-^{cu}(s,t;\lambda,\tilde \rho) U_0 =
      (A_{-} + B(s;\lambda,\tilde \rho))
      \Phi_-^{cu}(s,t;\lambda,\tilde \rho)U_0,
   \end{equation*}
   i.e.\ $\Phi_-^{cu}(\cdot,t;\lambda,\tilde \rho)U_0$ is a weak
   solution of \eqref{E:core equation'}.

   Similar calculations for the terms of the second equation of
   \eqref{E:dichotomyequation} show that
   $\Phi^s_-(\cdot,t;\lambda,\tilde \rho)U_0$ is a weak solution of
   \eqref{E:core equation'} on the interval $[t,s_1]$.

 Finally we check that the conditions of
 Definition~\ref{D:dichotomies} are satisfied.
   A similar computation as in~\eqref{E:dichotomy estimate} also shows
   that the estimates in (i) and (ii) of
   Definition~\ref{D:dichotomies} are satisfied for
   $\Phi_-^{cu}(s,t;\lambda,\rho)$ and $\Phi_-^s(s,t;\lambda,\tilde
   \rho)$ for any $\kappa^{cu}$ and $\kappa^{s}$ such that
   $-1=\tilde \kappa^{s}<\kappa^{s}<\kappa^{cu}< \tilde\kappa^{cu}<0$.
   {Since $\tilde \kappa^{cu}$ can be taken arbitrarily close to $1$,
   the same is true also for $\kappa^{cu}$.}
    Note that (iii)
   of Definition~\ref{D:dichotomies} is satisfied with
   $P_-^{cu}(s;\lambda,\tilde \rho) := \Phi^{cu}(s,s;\lambda,\tilde
   \rho)$ and $P_-^s(s;\lambda,\tilde \rho):=   \Phi_-^s(s,s;\lambda,\tilde \rho)$.
\end{proof}

To finish this section, we derive some more details about the
solutions of~\eqref{E:perturbed system} in the case when
$\lambda=\lambda_0$ and $\tilde \rho=0$. We are particularly
interested in the solutions on $J_-$, and we study the exact
growth/decay rate of solutions as $s\to -\infty$.  As we have seen
before, the space $X$ decouples into a direct sum of $4$-dimensional
pairwise orthogonal Fourier subspaces $X_k$, and that since $\theta$
is radially symmetric, the subspaces~$X_k$ are invariant both under the flow
of~\eqref{E:perturbed system} with $\tilde\rho=0$ and under the flow
of the asymptotic system~\eqref{E:-infty-system}.
{\begin{Lem}\label{L:centre}
  Let $e_j$, $j=1,\dots, 4$ be the standard basis of ${\mathbb C^4}$ and
  consider the four-dimensional invariant central space corresponding
  to $k=0$ of the unperturbed equation obtained when $\tilde\rho=0$
  and $\lambda=\lambda_0$ in~\eqref{E:perturbed system}. Then there exist
  two unique solutions $U_{0,j}(s)$ with $j=1,3$ such that
   \begin{equation*}
      \lim_{s\to -\infty} U_{0,j}(s) = e_j, \qquad j=1,3.
   \end{equation*}
   We may also pick two solutions $U_{0,j}$ with $j=2,4$ which grow
   algebraically as $s\to-\infty$ and satisfy
   \begin{equation*}
      \lim_{s\to-\infty}\frac {1}{s}U_{0,j}(s) = e_{j-1}, \qquad j=2,4.
   \end{equation*}
   The solutions $U_{0,j}$, $j=1,\dots,4$, are
   linearly independent.
\end{Lem}
\begin{proof}
   It is straightforward to check the assertions of the lemma, using \cite[Chapter 3.8]{CL55}.
\end{proof}}
{In Section \ref{S:adjoint} we will specify the solutions $U_{0,2}$ and $U_{0,4}$ using the adjoint system.}
\begin{Lem}\label{L:high Fourier modes}
  For every $k\in \mathbb Z\setminus\{0\}$, there exist solutions
  $U_{k,j}$ of \eqref{E:perturbed system} with $\tilde\rho=0$ and
  $\lambda=\lambda_0$ such that (together with the solutions specified
  in Lemma~\ref{L:centre} for $k=0$) we have
   \begin{equation*}
      \span\{U_{k,j}(s_1);\; k\in \mathbb Z,\, j=1,\dots, 4\} = X,
   \end{equation*}
   and for $s\to-\infty$,
   \begin{equation*}
      \begin{aligned}
         e^{|k| (s-s_1)}\, U_{k,1}(s) &\to
         (-1/k^2,1/|k|,0,0)^T  e^{ik\cdot}, \\
         e^{|k| (s-s_1)} \, U_{k,2}(s) &\to
         (0,0,-1/|k|,1)^T  e^{ik\cdot}, \\
         e^{-|k| (s-s_1)} \, U_{k,3}(s) &\to
         (1/k^2,1/|k|,0,0)^T  e^{ik\cdot}, \\
         e^{-|k| (s-s_1)} \, U_{k,4}(s) &\to
         (0,0,1/|k|,1)^T  e^{ik\cdot}.
      \end{aligned}
   \end{equation*}
\end{Lem}
\begin{proof}
  As seen in the beginning of this section, the
  system~\eqref{E:perturbed system} with $\tilde\rho=0$ and
  $\lambda=\lambda_0$ leaves the subspaces~$X_k$ invariant. The
  estimates on the matrix $B(s;\lambda_0,0)$ in~\eqref{E:estimate_B}
  now show that there are solutions of~\eqref{E:Fourier_unperturbed} (and hence
  of~\eqref{E:perturbed system}) which converge to the solutions of
  the system at infinity, see e.g.~\cite[Chapter~3.8]{CL55}.
  In Section~\ref{S:-infty} we have seen that \eqref{E:-infty-system}
  has two solutions in $X_k$ with decay rate  $e^{|k|s}$ and two
  with growth rate $e^{-|k|s}$ for $s\to-\infty$. A comparison
  with the eigenvectors of $\wA_{-}(k)$ in Section~\ref{S:-infty}, we
  obtain solutions $U_{k,j}$, with $k\in \mathbb Z\setminus\{0\}$ and
  $j=1,\dots,4$ with the desired properties.
\end{proof}

Next, we perturb the solutions $U_{0,j}(s)$ with $j=1,3$ described in
Lemma~\ref{L:centre} to solutions of \eqref{E:perturbed system}
for all sufficiently small potentials $\tilde\rho$.  First, we will
show that the four-dimensional central subspace corresponding to $k=0$
persists in~\eqref{E:perturbed system} as the intersection of the
ranges of $\Phi^{cs}_-({s_1,s_1};\lambda,\tilde\rho)$ and
$\Phi^{cu}_-({s_1,s_1};\lambda,\tilde\rho)$. Note that the difference
between the operators $A(s;\lambda,\tilde\rho)$ and $A(s;\lambda_0,0)$
in~\eqref{E:perturbed system} is
\begin{equation*}
A(s;\lambda,\tilde\rho)-A(s;\lambda_0,0) =
r'(s)^2(\lambda-\lambda_0-\tilde\rho)B_0 e^{2s}\,(\lambda-\lambda_0-\tilde\rho)B_0,
\end{equation*}
as $r'(s)^2=e^{2s}$ for $s \le s_1$ (see~\eqref{E:def B_0} for the
definition of $B_0$).

The function $e^\tau \rho(\tau)$ belongs to $L^2(J_-, H^{1/2})$. By Lemma 5,
$e^\tau \rho(\tau)u_1(\tau)$ also belongs to this space. Thus
$\| e^\tau \rho(\tau) B_0 U \|_X \in  L^2 (J_-)$ and
$\| e^{2\tau} \rho(\tau) B_0 U \|_X$ is the product of an $L^2$ function
and the exponentially decaying function $e^{\tau}$.

This allows us
to use the Gap Lemma as in \cite[\S4.3 and (4.12)]{S02} and
\cite[Proof of Lemma~4.1]{BSZ08} to show that \eqref{E:perturbed system}
has two linearly independent solutions
$U^{cb}_{0,j}(s;\lambda,\tilde\rho)$ for $j=1,3$ that converge to $e_j$
as $s\to-\infty$, and two other solutions which grow algebraically.
In fact, the results in these works show that any linear combination
of the bounded solutions $U^{cb}_{0,j}(s;\lambda,\tilde\rho)$ with $j=1,3$
can be found as a fixed point of the equation
\begin{eqnarray}\label{e:gl}
  U(s) & = & \Phi^{cu}_-(s;\lambda_0,0) U_0^{cb} +
  \int_{-\infty}^s \Phi^{cs}_-(s,\tau;\lambda_0,0)
  \,e^{2\tau}\,
(\lambda-\lambda_0-\tilde\rho(\tau))B_0 U(\tau)\, d\tau \\
  \nonumber &&
  - \int_{s}^{s_1} \Phi^{u}_-(s,\tau;\lambda_0,0) \,e^{2\tau}\,
  (\lambda-\lambda_0-\tilde\rho(\tau))B_0  U(\tau)\, d\tau,
\end{eqnarray}
where $U_0^{cb}$ belongs to the unperturbed bounded central subspace
spanned by $U_{0,j}(s_1)$ for $j=1,3$. We denote the fixed point
by $U^{cb}_-(s;\lambda,\tilde\rho,U_0^{cb})$ and
write
\begin{equation}\label{e:ipc} P^{cb}_-(s_1;\lambda,\tilde\rho)
  U_0^{cb} := U^{cb}_-(s_1;\lambda,\tilde\rho,U_0^{cb}) = U_0^{cb} +
  \int_{-\infty}^{s_1} \Phi^{cs}_-(s_1,\tau;\lambda_0,0)\,e^{2\tau}\,
  (\lambda-\lambda_0-\tilde\rho(\tau))B_0
  U^{cb}_-(\tau;\lambda,\tilde\rho,U_0^{cb})\, d\tau.
\end{equation}
Similarly, we can use \eqref{E:dichotomyequation} to describe the
solutions of \eqref{E:perturbed system} with exponential decay
as $s\to-\infty$ by
\begin{equation}\label{e:ipu}
  P_-^u(s_1;\lambda,\tilde \rho) = P_-^u(s_1;\lambda_0,0) +
  \int_{-\infty}^{s_1} \Phi_-^{cs} (s_1,\tau;\lambda_0,0)\,e^{2\tau}\,
  (\lambda-\lambda_0-\tilde\rho(\tau)) B_0
  \Phi_-^u(\tau,s_1;\lambda,\tilde \rho)\, d\tau.
\end{equation}
These results will be used later to characterize eigenfunctions of the
perturbed operator.
\end{section}


\begin{section}{{Dichotomies for the far field}}\label{S:s large}
  The method of Section~\ref{S:r_small} is not available for
  determining the dichotomies for $s$ large.  Going back
  to~\eqref{E:eigenvalue_equation}, we observe that $\theta$ and
  $\rho$ have support in a ball with radius $r_1$, and thus for $r\geq
  r_1$ the eigenvalue problem~\eqref{E:eigenvalue_equation} reduces to
  $(\Delta^2 - \lambda) u = 0$, which can be factorized:
\begin{equation*}
   (\Delta - \sqrt{\lambda})(\Delta + \sqrt{\lambda}) u =
   (\Delta + \sqrt{\lambda})(\Delta - \sqrt{\lambda}) u=0.
\end{equation*}
Expanding $u(r,\varphi)$ as a Fourier series in the angular variable
$\varphi$, we see that the Fourier coefficients $\widehat u_k$ satisfy
the differential equations
\begin{equation}\label{E:factorised}
  \left(\frac{\partial^2}{\partial r^2} + \frac{1}{r}
    \frac{\partial}{\partial r} - \frac{k^2}{r^2} - \sqrt{\lambda}\right)
  \left(\frac{\partial^2}{\partial r^2} + \frac{1}{r}
    \frac{\partial}{\partial r} - \frac{k^2}{r^2} +
    \sqrt{\lambda}\right) \widehat u_k = 0.
\end{equation}
For $k$ fixed, this is a fourth order linear ODE, so it has a
four-dimensional space of solutions. The general solution can then be
obtained as a linear combination of the solutions of
\begin{equation}\label{E:modified Bessel}
  \left(\frac{\partial^2}{\partial r^2} + \frac{1}{r}
    \frac{\partial}{\partial r} - \frac{k^2}{r^2} -
    \sqrt{\lambda}\right) \widehat u_k = 0
\end{equation}
and the solutions of
\begin{equation}\label{E:Bessel}
  \left(\frac{\partial^2}{\partial r^2} + \frac{1}{r}
    \frac{\partial}{\partial r} - \frac{k^2}{r^2} +
    \sqrt{\lambda}\right)\widehat u_k = 0,
\end{equation}
so that the general solution of~\eqref{E:factorised} is given by
\begin{equation*}
   \widehat u_k(r) = C_1 I_k(\lambda^{1/4} r) + C_2 K_k(\lambda^{1/4}
   r) + C_3 J_k(\lambda^{1/4} r) + C_4 Y_k(\lambda^{1/4} r),
\end{equation*}
where $J_k$ and $Y_k$ are Bessel functions of the first and second
kind, respectively, which satisfy equation \eqref{E:Bessel}, and $I_k$ and $K_k$ are
modified Bessel functions of the first and second kind, respectively, which satisfy
equation \eqref{E:modified Bessel}.

For $r\geq r_3=s_3$, we have $s=r$. Thus, we can define the
systems corresponding to equations \eqref{E:modified Bessel} and
\eqref{E:Bessel} with the variable~$s$ for $s\geq s_3$ as
\begin{equation}\label{E:modified Bessel-system}
  \begin{aligned}
    u_1' &= u_2, \\
    u_2' &= \left(\frac{k^2}{s^2} + \sqrt{\lambda}\right) u_1
    -\frac{1}{s} u_2,
  \end{aligned}
\end{equation}
and
\begin{equation}\label{E:Bessel-system}
  \begin{aligned}
    u_1' &= u_2, \\
    u_2' &= \left(\frac{k^2}{s^2} - \sqrt{\lambda}\right) u_1
    -\frac{1}{s} u_2,
  \end{aligned}
\end{equation}
respectively.  We will consider these systems for $s\geq r_1$ and
define $\phi_k(s,t)$ and $\psi_k(s,t)$ to be the evolution
operators corresponding to the systems~\eqref{E:modified
  Bessel-system} and~\eqref{E:Bessel-system}, respectively.
After deriving dichotomies for those systems, we will derive
dichotomies for the original system~\eqref{E:perturbed system}.

To derive dichotomies for~\eqref{E:modified Bessel-system}
and~\eqref{E:Bessel-system}, we introduce for $s\geq r_1$ the
function spaces $\tilde X^s = H^1(S^1) \times L^2(S^1)$, $\tilde Y^s = H^2(S^1)\times H^1(S^1)$ and
  $\tilde Z^s = H^3(S^1)\times H^2(S^1)$ with norms
\begin{equation*}
   \begin{aligned}
      \|\underline u\|_{\tilde X^s}^2 &:= \frac{1}{s^2}
      \|u_1\|_{H^1(S^1)}^2 + \|u_1\|_{L^2(S^1)}^2 + \|u_2\|_{L^2(S^1)}^2, \\
      \|\underline u\|_{\tilde Y^s}^2 &:= \frac{1}{s^2}
      \|u_1\|_{H^2(S^1)}^2 + \|u_1\|_{H^1(S^1)}^2 + \|u_2\|_{H^1(S^1)}^2, \\
      \|\underline u\|_{\tilde Z^s}^2 &:= \frac{1}{s^2}
      \|u_1\|_{H^3(S^1)}^2 + \|u_1\|_{H^2(S^1)}^2 + \|u_2\|_{H^2(S^1)}^2.
   \end{aligned}
\end{equation*}

We decompose the spaces $\tilde X^s$, $\tilde Y^s$, and
$\tilde Z^s$ into their Fourier subspaces
\begin{equation}\label{E:X^s-decomposition}
\tilde X^s = \overline{\oplus_{k\in \mathbb Z} \tilde X_k^s}, \quad
\tilde Y^s = \overline{\oplus_{k\in \mathbb Z} \tilde Y_k^s},
\quad\mbox{and}\quad
\tilde Z^s = \overline{\oplus_{k\in \mathbb Z} \tilde Z_k^s},
\end{equation}
where
\begin{equation*}
   \tilde X_k^s = \tilde Y_k^s = \tilde Z_k^s =
\{(a e^{ik\cdot},b e^{ik\cdot})^T;\; a, b \in \mathbb C\},
\end{equation*}
and the completion in \eqref{E:X^s-decomposition} is in the
respective norms of $\tilde X^s$, $\tilde Y^s$ and $\tilde Z^s$.  The
norms on $\tilde X_k^s$, $\tilde Y_k^s$ and $\tilde Z_k^s$ are given
by the restriction of the norms of $\tilde X^s$, $\tilde Y^s$, and
$\tilde Z^s$ respectively, and so
\begin{equation}\label{E:X^s_k norms}
   \begin{aligned}
      \|(a e^{ik\cdot},b e^{ik\cdot})^T\|^2_{\tilde X_k^s} &=
      \left(1 + \frac{k^2}{s^2}\right) |a|^2 + |b|^2; \\
      \|(a e^{ik\cdot},b e^{ik\cdot})^T\|^2_{\tilde Y_k^s} &      \left(1+k^2\right)\,
      \|(a e^{ik\cdot},b e^{ik\cdot})^T\|^2_{\tilde X_k^s};\\
      \|(a e^{ik\cdot},b e^{ik\cdot})^T\|^2_{\tilde Z_k^s} &      \left(1+k^2\right)^2\,
      \|(a e^{ik\cdot},b e^{ik\cdot})^T\|^2_{\tilde X_k^s}.
   \end{aligned}
\end{equation}

For each $\epsilon\in (0,\lambda^{1/4})$, we now prove the existence of a
{time-dependent} exponential dichotomy for
\eqref{E:modified Bessel-system} with constant $K>0$ and rates
$\kappa^s=-(\lambda^{1/4}-\epsilon)$ and $\kappa^u=(\lambda^{1/4} -
\epsilon)$ that are independent of $k$.  For
\eqref{E:Bessel-system}, we will show that the evolution operator always acts in
the center-unstable manifold and derive that its growth can
be bounded by any exponential.
\begin{Lem}\label{L:modified Bessel}
  There exists an $\epsilon_0>0$ such that for any $\epsilon\in
  (0,\epsilon_0)$ there exists a $K>0$ such that for any
  $k\in\mathbb{Z}$ and $\lambda\in (\lambda_0/2,2\lambda_0)$
  there exists a {time-dependent} exponential dichotomy of~\eqref{E:modified
    Bessel-system} on $J_+$ so that $\phi^s_k(s,t;\lambda)$ and
  $\phi^u_k(s,t;\lambda)$
   satisfy
   \begin{equation*}
      \begin{array}{r l l}
         \|\phi_k^{s}(s,t;\lambda)\|_{\mathcal L(\tilde X_k^t,\tilde
           X_k^s)}
         =\|\phi_k^{s}(s,t;\lambda)\|_{\mathcal L(\tilde Y_k^t,\tilde
           Y_k^s)}
         &\le K e^{-(\lambda^{1/4}-\epsilon)(s-t)} \qquad & s\ge t\ge r_1, \\
         \|\phi_k^{u}(s,t;\lambda)\|_{\mathcal L(\tilde X_k^t,\tilde
           X_k^s)} = \|\phi_k^{u}(s,t;\lambda)\|_{\mathcal L(\tilde Y_k^t,\tilde
           Y_k^s)} &\le K e^{-(\lambda^{1/4}-\epsilon)(t-s)} \qquad &
         t\ge s \ge r_1.
      \end{array}
   \end{equation*}
\end{Lem}
\begin{proof}
  Let $(u_1,u_2)^T$ satisfy equation \eqref{E:modified Bessel-system}.
  To get estimates which are uniform in~$k$, we follow \cite{aS98}, and let
   \begin{equation*}
      \tilde u_1(s) :=
      \left(\sqrt{\lambda} + \frac{k^2}{s^2}\right)^{1/2} u_1(s).
   \end{equation*}
   Note that
   \begin{equation*}
   \min(1,\sqrt{\lambda_0/2})\|(u_1,u_2)^T\|_{\tilde X_k^s} \le
   \|(\tilde u_1,u_2)^T\|_{\mathbb C^2} \le
   \max(1,\sqrt{2\lambda_0}) \|(u_1,u_2)^T\|_{\tilde X_k^s},
   \end{equation*}
   and that the constants above are independent of $k$ and $\lambda$.
   This shows that, when using the new variables $\tilde u_1$ and
   $u_2$, we can use the standard norm in $\mathbb C^2$.

   Next, we rewrite the system \eqref{E:modified Bessel-system} in the
   new variables $\tilde u_1$, $u_2$:
   \begin{equation*}
      \begin{aligned}
         \tilde u_1' &= \left(\sqrt{\lambda} +
           \frac{k^2}{s^2}\right)^{1/2} u_2 -
         \frac{k^2}{s^3} \left(\sqrt{\lambda} +
         \frac{k^2}{s^2}\right)^{-1} \tilde u_1, \\
         u_2' &= -\frac{1}{s} u_2 + \left(\sqrt{\lambda} +
           \frac{k^2}{s^2}\right)^{1/2} \tilde u_1,
      \end{aligned}
   \end{equation*}
   Now, we change the independent variable by making the substitution
   $d\tau/ds = (\sqrt{\lambda} + k^2/s^2)^{1/2}$. We write $s(\tau)$
   to describe the dependence of $s$ on $\tau$. We then obtain (where
   $'$ now denotes differentiation with respect to $\tau$)
   \begin{equation}\label{E:a}
      \begin{aligned}
         \tilde u_1' &= -  \frac{k^2}{s(\tau)^3} \left(\sqrt{\lambda} + \frac{k^2}{s(\tau)^2}\right)^{-3/2} \tilde u_1 + u_2, \\
         u_2' &= \tilde u_1 - \frac{1}{s(\tau)} \left(\sqrt{\lambda} + \frac{k^2}{s(\tau)^2}\right)^{-1/2} u_2. \\
      \end{aligned}
   \end{equation}
   Noting that $s(\tau)\to \infty$ as $\tau\to\infty$ we find that the
   limiting system at $+\infty$ is
   \begin{equation}\label{E:system at +infinity}
      \begin{aligned}
         \tilde u_1' &= u_2, \\
         u_2' &= \tilde u_1,
      \end{aligned}
   \end{equation}
   which is independent of $k$. The matrix associated with this system
   has eigenvalues $\pm 1$.  Hence equation \eqref{E:system at
     +infinity} possesses exponential dichotomies with $\kappa^u  = -\kappa^s = 1$. To get estimates for the perturbated
   system~\eqref{E:a},  we will use the estimates
   \begin{equation*}
      \begin{aligned}
        \left|\frac{k^2}{s(\tau)^3} \left(\sqrt{\lambda} +
            \frac{k^2}{s(\tau)^2}\right)^{-3/2}\right| &= \frac{1}{s(\tau)} \frac{1}{\sqrt{\sqrt{\lambda} +
          k^2/s(\tau)^2}} \frac{k^2/s(\tau)^2}{\sqrt{\lambda} +
          k^2/s(\tau)^2} \le
        \frac{1}{\lambda^{1/4}s(\tau)}, \\
        \left|\frac{1}{s(\tau)} \left(\sqrt{\lambda} +
            \frac{k^2}{s(\tau)^2}\right)^{-1/2}\right| &\le
        \frac{1}{\lambda^{1/4}s(\tau)}.
      \end{aligned}
   \end{equation*}
   This estimate is uniform in $\lambda$ in a neighbourhood of
   $\lambda_0$. The roughness theorem for exponential dichotomies
   \cite[Chapter 4]{wC78} now guarantees the existence of an exponential
   dichotomy also for the system \eqref{E:a}, and we denote the
   corresponding evolution operators by $\tilde \phi_k^s(\sigma,\tau)$
   and $\tilde \phi_k^u(\sigma,\tau)$.  For each positive $\tilde
   \epsilon$ sufficiently small, there exists a $K\ge 0$ such that
   \begin{equation*}
      \begin{array}{r l l}
         \|\tilde \phi_k^s(\sigma,\tau)\|_{\mathcal L(\mathbb C^2)}&\le K e^{-(1-\tilde \epsilon)(\sigma - \tau)}\quad & \sigma\ge \tau, \\
         \|\tilde \phi_k^u(\sigma,\tau)\|_{\mathcal L(\mathbb C^2)}&\le K e^{-(1-\tilde \epsilon)(\tau - \sigma)}\quad & \sigma\le \tau.
      \end{array}
   \end{equation*}
   Moreover, $K$ does not depend on $\lambda$ in a neighbourhood of $\lambda_0$ or on $k\in\mathbb{Z}$.

   It remains to translate this result back to the $s$ variable. We
   write $s=s(\sigma)$ and $t=s(\tau)$. Note that $ds/d\tau\le
   \lambda^{1/4}$, and so by the chain rule we have for $s>t$
   \begin{equation*}
      \begin{aligned}
         \|\phi_k^s(s,t)\|_{\mathcal L(\tilde X_k^t,\tilde X_k^s)} &\le C \|\tilde \phi_k^s(\sigma,\tau)\|_{\mathcal L(\mathbb C^2)}
         \le K e^{-(1 - \tilde \epsilon)(\sigma - \tau)} \\
         &\le K e^{-(1-\tilde \epsilon)\lambda^{1/4}(s-t)} \le K e^{-(\lambda^{1/4} - \epsilon)(s-t)},
      \end{aligned}
   \end{equation*}
   where we have put $\epsilon = \tilde \epsilon \lambda^{1/4}$. A
   similar calculation proves that for $t>s$
   \begin{equation*}
      \|\phi_k^u(s,t)\|_{\mathcal L(\tilde X_k^t,\tilde X_k^s)} \le K e^{-(\lambda^{1/4} - \epsilon)(t-s)}.
   \end{equation*}

   The estimates for $\tilde Y^s_k$ also follow from these estimates,
   since it is only a matter of multiplying both sides of the inequalities by a factor $(1+k^2)$.
\end{proof}

\begin{Lem}\label{L:Bessel}
   Let $\epsilon>0$ be given.
   Then there exists a $K>0$ such that for any $k\in\mathbb{Z}$ and
   $\lambda\in (\lambda_0/2,2\lambda_0)$ we have
   \begin{equation*}
     \|\psi_k(s,t;\lambda)\|_{\mathcal L(\tilde X_k^t,\tilde X_k^s)} = \|\psi_k(s,t;\lambda)\|_{\mathcal L(\tilde Y_k^t,\tilde Y_k^s)}
     \le K e^{\epsilon(t-s)} \quad\mbox{for}\quad t\ge s\ge r_1.
   \end{equation*}
\end{Lem}
\begin{proof}
  First we note a scaling invariance in~\eqref{E:Bessel-system}. If
  $(\bar u_1(s,t),\bar u_2(s,t))$ is a solution of~\eqref{E:Bessel-system}
  with $\lambda=1$, then $(u_1(s,t),u_2(s,t)) =
  (\bar  u_1(\lambda_1^{1/4}s,\lambda_1^{1/4}t),
  \lambda_1^{1/4}\bar u_2(\lambda_1^{1/4}s,\lambda_1^{1/4}t))$ is a
  solution of~\eqref{E:Bessel-system} with $\lambda=\lambda_1$.  So it
  is sufficient to prove the estimate in the Lemma in case
  $\lambda=1$. Using the explicit expressions for the solution in
  terms of Bessel function, it follows that, for $\lambda=1$,
  $\psi_k(s,t)$ is given by
   \begin{equation*}
      \begin{aligned}
         \psi_k(s,t) &= \left(
         \begin{matrix}
           J_k(s) & Y_k(s) \\
           J_k'(s) & Y_k'(s)
         \end{matrix}
         \right)\left(
         \begin{matrix}
           J_k( t) & Y_k( t) \\
            J_k'( t) &  Y_k'(t)
         \end{matrix}
         \right)^{-1} \\
         &= \frac{\pi t}{2} \left(
         \begin{matrix}
             J_k( s) & Y_k( s) \\
             J_k'( s) & Y_k'( s)
         \end{matrix}
         \right)\left(
         \begin{matrix}
            Y_k'( t) & -Y_k( t) \\
            - J_k'( t) & J_k(t)
         \end{matrix}
         \right) \\
         &= \frac{\pi t}{2} \left(
         \begin{matrix}
            a_k(s,t) & b_k(s,t) \\
            c_k(s,t) & d_k(s,t)
         \end{matrix}
         \right),
      \end{aligned}
   \end{equation*}
   where
   \begin{equation*}
      \begin{aligned}
        a_k(s,t) &:=  (J_k( s) Y_k'(t) - Y_k( s) J_k'( t)), \\
        b_k(s,t) &:= -J_k( s) Y_k( t) + Y_k(s) J_k( t), \\
        c_k(s,t) &:=  (J_k'( s) Y_k'( t) -  Y_k'( s) J_k'( t)), \\
        d_k(s,t) &:=  (-J_k'( s) Y_k( t) + Y_k'( s) J_k( t)),
      \end{aligned}
   \end{equation*}
   and we have used that the Wronskian of $J_k(t)$ and $Y_k(t)$ is
   $\frac{2}{\pi t}$ \cite[(9.1.16)]{AS72}.
   Writing $\underline u= (u_1,u_2)^T\in \widetilde X_k^t$, we have
   \begin{equation*}
      \| \psi_k(s,t) \underline u\|_{\widetilde X_k^s}^2 = \frac{\pi^2 t^2}{4}\left( \left(1 + \frac{k^2}{s^2}\right) (a_k(s,t) u_1 + b_k(s,t) u_2)^2 + (c_k(s,t) u_1 + d_k(s,t) u_2)^2\right),
   \end{equation*}
   and so we need to show that there exists a $K>0$ such that for every $u_1, u_2\in \mathbb R$, $k\in \mathbb Z$ and
   $t\ge s\ge r_1$,
   \begin{equation*}
      t^2\left( \left(1 + \frac{k^2}{s^2}\right) (a_k(s,t) u_1 + b_k(s,t) u_2)^2 + (c_k(s,t) u_1 + d_k(s,t) u_2)^2\right) \le K^2 e^{2\epsilon (t-s)}\left(\left(1 + \frac{k^2}{t^2}\right)u_1^2 + u_2^2\right).
   \end{equation*}

    By choosing $u_1$ and $u_2$ appropriately, we note that this inequality holds if and only if the following
    two inequalities hold for some $K>0$ and all
   $k\in \mathbb Z$ and $t>s>r_1$:
   \begin{equation*}
      \begin{aligned}
         \left(\left(1 + \frac{k^2}{s^2}\right)a_k(s,t)^2 + c_k(s,t)^2\right) t^2 &\le K^2 e^{2\epsilon (t-s)}\left(1 + \frac{k^2}{t^2}\right), \\
         \left(\left(1 + \frac{k^2}{s^2}\right)b_k(s,t)^2 + d_k(s,t)^2\right) t^2 &\le K^2 e^{2\epsilon (t-s)}.
      \end{aligned}
   \end{equation*}
   To simplify further, we note that the above two inequalities hold if  there exists a constant $K>0$ such that
   for $t\ge s \ge r_1$,
   \begin{equation}\label{E:4 inequalities}
      \begin{aligned}
         \sqrt{1 + \frac{k^2}{s^2}} |a_k(s,t)| t \le K e^{\epsilon(t-s)}\sqrt{1 + \frac{k^2}{t^2}}, \\
         \sqrt{1 + \frac{k^2}{s^2}} |b_k(s,t)| t \le K e^{\epsilon(t-s)}, \\
         |c_k(s,t)| t \le K e^{\epsilon(t-s)}\sqrt{1 + \frac{k^2}{t^2}}, \\
         |d_k(s,t)| t \le K e^{\epsilon(t-s)}.
      \end{aligned}
   \end{equation}

   First we will  prove the fourth inequality of \eqref{E:4 inequalities}. Let
   \begin{equation*}
      f_k(s,t) := e^{-\epsilon(t-s)}t d_k(s,t).
   \end{equation*}
   We need to show that $f_k(s,t)$ is uniformly bounded for $k\in
   \mathbb Z$ and $t\ge s \ge r_1$. Since $J_{-k}= (-1)^k J_k$ and
   $Y_{-k}=(-1)^k Y_k$, it is sufficient to consider $k\in \mathbb N$.

 We start with showing that $f_k$ is bounded for $k\in\mathbb{N}$ fixed.
In the slightly smaller sector $t\ge (1+\delta)s\ge (1+\delta)r_1$ (where
$\delta>0$ is arbitrary), we have $|f_k(s,t)|\to 0$ as $s^2+t^2\to
\infty$, {or equivalently, as $t\to \infty$}. Indeed, $Y_k(s)$, $J_k(s)$, $Y_k'(s)$ and $J_k'(s)$ are
bounded by a constant $C_k$ for $s\ge r_1$ \cite[(9.2.1)]{AS72}, and so
   \begin{equation}\label{E:f_k k-dep estimate}
     |f_k(s,t)|\le C_k^2 e^{-\epsilon(t-s)}t
     \le {C_k^2 t e^{-\epsilon \delta t/(1+\delta)}} \to 0
   \end{equation}
   as {$t\to \infty$}. Furthermore, $\sqrt{s}Y_k(s)$,
   $\sqrt{s}J_k(s)$, $\sqrt{s}Y_k'(s)$ and $\sqrt{s}J_k'(s)$ are
   bounded by a constant $D_k$ for $s\ge r_1$ \cite[(9.2.1)]{AS72}, and so
   for $r_1\le s\le t\le (1+\delta) s$ we have
\[
|f_k(s,t)| \le e^{-\epsilon(t-s)}\sqrt{(1+\delta) s t} \,d_k(s,t)
\leq \sqrt{1+\delta} \,D_k^2.
\]
   Altogether this implies that $f_k(s,t)$
   is bounded in the whole sector $t\ge s\ge r_1$ by a constant, possibly
   depending on $k$.

To show that in fact $f_k(s,t)$ is  bounded by a
$k$-independent constant, we consider $s\geq r_1$ as being fixed for the
moment.  First we note that in \eqref{E:f_k k-dep estimate}, we proved that
 for fixed $s$ and $k$, $f_k(s,t)\to 0$ as $t\to\infty$.
Thus for $s\geq r_1$ fixed, the
function $f_k(s,t)$ attains its maximum in an interior point $t>s$ or
at the boundary $t=s$.  We use a method by L. Landau \cite{ljL00} to
analyze the behaviour of $f_k(s,\cdot)$ at its critical points.  At
a critical point, we have $\partial f_k/\partial t(s,t) = 0$.  By
equation (11) of \cite{ljL00}, at the points where $\partial
f_k/\partial t = 0$ we have
   \begin{equation*}
      \frac{\partial}{\partial k} f_k(s,t)^2 =
      2 t \frac{f_k(s,t)^2}{e^{-2\epsilon}(t-s) t^2}
      \frac{\partial}{\partial t}\left( e^{-2\epsilon(t-s)} t^2 A_k(t)\right),
   \end{equation*}
   where $A_k(t) = \int_0^\infty K_0(2t\cosh \tau) e^{-2k\tau}\,
   d\tau$, and $K_0$ is a modified Bessel function of the second kind,
   satisfying
   \begin{equation*}
      K_0(x) = \int_0^\infty e^{-x\cosh \tau}\, d\tau.
   \end{equation*}
   In particular (since $2t f_k(s,t)^2/e^{2\epsilon(t-s)}t^2 >0$),
   $f_k(s,t)^2$ is decreasing in $k$ at a point where $\partial
   f_k/\partial t=0$ if and only if $e^{-2\epsilon(t-s)} t^2 A_k(t)$
   is decreasing in $t$.  Note that $A_k(t)$ is monotonically
   decreasing for $t>0$, and since $e^{-2\epsilon(t-s)} t^2$ is
   monotonically decreasing for $t>1/\epsilon$, we conclude that
   $|f_k(s,\cdot)|$ is monotonically decreasing in $k$ at its critical
   points for $t>1/\epsilon$.

   {
   Note that in the case when $r_1>1/\epsilon$,} we have proved that
   if the maximum of $f_k(s,\cdot)$ occurs for $t>s\ge r_1$, then the
   maximum is decreasing in~$k$, and hence stays bounded as $k$ increases.
   At the boundary $t=s\ge r_1$, we have $f_k(s,s)=2/\pi$, which is independent of
   $k$. As each function $f_k(s,t)$ is bounded, in particular
   $f_1(s,t)$ is bounded, it follows that $f_k(s,t)$ is bounded in the
   whole sector $t\ge s\ge r_1$ by a $k$-independent constant for all
   $k\in\mathbb{Z}$. {This shows that $f_k$ is uniformly bounded in the case
   when $r_1>1/\epsilon$.}

   {When $r_1<1/\epsilon$,}
   we also need to estimate $f_k(s,t)$ in the triangle
   $1/\epsilon>t>s\geq r_1$. Here we use the estimate
   $e^{-\epsilon(t-s)}t
   \le 1/\epsilon$. It follows that $|f_k(s,t)|\le |g_k(s,t)|$, where
   \begin{equation*}
     g_k(s,t) :=     \frac 1{\epsilon}
     (Y_k'(s)J_k(t) - J_k'(s) Y_k(t)).
   \end{equation*}
   Applying Section 3 of \cite{ljL00} we conclude that $g_k(s,t)^2$ is
   decreasing in $k$ at the points where $\partial g_k(s,t)/\partial t
   = 0$. Furthermore, $g_k(s,t)\to0$, for $t\to\infty$ and $g_k(s,s) = \frac 2{\pi\epsilon s}\leq \frac 2{\pi\epsilon r_1}$ for $s\geq
   r_1$.
The proof of the fourth inequality of \eqref{E:4
     inequalities} is complete.

   Next, we prove the second inequality of \eqref{E:4 inequalities}.
   By \cite[(9.1.27)]{AS72} we have
   \begin{equation*}
      \begin{aligned}
         \left(1 + \frac{k}{s}\right) Y_k(s) &=
         Y_k(s) + \frac{1}{2}(Y_{k-1}(s) + Y_{k+1}(s)), \\
         \left(1 + \frac{k}{s}\right) J_k(s) &=
         J_k(s) + \frac{1}{2}(J_{k-1}(s) + J_{k+1}(s)).
      \end{aligned}
   \end{equation*}
   Note that
   \begin{equation*}
      \left(1 + \frac{k^2}{s^2}\right)^{1/2} \le
      1 + \frac{k}{s} \le \sqrt{2} \left(1+\frac{k^2}{s^2}\right)^{1/2},
   \end{equation*}
   and so the second inequality of \eqref{E:4 inequalities} is equivalent to
   \begin{equation*}
    \left|  t\left( \bigl( Y_k(s) + \frac{1}{2} Y_{k-1}(s) +
        \frac{1}{2}Y_{k+1}(s)\bigr)J_k(t) - \bigl(J_k(s) + \frac{1}{2}
      J_{k-1}(s) + \frac{1}{2}J_{k+1}(s)\bigr)Y_k(t)\right)\right| \le
    K e^{\epsilon(t-s)}.
   \end{equation*}

   To prove the inequality, we use the same method as above, with
   \begin{equation*}
      f_k(s,t) := t e^{-\epsilon(t-s)}\left(\bigl(Y_k(s) +
        \frac{1}{2} Y_{k-1}(s) + \frac{1}{2}Y_{k+1}(s)\bigr)J_k(t) -
        \bigl(J_k(s) + \frac{1}{2} J_{k-1}(s) +
        \frac{1}{2}J_{k+1}(s)\bigr)Y_k(t)\right)
   \end{equation*}
   and
   \begin{equation*}
      g_k(s,t) := t e^{-s}\frac{e^{-1}}{\epsilon}\left(\bigl( Y_k(s) +
        \frac{1}{2} Y_{k-1}(s) + \frac{1}{2}Y_{k+1}(s)\bigr)J_k(t) -
        \bigl(J_k(s) + \frac{1}{2} J_{k-1}(s) +
        \frac{1}{2}J_{k+1}(s)\bigr)Y_k(t)\right).
   \end{equation*}
   Then $|f_k(s,t)|\to 0$ as $s^2+t^2\to \infty$ in the sector $t\ge s\ge
   r_1$ and as above, $f_k(s,t)^2$ is decreasing in $k$ at the points
   where $\partial f_k/\partial t=0$ if $t> 1/\epsilon$.  In the
   triangle $1/\epsilon\ge t\ge s\ge r_1$, $|f_k(s,t)|\le |g_k(s,t)|$ and
   $g_k$ is decreasing in $k$ at the points where $\partial
   g_k/\partial t= 0$. Finally, at the boundary points where $s=t\ge
   r_1$, we have $f_k(s,s)=g_k(s,s)=0$. We conclude that the second
   inequality is valid in the whole sector $t\ge s\ge r_1$.

   For the third inequality of \eqref{E:4 inequalities}, we use the
   last identity of \cite[(9.1.27)]{AS72} which shows that the
   inequality is equivalent to the two inequalities
   \begin{equation*}
      \begin{aligned}
         \frac{k}{t} |J_k'(s) Y_k(t) - Y_k'(s) J_k(t)| t & \le
         K e^{\epsilon (t-s)} \left(1 + \frac{k^2}{t^2}\right)^{1/2}, \\
         |J_k'(s) Y_{k+1}(t) - Y_k'(s) J_{k+1}(t)|t &\le
         K e^{\epsilon (t-s)}\left(1 + \frac{k^2}{t^2}\right)^{1/2}.
      \end{aligned}
   \end{equation*}
   The first of these inequalities follows from the fourth inequality
   of \eqref{E:4 inequalities}, and the second can be handled as in
   the proof of the second and fourth inequalities after noting that
   at the boundary where $t=s\ge r_1$ we have
   \begin{equation*}
      |J_k'(s)Y_{k+1}(s) - Y_k'(s)J_{k+1}(s)| =
      \frac{k}{s}|J_k'(s) Y_k(s) - Y_k'(s) J_k(s)| = \frac{2k}{\pi s^2},
   \end{equation*}
   thus $f_k(s,s)= \frac{2}{\pi\sqrt{1+s^2/k^2}} \leq \frac 2\pi$.  We
   omit the details.

   It remains to prove the first inequality of \eqref{E:4
     inequalities}. This can be handled as the second inequality by
   using
   \begin{equation*}
      J_k(s) Y_k'(t) - Y_k(s)J_k'(t) = \frac{k}{t}(J_k(s)Y_k(t) -
      Y_k(s)J_k(t)) + (Y_k(s)J_{k+1}(t) - J_k(s)Y_{k+1}(t))
   \end{equation*}
   and splitting the inequality up into the two inequalities
   \begin{equation*}
      \begin{aligned}
         \left(1 + \frac{k^2}{s^2}\right)^{1/2} \frac{k}{t}
         |J_k(s)Y_k(t) - Y_k(s) J_k(t)|t &\le K e^{\epsilon(t-s)}
         \left(1 + \frac{k^2}{t^2}\right)^{1/2}, \\
         |Y_k(s) J_{k+1}(t) - J_k(s)Y_{k+1}(t)|
         \left(1 +  \frac{k^2}{s^2}\right)^{1/2} t &\le
         K e^{\epsilon(t-s)}\left(1 + \frac{k^2}{t^2}\right)^{1/2}.
      \end{aligned}
   \end{equation*}
   The first inequality follows directly from the second inequality of
   \eqref{E:4 inequalities}, and the second inequality is proved in
   the same way as the second inequality of \eqref{E:4 inequalities},
   except that at the boundary where $t=s$ we have
   \begin{equation*}
      \left|Y_k(s) J_{k+1}(s) - J_k(s) Y_{k+1}(s)\right| = 2/\pi s.
   \end{equation*}
   The details are omitted.

 It is also clear that the estimates above are uniform in $\lambda$
 for $\lambda\in (\lambda_0/2,2\lambda_0)$.

   The estimates for $\tilde Y^s_k$ follow by the same estimates, since
   it is just a matter of multiplying each side of the inequalities by
   the factor $(1+k^2)$.
\end{proof}

We are ready to prove that there exist {time-dependent} exponential dichotomies on
$J_+$ for the full system~\eqref{E:perturbed system}.  First we define
the spaces $\mathcal X^s:= \tilde X^s\times \tilde X^s$ and $\mathcal
Y^s:= \tilde Y^s\times \tilde Y^s$. Note that $\mathcal Y^s\subset X
\subset\mathcal X^s$.  As before, we can decompose those spaces into
$\mathcal X^s = \overline{\oplus_{k\in \mathbb Z}\mathcal X_k^s}$ and
$\mathcal Y^s = \overline{\oplus_{k\in \mathbb Z}\mathcal Y_k^s}$ with
$\mathcal X_k^s=\tilde X_k^s\times \tilde X_k^s$ and $\mathcal
Y_k^s=\tilde Y_k^s\times \tilde Y_k^s$.

For $s>s_3=r_3$, we have that $r=s$ and hence, for $s\in[s_3,\infty)$, the
system~\eqref{E:perturbed system} reduces to
\begin{equation}\label{E:far field system}
   U' = A(s;\lambda) U \mbox{~~with~~} A(s;\lambda) =
   \begin{pmatrix}
     0&1&0&0\\
     - \frac{1}{s^2} \partial^2 & - \frac{1}{s} & 1 &0\\
     0&0&0&1\\
     \lambda &0&- \frac{1}{s^2} \partial^2&- \frac{1}{s}
   \end{pmatrix}.
\end{equation}
We consider the system~\eqref{E:far field system} for
$s\in[r_1,\infty)$ and record that \eqref{E:far field system} and \eqref{E:perturbed system} coincide on the smaller
interval~$[s_3,\infty)$. The exponential dichotomy
for~\eqref{E:perturbed system} on the whole interval $J_+$ will follow
from the fact that the systems~\eqref{E:perturbed system}
and~\eqref{E:far field system} are linked by the smooth transformation $r(s)$ of the independent variable on the compact interval $[s_1,s_3]$. Thus, since we are
using $r=s$ in~\eqref{E:far field system}, we see that if $\widetilde{U}(s)$ is a
solution of~\eqref{E:far field system} for $s\in[r_1,\infty)$, then
$U(s)=\mathrm{diag}(1,r'(s),1,r'(s))\,\widetilde{U}(r(s))$ is a solution
of~\eqref{E:perturbed system} for all $s\in[s_1,\infty)$. Recall that
there are constants $0<c<C$ such that $c\leq r'(s)\leq C$ for all
$s\in\mathbb{R}$, thus dichotomy results for $\widetilde U$ will
immediately give similar dichotomy results for $U$.

It is easy to check (similarly to Lemma~\ref{L:closed}) that
$A(s;\lambda):\mathcal{X}^s \to  \mathcal{X}^s$ is closed and
densely defined with domain $\mathcal{Y}^s$ and that
$A(s;\lambda):\mathcal{Y}^s \to \mathcal{Y}^s$ is closed and densely
defined with domain $\tilde Z^s\times \tilde Z^s$.

The Fourier coefficients of $U = (u_1,u_2,u_3,u_4)^T$ satisfy the system
\begin{equation}\label{E:system at infinity}
   \begin{aligned}
      u_1' &= u_2, \\
      u_2' &= \frac{k^2}{s^2} u_1 - \frac{1}{s} u_2 + u_3, \\
      u_3' &= u_4, \\
      u_4' &= \lambda u_1 + \frac{k^2}{s^2} u_3 - \frac{1}{s} u_4,
   \end{aligned}
\end{equation}
of ODEs, where we omit the subscript $k$.
We denote by $\Phi_k(s,t)$ the evolution operator corresponding to the system \eqref{E:system at infinity} and
consider this evolution operator in either $\mathcal{X}^s$ or
$\mathcal{Y}^s$.

Now we can use the earlier dichotomy results to show the
existence of a uniform ($s$-dependent) exponential dichotomy for the system
\eqref{E:system at infinity} and hence for~\eqref{E:far field system}.

\begin{Lem}\label{L:far field dichotomy}
  There exists an $\epsilon_0>0$ such that for any $\epsilon\in
  (0,\epsilon_0)$ there exists a $K>0$ such that for any $\lambda\in
  (\lambda_0/2,2\lambda_0)$ there exists an $s$-dependent exponential dichotomy
  of~\eqref{E:far field system} on $J_+$ such that the evolution
  operators $\Phi^{s}_+(s,t;\lambda)$ and $\Phi^{cu}_+(s,t;\lambda)$
  solve \eqref{E:far field system} and
   \begin{equation}\label{E:far field dichotomy}
      \begin{array}{r l l}
        \|\Phi_+^s(s,t;\lambda)\|_{\mathcal L(\mathcal X^t,\mathcal
          X^s)}    =     \|\Phi_+^s(s,t;\lambda)\|_{\mathcal L(\mathcal Y^t,\mathcal Y^s)}
        &\le K e^{-(\lambda^{1/4}-\epsilon)(s-t)}, \quad &s\ge t\ge r_1, \\
        \|\Phi_+^{cu}(s,t;\lambda)\|_{\mathcal L(\mathcal X^t ,
          \mathcal X^s)}=
        \|\Phi_+^{cu}(s,t;\lambda)\|_{\mathcal L(\mathcal Y^t ,
          \mathcal Y^s)}
        &\le K e^{ \epsilon(t-s)}, \quad &t\ge s\ge r_1.
      \end{array}
   \end{equation}
The dichotomy is smooth in~$\lambda$ for $\lambda$ near $\lambda_0$.
\end{Lem}
\begin{proof}
  Since $\mathcal X^s= \overline {\oplus \mathcal X_k^s}$ and each
  $\mathcal X_k^t$ is mapped into $\mathcal X_k^s$ under the flow of
  \eqref{E:far field system}, we write for $U\in \mathcal X^t$
\begin{equation}\label{E:far field dich def}
   \begin{aligned}
      \Phi^s_+(s,t) U &:= \sum_{k\in \mathbb Z} \Phi_k^s(s,t)\widehat
      U_k e^{ik\cdot}, \quad &s\ge t\ge r_1, \\
      \Phi^{cu}_+(s,t) U &:= \sum_{k\in \mathbb Z}
      \Phi_k^{cu}(s,t)\widehat U_k e^{ik\cdot}, \quad &t\ge s\ge r_1.
   \end{aligned}
\end{equation}
Moreover, for each $k\in\mathbb{Z}$ the evolution operator
$\Phi_k(s,t)$ associated with~\eqref{E:system at infinity} can be
expressed in terms of $\phi_k(s,t)$ and $\psi_k(s,t)$.  Indeed, it can
be seen that
\begin{equation*}
   \Phi_k(s,t) = \frac{1}{2}\left(
   \begin{matrix}
      \phi_k(s,t) + \psi_k(s,t) & \frac{1}{\sqrt{\lambda}}(\phi_k(s,t) - \psi_k(s,t)) \\
      \sqrt{\lambda}(\phi_k(s,t) - \psi_k(s,t)) & \phi_k(s,t) + \psi_k(s,t)
   \end{matrix}
   \right).
\end{equation*}

Similarly,
\begin{equation*}
   \Phi^s_k(s,t) = \frac{1}{2}\left(
   \begin{matrix}
      \phi_k^s(s,t) & \frac{1}{\sqrt{\lambda}}\phi_k^s(s,t)  \\
      \sqrt{\lambda}\phi_k^s(s,t)  & \phi_k^s(s,t)
   \end{matrix}
   \right)
\end{equation*}
and
\begin{equation*}
\Phi_k^{cu}(s,t) = \frac{1}{2}\left(
   \begin{matrix}
      \phi_k^u(s,t) + \psi_k(s,t) & \frac{1}{\sqrt{\lambda}}(\phi_k^u(s,t) - \psi_k(s,t)) \\
      \sqrt{\lambda}(\phi_k^u(s,t) - \psi_k(s,t)) & \phi_k^u(s,t) + \psi_k(s,t)
   \end{matrix}
   \right).
\end{equation*}

Introducing the temporary notation $\widehat U_k=(\underline
u,\underline v)^T\in \mathcal X_k^t$, where $\underline u(u_1,u_2)^T$ and $\underline v=(u_3,u_4)^T$, we have
\begin{equation}\label{E:underline u}
  \|\Phi_k^s(s,t) \widehat U_k\|_{\mathcal X_k^s}^2 = \frac{1}{4}(1 + \lambda) \left\|\phi_k^s(s,t) \underline u + \frac{1}{\sqrt{\lambda}}\phi_k^s(s,t) \underline v\right\|_{\tilde X_k^s}^2.
\end{equation}
{Since we will take the supremum over all $\widehat U_k\in \mathcal X_k^t$ such that
$\|\widehat U_k\|_{\mathcal X_k^t}=1$, we may without loss of generality
assume that
$\underline v = \sqrt{\lambda}\underline u$, since all other choices will result in a smaller
value of the right hand side of \eqref{E:underline u}.}  For any such
$\underline u$ and $\underline v\in \tilde X_k^t$, the condition that
$\|\widehat U_k\|_{\mathcal X_k^t}^2=1$ implies that $\|\underline u\|_{\tilde
  X_k^t}^2 ( 1 + \lambda)=1$.  We therefore have
\begin{equation*}
   \begin{aligned}
     \sup_{\|\widehat U_k\|_{\mathcal X_k^t}^2 = 1} \|\Phi_k^s(s,t)
     \widehat U_k\|_{\mathcal X_k^s}^2 &= \sup_{\|\underline
       u\|_{\tilde X_k^t}^2 = 1/(1+\lambda)} \frac{1+
       \lambda}{2}\|\phi_k^s(s,t)
     \underline u\|_{\tilde X_k^s}^2 \\
     &= \frac{1}{2} \sup_{\|\underline u\|_{\tilde X_k^t}^2 = 1}
     \|\phi_k^s(s,t)\underline u\|_{\tilde X_k^s}^2,
   \end{aligned}
\end{equation*}
which shows that
\begin{equation*}
   \|\Phi_k^s(s,t)\|_{\mathcal L(\mathcal X_k^t,\mathcal X_k^s)} {=} \frac{1}{\sqrt{2}} \|\phi_k^s(s,t)\|_{\mathcal L(\tilde X_k^t,\tilde X_k^s)}.
\end{equation*}
Likewise,
\begin{equation*}
   \begin{aligned}
      \|\Phi_k^{cu}(s,t)\|_{\mathcal L(\mathcal X_k^t,\mathcal X_k^s)}
      \le \|\Phi_k^{c}(s,t)\|_{\mathcal L(\mathcal X_k^t,\mathcal X_k^s)} +
      \|\Phi_k^u(s,t)\|_{\mathcal L(\mathcal X_k^t,\mathcal X_k^s)},
   \end{aligned}
\end{equation*}
where
\begin{equation*}
  \Phi_k^{c}(s,t) =  \begin{pmatrix}
    \psi_k(s,t) & -\frac{\psi_k(s,t)}{\sqrt\lambda} \\
    -\sqrt\lambda\psi_k(s,t) & \psi_k(s,t)
  \end{pmatrix}
\mbox{~~and~~}  \Phi_k^u(s,t) = \Phi_k^{cu}(s,t) - \Phi_k^{c}(s,t).
\end{equation*}
As above,
\begin{equation*}
   \begin{aligned}
      \|\Phi^c_k(s,t)\|_{\mathcal L(\mathcal X_k^t,\mathcal X_k^s)} &{=}\frac{1}{\sqrt{2}} \|\psi_k(s,t)\|_{\mathcal L(\tilde X_k^t,\tilde X_k^s)}^2, \\
      \|\Phi^u_k(s,t)\|_{\mathcal L(\mathcal X_k^t,\mathcal X_k^s)} &{=}\frac{1}{\sqrt{2}}\|\phi_k^u(s,t)\|_{\mathcal L(\tilde X_k^t,\tilde X_k^s)}^2.
   \end{aligned}
\end{equation*}
From Lemmas~\ref{L:modified Bessel} and~\ref{L:Bessel} it follows that
\begin{equation}\label{E:far field k-dichotomy}
   \begin{array}{r l l}
      \|\Phi_k^s(s,t)\|_{\mathcal L(\mathcal X_k^t,\mathcal
        X_k^s)}&\le K e^{-(\lambda^{1/4}-\epsilon)(s-t)}, \quad &s\ge
      t\ge r_1, \\
      \|\Phi_k^{cu}(s,t)\|_{\mathcal L(\mathcal X_k^t,\mathcal
        X_k^s)}&\le K e^{ \epsilon(t-s)}, \quad &t \ge s \ge r_1.
   \end{array}
\end{equation}
By \eqref{E:far field dich def} and \eqref{E:far field k-dichotomy} we
see that for $s\ge t\ge r_1$
\begin{equation*}
   \begin{aligned}
     \|\Phi_+^s(s,t)\|^2_{\mathcal L(\mathcal X^t,\mathcal X^s)} &=     \sup_{\|U\|_{\mathcal X^t}=1} \|\Phi^s(s,t) U\|^2_{X^s}
     = \sup_{\|U\|_{\mathcal X^t}=1}\Bigl\|\sum_{k\in \mathbb Z}
     \Phi_k^s(s,t)\widehat U_k e^{ik\cdot}\Bigr\|^2_{\mathcal X^s} \\
     &= \sup_{\|U\|_{\mathcal X^t}=1}\sum_{k\in \mathbb Z}
     \left\|\Phi_k^s(s,t)\widehat U_k\right\|^2_{\mathcal X_k^s}
     \le \sup_{\|U\|_{\mathcal X^t}=1}\sum_{k\in \mathbb Z}\left\|\Phi_k^s(s,t)\right\|_{\mathcal L(\mathcal X^t,\mathcal X^s)}^2\|\widehat U_k\|^2_{\mathcal X_k^t} \\
     &\le \sup_{\|U\|_{\mathcal X^t}=1}\sum_{k\in \mathbb Z} K^2
     e^{-2(\lambda^{1/4}-\epsilon)(s-t)}\|\widehat U_k\|^2_{\mathcal X_k^t} =K^2
     e^{-2(\lambda^{1/4}-\epsilon)(s-t)}.
   \end{aligned}
\end{equation*}
A similar calculation shows that $\Phi^{cu}_+$ satisfies the second
equation of \eqref{E:far field dichotomy}. The estimates for $\tilde Y^s_k$ follow by the same estimates since it is just a matter of
 multiplying each side of the inequalities by a factor $(1+k^2)^2$.

Finally, the smoothness in $\lambda$ follows from the implicit
function theorem~\cite[Corollary~3.1.11]{mB77}. As $A(s;\lambda)$
depends linearly on $\lambda$, we get that
$\Phi_+^{s/cu}(s,t;\lambda)$ satisfies
\[
\frac{d}{ds} \Phi_+^{s/cu}(s,t;\lambda) =
\left[ A(s;\lambda_0) + (\lambda-\lambda_0)\, B_0  \right]\,
\Phi_+^{s/cu}(s,t;\lambda) ,
\]
where $B_0$ is defined by~\eqref{E:def B_0}.

{For
$\lambda\in(\lambda_0/2,2\lambda_0)$ and sufficiently close to
$\lambda_0$, the pair
$(\Phi_+^{s}(s,t;\lambda),\Phi_+^{cu}(s,t;\lambda))$ satisfies the
fixed point equation}
\begin{equation}\label{e:ips}
   \begin{aligned}
      \Phi_+^{s}(s,t;\lambda) &= \Phi_+^{s}(s,t;\lambda_0) {+\Phi_+^s(s,t;\lambda_0) \Phi_+^{cu}(s_1,t;\lambda)}
  +(\lambda-\lambda_0)\,
  \left[-  \int_{r_1}^t \Phi_+^{s}(s,\tau;\lambda_0) \, B_0 \,
    \Phi_+^{cu}(\tau,t;\lambda) \, d\tau \right. \\ \nonumber
    & \left.+\int_t^s \Phi_+^{s}(s,\tau;\lambda_0) \, B_0 \,
    \Phi_+^{s}(\tau,t;\lambda) \, d\tau
    - \int_s^{\infty}
    \Phi_+^{cu}(s,\tau;\lambda_0) \, B_0 \, \Phi_+^{s}(\tau,t;\lambda) \,
    d\tau   \right],\\ \nonumber
& s\geq t \geq r_1;\\ \nonumber
\Phi_+^{cu}(s,t;\lambda) &= \Phi_+^{cu}(s,t;\lambda_0) {- \Phi_+^s(s,t;\lambda_0) \Phi_+^{cu}(s_1,t;\lambda)}
  +(\lambda-\lambda_0)\,
  \left[\int_{r_1}^s
    \Phi_+^{s}(s,\tau;\lambda_0) \, B_0 \, \Phi_+^{cu}(\tau,t;\lambda) \,
    d\tau\right. \\ \nonumber
     & \left.-\int_s^t \Phi_+^{cu}(s,\tau;\lambda_0) \, B_0 \,
    \Phi_+^{cu}(\tau,t;\lambda) \, d\tau
   +  \int_t^{\infty} \Phi_+^{cu}(s,\tau;\lambda_0) \, B_0 \,
    \Phi_+^{s}(\tau,t;\lambda) \, d\tau   \right],\\ \nonumber
& t\geq s \geq r_1.
   \end{aligned}
\end{equation}
This fixed point equation is considered as a mapping on $\hat
X^s\times \hat X^{cu}$, where $\hat X^{s/cu}$ are defined as
\begin{eqnarray*}
  \hat X^s &=& \left\{\Phi^s(s,t)\in \mathcal{L}
    (\mathcal{X}^t,\mathcal{X}^s) \,;\, \|\Phi^s(s,t)\|_{\hat X^s}
    =\sup_{s\geq t\geq r_1} e^{((\lambda_0/2)^{1/4}-\epsilon)(s-t)}
    \|\Phi^s(s,t)\|_{\mathcal{L}(\mathcal{X}^t,\mathcal{X}^s)}<\infty
  \right\}\\
  \hat X^{cu} &=& \left\{\Phi^{cu}(s,t)\in \mathcal{L}
    (\mathcal{X}^t,\mathcal{X}^s) \,;\, \|\Phi^{cu}(s,t)\|_{\hat X^{cu}}
    =\sup_{t\geq s\geq r_1} e^{-2\epsilon(t-s)} \|\Phi^{cu}(s,t)\|_{\mathcal{L}
      (\mathcal{X}^t,\mathcal{X}^s)}<\infty  \right\}.
\end{eqnarray*}
With the estimates derived before on $\Phi_+^s(s,t;\lambda)$ and
$\Phi_+^{cu}(s,t;\lambda)$, it is easy to check that the right-hand
side is well-defined in this space. With the implicit function
theorem, it follows immediately that the mapping is smooth for
$\lambda$ near $\lambda_0$.
\end{proof}
\end{section}


\begin{section}{The unperturbed adjoint system}\label{S:adjoint}
The dual space of $X$ is $X' = H^{-2}\times H^{-1}\times H^{-1} \times
L^2$ (using the $L^2$ pairing).  For the space $X'$, we make the
decomposition
\begin{equation*}
   X' = \overline{\oplus_{k\in \mathbb Z}X_k'},
\end{equation*}
where $X_k'$ are $4$-dimensional pairwise orthogonal subspaces $\span
\{\left((a,b,c,d)e^{ik\cdot}\right);\; a, b, c, d\in \mathbb C\}
\subset X'$.  For $W\in X_k'$ and $U\in X_k$, we have the pairing
\begin{equation}\label{E:dual pairing}
  \langle W, U\rangle :=  \overline w_1 u_1 + \overline w_2 u_2 + \overline w_3 u_3 + \overline w_4 u_4,
\end{equation}
where the bar denotes the complex conjugate. This means that we may
use the standard inner product on $\mathbb C^4$ when computing the
adjoint equation.

Similarly, the dual space of $\mathcal{X}:=H^{1}\times L^{2}\times
H^{1} \times L^2$ is $\mathcal{X}' = H^{-1}\times L^{2}\times H^{-1}
\times L^2$ (using the $L^2$ pairing) and we can make the same
decomposition as above, i.e.,
\begin{equation*}
   \mathcal{X}' = \overline{\oplus_{k\in \mathbb Z}\mathcal{X}_k'},
\end{equation*}
where $\mathcal{X}_k'$ are the same $4$-dimensional subspaces as
above but are now regarded as subspaces of $\mathcal{X}'$. For $W\in
\mathcal{X}_k'$ and $U\in \mathcal{X}_k$, the pairing is
again as in~\eqref{E:dual pairing}.

At the end of Section~\ref{S:r_small}, we have investigated the
solutions of the unperturbed linear system
$U'=A(s;\lambda_0,0)U$ on $J_-$. The adjoint unperturbed system is
\begin{equation}\label{E:adjoint}
   W' = -A(s;\lambda_0,0)^* W = -(A_-^*+B(s;\lambda_0,0)^*)\,W.
\end{equation}
Just as in the case of the unperturbed linear system itself, expanding
$W$ in a Fourier series shows that the Fourier spaces $X_k'$ are
invariant under the flow of the adjoint system~\eqref{E:adjoint}, and
the Fourier coefficients satisfy the adjoint equation
of~\eqref{E:Fourier_unperturbed}, i.e.,
\begin{equation}\label{E:Fourier_unperturbed_adjoint}
  \wW_k'(s) = -\left[\wA_-(k)^* + B(s;\lambda,0)^*\right]\,\wW_k(s).
\end{equation}
It is well known and straightforward to check that the pairing of a
solution of a linear system with a solution of its adjoint is
constant. For our systems, this means that any two solutions
$\wU_k(s)$ of~\eqref{E:Fourier_unperturbed} and $\wW_k(s)$
of~\eqref{E:Fourier_unperturbed_adjoint} satisfy
\begin{equation}\label{E:dual}
  \frac{d}{ds}\langle \wW_{k}(s),\wU_{k}(s)\rangle = 0,
\mbox{ and thus } \langle \wW_{k}(s),\wU_{k}(s)\rangle = \langle
\wW_{k}(s_1),\wU_{k}(s_1)\rangle \mbox{ for any } s\in \mathbb R.
\end{equation}

From \cite[p.\,56]{SS01} it follows that if a finite-dimensional
linear system has an exponential dichotomy on an interval~$J$ with
constants $K$, $\kappa^s$ and $\kappa^u$, then the adjoint system has
an exponential dichotomy on~$J$ with the dichotomy constants $K$,
$-\kappa^u$ and $-\kappa^s$. Furthermore, if we denote evolution operators
corresponding to the exponential dichotomy of the adjoint
system~\eqref{E:Fourier_unperturbed_adjoint} by
$\widehat\Phi^{s}_{k\pm}(s,t)$ and $\widehat\Phi^{u}_k(s,t)$, respectively, then
$\widehat\Phi_k^s(s,t)=\Phi^{u}_k(t,s)^*$ for $t\leq s$ with $s,t\in J$ and
$\widehat\Phi_k^{u}(s,t)=\Phi^{s}_k(t,s)^*$ for $s\leq t$ with $s,t\in J$.

 On $J_-$, the dichotomy constant $K$
 in~\eqref{E:unperturbed_dichotomy_k_-} is independent of $k$, and so we
 immediately get the following estimates about the solutions of the
 adjoint system in the Fourier spaces $X'_k$ with norm
 \begin{equation*}
    \|\wW_k\|^2_{X'_k}:=\|\wW_k {e^{ik\cdot}}\|^2_{X'}= \frac{|w_1|^2}{(k^2+1)^2} +
    \frac{|w_2|^2}{k^2+1} + \frac{|w_3|^2}{k^2+1} + |w_4|^2:
 \end{equation*}
\begin{Lem}\label{L:adjoint_upper_bound_J-}
  There exists a $K>0$ such that for every $k\in
  \mathbb{Z}\backslash\{0\}$ and $\wW_k\in X'_k$
\begin{equation*}
\begin{aligned}
  \|\widehat\Phi^s_{k}(s,t) \wW_k \|_{X'_{{k}}} &\leq K\, e^{-|k|(s-t)}\,
  \|\wW_k\|_{X'_{{k}}}, &\quad t\leq s\leq s_1,\\
  \|\widehat\Phi^u_{k}(s,t) \wW_k \|_{X'_{{k}}} &\leq K\, e^{|k|(s-t)}\,
  \|\wW_k\|_{X'_{{k}}}, &\quad s\leq t\leq s_1 .
\end{aligned}
\end{equation*}
Furthermore, for any solution $\wW_k(s)$ with $\wW_k(s_1)\in
\Ran(\widehat\Phi_{k{-}}^u(s_1,s_1))$, we have
\[
|w_{k}(s)|\leq Ke^{|k|(s-s_1)}\|\wW_k(s_1)\|_{X'}
\]
for all $s\leq s_1$, where $w_{k}(s)$ denotes the fourth component of $\wW_k(s)$.
\end{Lem}

Similar arguments give the dichotomy of the adjoint system
on~$J_+$. From~\eqref{E:far field k-dichotomy}, it follows that the
solutions of the linearised system in the Fourier spaces
$\Phi^s_k(s,t)$ and $\Phi^{cu}_k(s,t)$ have an exponential dichotomy
with a uniform constant~$K$. The dual space of $\mathcal{X}^s$ is
denoted by $(\mathcal{X}^s)'$, and for $s$ fixed, this space is equivalent to
$H^{-1}\times L^2\times H^{-1}\times L^2$. The dual Fourier space is
denoted by $(\mathcal{X}_k^s)'$, and its norm is
\[
\|\wW_k\|^2_{(\mathcal{X}_k^s)'}:=\|\wW_k {e^{ik\cdot}}\|^2_{(\mathcal{X}^s)'}
=\frac{s^2}{k^2+s^2} |w_1|^2 +|w_2|^2 + \frac{s^2}{k^2+s^2}|w_3|^2+
|w_4|^2.
\]
On $J_+$, we have the following estimates:
\begin{Lem}\label{L:adjoint_upper_bound_J+}
  For every $\epsilon>0$, there exists a $K>0$ such that for every $k\in
  \mathbb{Z}\backslash\{0\}$ and $\wW_k\in \mathcal{X}'_k$
\begin{equation*}
\begin{aligned}
  \|\widehat\Phi^{cs}_k(s,t) \wW_k \|_{(\mathcal{X}^s)'}
  &\leq K\, e^{\epsilon(s-t)}\,
  \|\wW_k \|_{(\mathcal{X}^t)'}, &\quad s_1\leq t\leq s,\\
  \|\widehat\Phi^{u}_k(s,t) \wW_k \|_{(\mathcal{X}^s)'}&\leq K\,
  e^{(\lambda_0^{1/4}-\epsilon)(s-t)}\,
  \|\wW_k\|_{(\mathcal{X}^t)'}, &\quad   s_1\leq s\leq t.
\end{aligned}
\end{equation*}
Moreover, for any solution $\wW_k(s)$ with $\wW_k(s_1)\in
\Ran(\widehat\Phi_k^{cs}(s_1,s_1))$, we have
\[
|w_{k}(s)|\leq
Ke^{\epsilon(s-s_1)}\|\wW_k(s_1)\|_{(\mathcal{X}^{s_1})'}
\]
for all $s\geq s_1$, where $w_{k}(s)$ denotes the fourth component of $\wW_k(s)$.
\end{Lem}

Next we look at the adjoint system associated with $k=0$.  In
Lemma~\ref{L:centre}, we have seen that the solutions space of
the linear system at $k=0$ is spanned by $U_{0,j}$, $j=1,\ldots,4$.
Now let $Z_{0,j}(s)$ be solutions of the adjoint
system~\eqref{E:Fourier_unperturbed_adjoint} with $k=0$ such that
$\{Z_{0,l}(s_1)\}_{l=1,\ldots,4}$ is a dual basis of
$\{U_{0,j}(s_1)\}_{j=1,\ldots,4}$ (i.e.,
$\langle Z_{0,l}(s_1),U_{0,j}(s_1)\rangle=\delta_{lj}$). With~\eqref{E:dual},
this implies $\delta_{lj}=\langle Z_{0,l}(s),U_{0,j}(s)\rangle$ for
any $s\leq s_1$.
The unbounded solutions $U_{0,2}$ and $U_{0,4}$ are not unique, but
they can be chosen such that
$Z_{0,2}$ and $Z_{0,4}$ are bounded on  $J_-$,
whereas $Z_{0,1}$ and $Z_{0,3}$ grow algebraically as $s\to -\infty$. {A
convenient choice for $Z_{0,2}$ and $Z_{0,4}$ is
\begin{equation}\label{E:adjoint 0}
   Z_{0,j}(z) := \frac{r(s)}{r'(s)}U_{0,j-1}^\perp(s),\qquad j=2, 4,
\end{equation}
where $U^\perp = (-u_4,u_3,-u_2,u_1)$ if $U=(u_1,u_2,u_3,u_4)^T$. It is easy to
check that $Z_{0,j}^{\perp}(s)$ are solutions of the adjoint system \eqref{E:adjoint}
using that $U_{0,j-1}(s)$ are solutions of the original system \eqref{E:perturbed system}.}

As will be shown below, a similar exponential dichotomy on $J_-$ as in
Lemma~\ref{L:adjoint_upper_bound_J-} also holds with norms in $\mathcal{X}'$.
\begin{Lem}\label{L:adjoint_upper_bound_J-2}
  There exists a $K>0$ such that, for every $k\in
  \mathbb{Z}\backslash\{0\}$ and $\wW_k\in (\mathcal{X}_k)'$, we have
\begin{equation*}
\begin{aligned}
  \|\widehat\Phi^s_k(s,t) \wW_k \|_{\mathcal{X}'} &\leq K\,
  e^{-|k|(s-t)}\, \|\wW_k \|_{\mathcal{X}'}, &\quad t\leq s\leq s_1,\\
  \|\widehat\Phi^u_k(s,t) \wW_k \|_{\mathcal{X}'} &\leq K\,
  e^{|k|(s-t)}\, \|\wW_k \|_{\mathcal{X}'}, &\quad s\leq t\leq s_1 .
\end{aligned}
\end{equation*}
For any solution $\wW_k(s)$ with $\wW_k(s_1)\in
\Ran(\widehat\Phi_k^u(s_1,s_1))$ and with $w_{k}(s)$ denoting the fourth component
of $\wW_k(s)$, we have $|w_{k}(s)|\leq
Ke^{|k|(s-s_1)}\|\wW_k(s_1)\|_{\mathcal{X}'}$ for all $s\leq s_1$.
\end{Lem}
\begin{proof}
  First we will show that the linear
  system~\eqref{E:Fourier_unperturbed} has an exponential dichotomy in
  $\mathcal{X}_k$. The proof is very similar to the one in
  section~\ref{S:r_small} with a slightly modified matrix $M_k$. Define the
  matrix $\wM_k$ whose columns consist of eigenvectors of $A_-(k)$ that are scaled different to those in $M_k$:
  \begin{equation*}
\wM_k = \left(
   \begin{matrix}
      -1/|k| & 0 & 1/|k| & 0 \\
      1 & 0 & 1 & 0 \\
      0 & -1/|k| & 0 & 1/|k| \\
      0 & 1 & 0 & 1
   \end{matrix}
   \right)    .
  \end{equation*}
As $\wM_k$ consists of eigenvectors of   $A_-(k)$, it follows
immediately that $A_-(k) = \wM_k D_k \wM_k^{-1}$. It is also
straightforward to verify that $\wM_k$ is a homeomorphism between
$\mathbb{C}^4$ and $\mathcal{X}_k$ with
$\|\wM_k\|_{\mathcal{L}(\mathbb{C}^4,\mathcal{X}_k)} \to \sqrt2$ as
$|k|\to\infty$.

Using the same ideas as in the proof
of~\eqref{E:unperturbed_dichotomy_k_-}, with $M_k$ replaced by $\wM_k$
and exploiting the observation that
\begin{equation*}
   |\wM_k^{-1}B(\tau/|k|;\lambda_0,0)\wM_k|\leq \frac {e^{2\tau/|k|}}{2|k|} \sup_{s\leq s_1} \{1,
   |\lambda-\tilde\theta(s)|\},
\end{equation*}
we find that
\begin{equation}\label{E:unperturbed_dichotomy_k_-2}
\begin{aligned}
  \|\Phi^s_k(s,t) \wU_k \|_{\mathcal X} &\leq K\, e^{-|k|(s-t)}\,
  \|\wU_k \|_{\mathcal X}, &\quad t\leq s\leq s_1,\\
  \|\Phi^u_k(s,t) \wU_k \|_{\mathcal X} &\leq K\, e^{|k|(s-t)}\,
  \|\wU_k \|_{\mathcal X}, &\quad s\leq t\leq s_1
\end{aligned}
\end{equation}
for any $\wU_k\in \mathcal{X}_k$ for some constant $K$ that is independent of $k$.
As $\mathcal{X}_k$ is finite-dimensional, this immediately implies the
estimates of the Lemma.
\end{proof}

{Finally we will show that, for large values of $k$, the solutions of
\eqref{E:Fourier_unperturbed_adjoint} are close to the solutions of
the asymptotic
system $\wW_k'(s)=-A_-(k)^* \wW_k(s)$.}  Recall that we
denote the spectral projection onto the eigenspace of $A_-(k)$ associated with the positive eigenvalue $|k|$ by $P_k^u$ and the complementary projection onto the eigenspace of
$A_-(k)$ associated with the negative eigenvalue $-|k|$ by $P_k^s$.

\begin{Lem}\label{L:adjoint_lower_bound}
  For every $\epsilon>0$, there exists an ${N}\in\mathbb{N}$ and a $\delta>0$ such
  that, for every $|k|>N$,
  we have
\begin{equation}\label{E:sol-close}
   \begin{aligned}
    \|\widehat\Phi^{u}_k(s,s_1)-
    e^{|k|(s-s_1)}(P_k^{s})^*\|_{\mathcal{L}(X'_k)} &\leq \epsilon\,
    e^{|k|(s-s_1)},&\mbox{for all } s_1-\delta\leq s\leq s_1,\\
    \|\widehat\Phi^{u}_k(s,s_1)- e^{|k|(s-s_1)}(P_k^{s})^*
    \|_{\mathcal{L}(\mathcal{X}'_k)} &\leq \epsilon\,
    e^{|k|(s-s_1)},&\mbox{for all } s\leq s_1.
  \end{aligned}
\end{equation}
  Thus, for $|k|>{N}$ and $\wW_k(s_1)\in \Ran(\widehat\Phi_k^{u}(s_1,s_1))$,
  we have
  \begin{equation*}
     \begin{aligned}
        \|\wW_k(s_1)- (P_k^{s})^* \wW_k(s_1)\|_{X'_k} &\leq \epsilon\|\wW_k(s_1)\|_{X'_k}, \\
        \|\wW_k(s_1)- (P_k^{s})^* \wW_k(s_1)\|_{\mathcal{X_k}'} &\leq   \epsilon\|\wW_k(s_1)\|_{\mathcal{X_k}'}.
     \end{aligned}
  \end{equation*}
\end{Lem}
\begin{proof}
The {evolution operator} $\widehat\Phi^{u}_k(s,s_1)$ satisfies
\begin{equation*}
\begin{aligned}
\widehat\Phi^{u}_k(s,s_1) ={}& e^{|k|(s-s_1)}(P_k^{s})^*
  - \int_{-\infty}^s e^{-|k|(s-t)}(P_k^u)^* B(t;\lambda_0,0)^*
  \widehat\Phi^{u}_k(t,s_1)\, dt \\
&{}+ \int_s^{s_1}   e^{|k|(s-t)}(P_k^s)^* B(t;\lambda_0,0)^*
  \widehat\Phi^{u}_k(t,s_1)\, dt.
\end{aligned}
\end{equation*}
From its definition~\eqref{E:B in J-}, we immediately see that
$\|B(t;\lambda_0,0)^*\|_{\mathcal{L}(X'_k)} \leq Ce^{2t}$ and
$\|B(t;\lambda_0,0)^*\|_{\mathcal{L}(\mathcal{X}'_k)}
\leq \frac{C}{\sqrt{k^2+1}}\,e^{2t}$ for some constant $C$,
independent of $k$. Hence, with the dichotomy estimates from
Lemma~\ref{L:adjoint_upper_bound_J-}, we get that for $s\leq s_1$ and
$k\in\mathbb{Z}\backslash\{0\}$
\begin{equation*}
\begin{aligned}
\|\widehat\Phi^{u}_k(s,s_1) -
e^{|k|(s-s_1)}(P_k^{s})^*\|_{\mathcal{L}(X'_k)}  &\leq
CK \int_{-\infty}^s e^{-|k|(s-t)}e^{2t}e^{|k|(t-s_1)}dt \\
&\qquad{}+CK \int_s^{s_1}  e^{|k|(s-t)}e^{2t}e^{|k|(t-s_1)}dt\\
&\leq \frac{CK}{2} e^{|k|(s-s_1)}\,\left(\frac{e^{2s_1}}{1+|k|} +
  e^{2s_1}-e^{2s} \right).
\end{aligned}
\end{equation*}
{It is now easy to see that we can choose $N>0$ large enough and $\delta>0$ small enough such that the
first inequality of \eqref{E:sol-close} is satisfied.}

With the dichotomy estimates from
Lemma~\ref{L:adjoint_upper_bound_J-2}, we get
\begin{equation*}
\begin{aligned}
\|\widehat\Phi^{u}_k(s,s_1) -
e^{|k|(s-s_1)}(P_k^{s})^*\|_{\mathcal{L}(\mathcal{X}'_k)}  &\leq
\frac{C}{\sqrt{k^2+1}}\,K \int_{-\infty}^s e^{-|k|(s-t)}e^{2t}e^{|k|(t-s_1)}dt \\
&\qquad{}+\frac{C}{\sqrt{k^2+1}}\,K \int_s^{s_1}
e^{|k|(s-t)}e^{2t}e^{|k|(t-s_1)}dt\\
&\leq \frac{CK}{2\sqrt{k^2+1}} e^{|k|(s-s_1)}\,\left(\frac{e^{2s_1}}{(1+|k|)} +
  e^{2s_1} \right) \\
  &{\leq \frac{CK e^{2s_1}}{\sqrt{k^2+1}} e^{|k|(s-s_1)}}.
\end{aligned}
\end{equation*}
{It follows that also the second inequality of \eqref{E:sol-close} is valid when $N$ is sufficiently large.}
\end{proof}
\end{section}


\begin{section}{{Matching the core and far field solutions}}\label{S:Lyapunov-Schmidt}
  In the next lemma we show that $u$ is an eigenfunction that belongs to an embedded eigenvalue $\lambda$ of $\mathcal L+ \rho$ if and only if $U$ is a
  solution of~\eqref{E:perturbed system} {such that $U(s_1)\in \Ran P_+^s(s_1;\lambda,\tilde\rho)\cap (\Ran P_-^u(s_1;\lambda,\tilde\rho)\oplus P_-^{cb}(s_1;\lambda,\tilde\rho))$.}
\begin{Lem}\label{L:system-eigenfunction}
  Let $u$ be an $L^2$ solution of \eqref{E:eigenvalue_equation}. Then
  the corresponding solution $U(s)$ of the system \eqref{E:perturbed
    system} is bounded in $X$ as $s\to -\infty$
  and decays exponentially with rate $\lambda^{1/4} -\epsilon$ as $s\to +\infty$
  in the sense that for any $\epsilon\in (0,\lambda^{1/4})$ there exists a constant $K>0$  such
  that
   \begin{equation}\label{E:exponential decay}
      {\|U(s)\|_{\mathcal X^s} \le K e^{-(\lambda^{1/4}-\epsilon) s}}
   \end{equation}
   for every $s\ge s_1$.
   Conversely, if $U$ is a weak solution of \eqref{E:perturbed system}
   such that $\|U(s)\|_{X}$ is bounded as $s\in J_-$
   and such that $\|U(s)\|_{\mathcal X^s}$ decays exponentially as
   $s\to +\infty$ {(with any decay rate)},
   then it corresponds to an $H^4$ solution $u$ of
   \eqref{E:eigenvalue_equation}.
\end{Lem}
\begin{proof}
  If $u$ is an eigenfunction of \eqref{E:perturbed system}, it belongs
  to $H^4(\mathbb R^2)$. By Lemma \ref{L:equation-system}, the
  associated solution $U(s)$ of \eqref{E:perturbed system} is bounded
  in $X$ as $s\in J_-$.  In $J_+$, for $s\geq s_3$, the
  system~\eqref{E:perturbed system} reduces to~\eqref{E:far field
    system} and the decaying solutions of this system are series
  formed by Bessel functions $K_k$, $J_k$ and $Y_k$. The decay of
  $J_k$ and $Y_k$ is asymptotic to $1/\sqrt{s}$ as $s\to+\infty$, and
  so these solutions do not give rise to $L^2$ solutions of
  \eqref{E:eigenvalue_equation}. It now follows from Lemma~\ref{L:far
    field dichotomy} that $U(s)$ decays exponentially as $s\to
  +\infty$, in the sense that \eqref{E:exponential decay} holds.

  Assume that $U$ is a bounded weak solution of \eqref{E:perturbed
    system} which decays exponentially as $s\in J_+$.  We need to show
  that the first component of $U$ which we denote by $u$ belongs to
  $H^4(\mathbb R^2)$ when regarded as a function of the two variables
  $(r,\varphi)$ in radial coordinates. As $U$ is a weak solution
  of~\eqref{E:perturbed system}, by Definition~\ref{D:weak sol}, $U\in
  L^2_{loc}(J;Y)\cap H^1_{loc}(J;X)$. Also, $\|U\|_X$ is bounded on
  $J_-$, and and hence $U\in L^\infty_{loc}(\mathbb{R}_-;X)$.  From
  Lemma~\ref{L:equation-system} we know that $u\in
  H^4_{loc}(\mathbb{R}^2)$, so we only need to worry about the decay
  properties of $u$ as $r\to \infty$ (i.e.\ as $s\to\infty$).  From
  Lemma~\ref{L:far field dichotomy} and the definition of $\mathcal
  X^s$ it follows that for any $0<\epsilon<\epsilon_0$ and $s\ge s_1$
  \begin{equation*}
     \left\|u(s)\right\|_{L^2(S^1)}\le
     K e^{-(\lambda^{1/4}-\epsilon) (s-s_1)}\|U(s_1)\|_{\mathcal X^{s_1}}
   \end{equation*}
   and so it is clear that $u\in L^2(\mathbb R^2)$. From
   \eqref{E:eigenvalue_equation} it then follows that $u\in
   H^4(\mathbb R^2)$.
\end{proof}

Recall that $u_*$ is the radially symmetric eigenfunction associated with the embedded
eigenvalue~$\lambda_0$ when $\tilde\rho=0$.  Let $U_*$ be the associated solution of \eqref{E:perturbed
  system} with $\tilde \rho=0$ and $\lambda=\lambda_0$, defined for
$s\in \mathbb R$, i.e.\ $U_* = (u_*,u_*',\Delta u_*, (\Delta u_*)')^T$.
Define $\mathcal X := H^1\times L^2\times H^1\times L^2$ and recall
that $X=H^2\times H^1\times H^1\times L^2$.  Let
\begin{equation*}
   \begin{aligned}
      E^s_+ :&= \{U\in \mathcal X;\; P^s_+(s_1;\lambda_0,0)U = U\}, \\
      E^u_- :&= \{U\in X;\; P^u_-(s_1;\lambda_0,0)U = U\}, \\
      E^{cb}_- :&= \span\{U_{0,1}(s_1),U_{0,3}(s_1)\}\subset X, \\
   \end{aligned}
\end{equation*}
where $U_{0,j}$ are defined in Lemma~\ref{L:centre}. Roughly
speaking, $E^s_+$ and $E^{u}_-$ consist of the initial values at $s_1$
of solutions of \eqref{E:perturbed system} with $\tilde \rho=0$ and
$\lambda=\lambda_0$ which decay exponentially as $s\to \infty$ and as
$s\to -\infty$, respectively, and $E^{u}_-\oplus E_-^{cb}$ consists of
the bounded solutions on $J_-$. Note that the norm of $\mathcal X$ is used for
$E_+^s$, while the norm of $X$ is used for $E_-^u$ and $E_-^{cb}$.  We
have $E^s_+ \cap (E^{u}_-\oplus E^{cb}_-) = \span\{U_*(s_1)\}$ since
$\lambda_0$ is an eigenvalue of $\mathcal L$ with multiplicity $1$.

Next, we introduce a new Hilbert space $\overline X$ such that
$X\subset \overline X\subset \mathcal X$, and special solutions
$V_{k,j}$, $k\in \mathbb Z$, $j=1,\dots, 4$, such that
$\{V_{k,j}(s_1)\}$ is a basis for $\overline X$.

We have seen that the unperturbed system \eqref{E:perturbed system}
decouples when $\tilde \rho=0$ so that the subspaces $X_k$ and
$\mathcal{X}_k$ are invariant under the flow of \eqref{E:perturbed
  system} with $\tilde\rho=0$.  For $k=0$, we pick
$V_{0,1}(s_1) := U_*(s_1)\in E_-^{cb} \cap E_+^s$.  Note that there are no other
solutions in $X_0$ which decay exponentially as $s\to \infty$.  We
pick $V_{0,2}(s_1)\in E_-^{cb}$ such that $V_{0,2}(s_1)\notin E_+^s$
and $E_-^{cb}= \span\{V_{0,1}(s_1),V_{0,2}(s_1)\}$ (thus
$\span\{V_{0,1}(s_1),V_{0,2}(s_1)\}=\span\{U_{0,1}(s_1),U_{0,3}(s_1)\}$).
{Next, we will choose $V_{0,3}(s_1)$ and $V_{0,4}(s_1)\in X_0$ such that
they belong to the span of $U_{0,2}(s_1)$ and $U_{0,4}(s_2)$.
In order to do this, we introduce a dual basis $W_{0,j}(s_1)$ of $V_{0,j}(s_1)$ and choose
$W_{0,3}(s_1):=  U_*^\perp(s_1)$, where $U_*^\perp(s_1) = (-u_4(s_1),u_3(s_1),-u_2(s_1),u_1(s_1))$
{and} 
$u_j(s_1)$ are the components of $U_*(s_1)$, while $W_{0,4}(s_1)$ is any other vector so that
$\span\{W_{0,3}(s_1),W_{0,4}(s_1)\}=\span \{Z_{0,2}(s_1),Z_{0,4}(s_1)\}$. The remaining
vectors $W_{0,1}(s_1)$, $W_{0,2}(s_1)$, $V_{0,3}(s_1)$ and $V_{0,4}(s_1)$ are determined
by the conditions that $\{W_{0,j}(s_1);\; j=1,\dots, 4\}$ and $\{V_{0,j}(s_1);\; j=1,\dots, 4\}$
are dual bases:
\begin{equation*}
   \langle W_{0,j}(s_1),V_{0,l}(s_1)\rangle = \delta_{jl}, \qquad l=1,\dots, 4.
\end{equation*}}

 We use the notation
$E_-^{ca}:= \span\{V_{0,3}(s_1),V_{0,4}(s_1)\}$. We normalise the vectors
such that $\|V_{0,j}(s_1)\|_{X}=1$ for $j=1,\ldots,4$ and note that, for
$V\in X_0$, we have $\|V\|_X=\|V\|_{\mathcal{X}}$.
   Define
$W_{0,j}(s)$ so that it satisfies the adjoint system~\eqref{E:adjoint}
(and hence~\eqref{E:Fourier_unperturbed_adjoint} with $k=0$) and passes
through $W_{0,j}(s_1)$ for $s=s_1$. From~\eqref{E:dual} and the
relation above, it follows immediately that $\langle
W_{0,j}(s),V_{0,l}(s)\rangle = \delta_{jl}$ for all $s\leq s_1$. Furthermore, from
\eqref{E:adjoint 0}, we conclude that $W_{0,3}(s)$ and $W_{0,4}(s)$ are
bounded solutions of the adjoint system on $J_-$.

Next we consider $k\neq 0$. The spaces $X_k$ and $\mathcal{X}_k$ are
four-dimensional, and $E_+^s\cap X_k$ and $E_+^s\cap \mathcal{X}_k$ are one-dimensional,
 $E_-^u \cap X_k$ and $E_-^u\cap \mathcal{X}_k$ are two-dimensional, and
$E_-^u\cap E^s_+\cap X_k=\{0\}=E_-^u\cap E^s_+\cap
\mathcal{X}_k$, since the multiplicity of the eigenfunction $U_*$
is~$1$. Using this, we define base vectors in $X_k$ and
$\mathcal{X}_k$ as follows:  For $k\ne 0$ we pick $V_{k,1}(s_1)\in
E_+^{s}$ (and hence $V_{k,1}(s_1)\notin E_-^{u}$). We also pick
$V_{k,2}(s_1)$ and $V_{k,3}(s_1)$ so that they belong to $E_-^u$ (and
hence do not belong to $E_+^{s}$). Thus
$\{V_{k,1}(s_1), V_{k,2}(s_1), V_{k,3}(s_1)\}$ span a three-dimensional subspace in the four-dimensional spaces $X_k$ and
$\mathcal X_k$. We normalize the solutions
$V_{k,j}$ such that for $k\in \mathbb{Z}\backslash\{0\}$:
$\|V_{k,1}(s_1)\|_{\mathcal X}=1$ and $\|V_{k,j}(s_1)\|_{X}=1$ for
$j=2,3$.
Hence there exists a unique ({up to multiplication by a unimodular constant}) vector~$W_{k,4}\in
\mathcal X_k'$ such that
\[
\langle W_{k,4}\,,\,V_{k,j}\rangle =0  \mbox{ for } j=1,2,3
\mbox { and } \|W_{k,4}(s_1)\|_{\mathcal X_k'} =1.
\]
Then {$\langle W_{k,4}(s_1),V\rangle=0$ for
$V\in \Ran(\Phi_{k,+}^s(s_1,s_1))+\Ran(\Phi_{k,-}^u(s_1,s_1))$} and hence
$W_{k,4}(s_1)\in\Ran(\Phi_{k,+}^{cu}(s_1,s_1)^*) \cap
\Ran(\Phi_{k,-}^{s}(s_1,s_1)^*)=\Ran(\widehat\Phi_{k,+}^{cs}(s_1,s_1)) \cap
\Ran(\widehat\Phi_{k,-}^{u}(s_1,s_1))$.
We take the
one remaining solution in $X_k$ and $\mathcal{X}_k$ such that
$\langle W_{k,4}(s_1)\,,\,V_{k,4}(s_1)\rangle=1$ and
$\|V_{k,4}(s_1)\|_{\mathcal{X}} = 1$.
Then $V_{k,4}(s_1)\notin E_-^{u} \cup E_+^{s}$ as {$\langle W_{k,4}(s_1), V\rangle=0$ for  $V\in E_-^{u} + E_+^{s}$}.

Let $\overline X$ be defined by
\begin{equation*}
   \overline X := \biggl\{U = \sum_{
   \begin{smallmatrix}
      k\in \mathbb Z \\
      j=1,\dots,4
    \end{smallmatrix}} a_{k,j} V_{k,j}(s_1) \in \mathcal X;\;
  \|U\|_{\overline X}^2 := \biggl\| \sum_{\begin{smallmatrix}
       k\in \mathbb Z \\
       j=1, 4
   \end{smallmatrix}} a_{k,j}
  V_{k,j}(s_1) \biggr\|_{\mathcal X}^2 + \biggl\| \sum_{
   \begin{smallmatrix}
       k\in \mathbb Z \\
       j=2, 3
   \end{smallmatrix}}
   a_{k,j} V_{k,j}(s_1)\biggr\|_{X}^2 <\infty\biggr\}.
\end{equation*}
Note that $\overline X$ is the direct sum of two Hilbert spaces, which
are closed subspaces of $\mathcal X$ and $X$, respectively.  It
follows that $\overline X$ is a Hilbert space.

Note that $E_+^s$ and $E_-^u$ are both closed subspaces of $\overline
X$: Indeed,
{\begin{equation*}
   \begin{aligned}
      E_+^s &=\cl_{\overline X} \span\{V_{k,1}(s_1)\}_{k\in \mathbb Z}= \cl_{\mathcal X} \span \{V_{k,1}(s_1)\}_{k\in \mathbb Z} \\
      E_-^{u} &= \cl_{\overline X} \span \{V_{k,2}(s_1),V_{k,3}(s_1)\}_{k\in\mathbb Z\setminus \{0\}} =\cl_{X} \span
    \{V_{k,2}(s_1),V_{k,3}(s_1)\}_{k\in\mathbb Z\setminus \{0\}}.
   \end{aligned}
\end{equation*}}
It is clear that $E_-^{cb}$ is a closed subspace of
$\overline X$ since it is finite-dimensional.

Define $\iota:E_+^s\times (E_-^{u} \oplus E_-^{cb}) \times \mathbb
R\times \tilde {\mathcal R}\to \overline{X}$ by
\begin{equation}\label{E:iota-def}
   \begin{aligned}
      \iota(U_0^s,U_0^{u} + U_0^{cb};\lambda,\tilde \rho) :&= P_-^{u}(s_1;\lambda,\tilde \rho) U_0^{u} +
      P_-^{cb}(s_1;\lambda,\tilde \rho) U_0^{cb} -
      P^s_+(s_1;\lambda,\tilde \rho) U_0^s \\ &=       P_-^{u}(s_1;\lambda,\tilde \rho) U_0^{u} +
      P_-^{cb}(s_1;\lambda,\tilde \rho) U_0^{cb} -
      P^s_+(s_1;\lambda,0) U_0^s,
   \end{aligned}
\end{equation}
where we recall the definition (\ref{e:ipc}) and note that the last
equality holds since $\tilde \rho(s)=0$ for $s\ge s_1$. We will see
that the range of $\iota$ is a subspace of $\overline X$, and that
$\iota$ is smooth into this space. For this we need a more explicit
formula for $\iota$.  By \eqref{e:ipc}, \eqref{e:ipu} and
\eqref{e:ips} evaluated at $s=t=s_1$, we have (see~\eqref{E:def B_0}
for the definition of~$B_0$)
\begin{eqnarray}
  P^{cb}_-(s_1;\lambda,\tilde\rho) U_0^{cb} & = & U_0^{cb} +
  \int_{-\infty}^{s_1} \Phi^{cs}_-(s_1,\tau;\lambda_0,0)
  r^\prime(\tau)^2(\lambda-\lambda_0-\tilde\rho(\tau))B_0
  U^{cb}_-(\tau;\lambda,\tilde\rho,U_0^{cb})\, d\tau \\
  P_-^u(s_1;\lambda,\tilde \rho) & = & P_-^u(s_1;\lambda_0,0) +
  \int_{-\infty}^{s_1} \Phi_-^{cs} (s_1,\tau;\lambda_0,0)
  r^\prime(\tau)^2(\lambda-\lambda_0-\tilde\rho(\tau))B_0
  \Phi_-^u(\tau,s_1;\lambda,\tilde \rho)\, d\tau \\ \label{E:projection+}
  P_+^s(s_1;\lambda,0) & = & P_+^s(s_1;\lambda_0,0) -
  \int_{s_1}^\infty \Phi_+^{cu}(s_1,\tau;\lambda_0,0)r^\prime(\tau)^2
  (\lambda-\lambda_0)B_0 \Phi_+^{s}(\tau,s_1;\lambda,0)\, d\tau,
\end{eqnarray}
where we recall that $r(\tau)=e^\tau$ for $\tau<s_2$, so that
$r'(\tau)^2 = e^{2\tau}$ in this interval.
Hence we may write
\begin{eqnarray}\label{E:iota-def2}
  \iota(U_0^s,U_0^{u}+U_0^{cb};\lambda,\tilde \rho) = U_0^{u} +
  U_0^{cb} - U_0^s + \int_{s_1}^\infty
  \Phi_+^{cu}(s_1,\tau;\lambda_0,0) r^\prime(\tau)^2
  (\lambda-\lambda_0)B_0 \Phi_+^{s}(\tau,s_1;\lambda,0) U_0^s\, d\tau
  \\ \nonumber
  + \int_{-\infty}^{s_1} \Phi^{cs}_-(s_1,\tau;\lambda_0,0)
  e^{2\tau}(\lambda-\lambda_0-\tilde\rho(\tau))B_0 \left[
    U^{cb}_-(\tau;\lambda,\tilde\rho,U_0^{cb})+\Phi_-^u(\tau,s_1;\lambda,\tilde
    \rho)U_0^u\right]\, d\tau.
\end{eqnarray}

\begin{Lem}\label{L:iota-smooth}
  The map $\iota:E_+^s\times (E_-^u\oplus E_-^{cb}) \times \mathbb R \times
  \tilde{\mathcal R} \to \overline X$ is smooth.
\end{Lem}
\begin{proof}
  {We have seen in Theorem \ref{T:dich-nonsmooth} that
  $(U_0^u,\lambda,\tilde \rho)\mapsto P_-^u(s_1;\lambda,\tilde \rho)U_0^u$ and
  $(U_0^{cb},\lambda,\tilde \rho)\mapsto P_-^{cu}(s_1;\lambda,\tilde \rho)U_0^{cb}$
  are smooth as functions from $E_-^u\times \mathbb R \times \widetilde {\mathcal R}$ to
  $X\subset \overline X$ and from $E_-^{cb}\times \mathbb R \times \widetilde {\mathcal R}$ to
  $X\subset \overline X$, respectively. }

  {Hence it suffices to prove that $P_+^s(s_1;\lambda,0)$ is smooth from
  $E_+^s\times \mathbb R$ to $\overline X$. We do this by studying the terms of \eqref{E:projection+}
  separately. It is clear that $U_0^s\in \overline
  X$.}

  {Next, we study the integral term, and }
  note that by Lemma~\ref{L:far field dichotomy}, for
  $\tau\ge s_1$, $\Phi_+^s(\tau,s_1;\lambda,0):\mathcal X^{s_1}\to
  \mathcal X^\tau$ with norm bounded by $K
  e^{(\lambda^{1/4}-\epsilon)(\tau-s_1)}$ and
  $\Phi_+^{cu}(s_1,\tau;\lambda_0,0):\mathcal Y^\tau\to \mathcal
  Y^{s_1}$ with norm bounded by $Ke^{\epsilon(\tau-s_1)}$. Recall that
  $\mathcal X^s= H^1\times L^2\times H^1\times L^2$ and
  $\mathcal Y^s = H^2\times H^1\times H^2\times H^1$ with $s$-dependent norms. Thus
  $B_0:\mathcal X^\tau\to \mathcal Y^\tau$ is bounded with norm~$\tau$.
  Using the exponential estimates of Lemma~\ref{L:far field dichotomy}
  and that $\mathcal Y^{s_1}\subset X\subset\overline X\subset\mathcal
  X^{s_1}$, there exists a constant $C>0$ such that for any
  $\epsilon>0$ sufficiently small
  \begin{equation*}
    \begin{aligned}
      \left\| \int_{s_1}^\infty \Phi_+^{cu}(s_1,\tau;\lambda_0,0)
        r^\prime(\tau)^2B_0 \Phi_+^{s}(\tau,s_1;\lambda,0) U_0^s\,
        d\tau \right\|_{\overline X} &\le C\,\left\| \int_{s_1}^\infty
        \Phi_+^{cu}(s_1,\tau;\lambda_0,0) B_0
        \Phi_+^{s}(\tau,s_1;\lambda,0) U_0^s\, d\tau
      \right\|_{\mathcal Y^{s_1}} \\
      &\le CK^2 e^{(\lambda^{1/4}-2\epsilon)s_1} \int_{s_1}^\infty
      \tau \,e^{-(\lambda^{1/4}-2\epsilon)\tau} d\tau
      \|U_0^s\|_{\mathcal X^{s_1}} \\
      &\le \frac{\tilde{C}K^2}{(\lambda^{1/4} - 2\epsilon)^2}
      \|U_0^s\|_{\overline X},
    \end{aligned}
  \end{equation*}
for some constant $\tilde C$.
  To show that the integral term is smooth in $\lambda$ into
  $\overline X$, it suffices to prove that for $n\ge 1$
  \begin{equation*}
    \int_{s_1}^\infty \Phi_+^{cu}(s_1,\tau;\lambda_0,0) r^\prime(\tau)^2B_0
    \frac{d^n}{d\lambda^n}\Phi_+^s(\tau,s_1;\lambda,0) U_0^s\, d\tau
   \end{equation*}
   belongs to $\overline X$. This follows since, by Lemma~\ref{L:far field
     dichotomy}, $\frac{d^n}{d\lambda^n}\Phi_+^s(\tau,s_1;\lambda,0)$
   satisfies a similar exponential decay estimate as
   $\Phi_+^s(\tau,s_1;\lambda,0)$ does. It follows as above that the
   integral in question converges in $\mathcal Y^{s_1}\subset\overline
   X$.  Smoothness in $U_0^s$ is immediate, since $\iota$ is bounded
   and linear in $U_0^s$ into $\overline X$.
\end{proof}

\begin{Lem}\label{L:iota-eq}
  The operator $\Delta^2 + \theta + \rho$ has an embedded eigenvalue
  $\lambda>0$ if and only if there exist $U_0^s\in E_+^s$, $U_0^{u}\in
  E_-^{u}$ and $U_0^{cb}\in E_-^{cb}$ such that
  \begin{equation}\label{E:iota-eq}
    \iota(U_0^s,U_0^{u} + U_0^{cb};\lambda,\tilde \rho) = 0.
  \end{equation}
\end{Lem}
\begin{proof}
  If $\lambda$ is an eigenvalue of $\mathcal L + \rho$, then by
  Lemma~\ref{L:system-eigenfunction}, the corresponding solution of
  the system \eqref{E:perturbed system} is bounded as $s\to -\infty$
  and decays exponentially as $s\to +\infty$.  Hence there exists a
  solution of \eqref{E:perturbed system} with initial condition
  \begin{equation*}
    P^s_+(s_1;\lambda,\tilde \rho) U_0^s = P_-^{u}(s_1;\lambda,\tilde \rho) U_0^{u} + P_-^{cb}(s_1;\lambda,\tilde \rho) U_0^{cb}
  \end{equation*}
  at $s=s_1$, i.e.\ \eqref{E:iota-eq} holds.

  Conversely, suppose that \eqref{E:iota-eq} is satisfied for some
  $(U_0^s,U_0^{u}+U_0^{cb},\lambda,\tilde \rho)\in E_+^s\times
  (E_-^u\oplus E_-^{cb})\times \mathbb R\times\tilde{\mathcal
    R}$. Then there exists a solution of \eqref{E:perturbed system} with
  initial condition $U_0^s = U_0^{u}+ U_0^{cb}$.  By
  Lemma~\ref{L:system-eigenfunction}, this implies that $\lambda$ is
  an eigenvalue of $\mathcal L + \rho$.
\end{proof}

\begin{Lem}\label{L:decouple}
  The subspaces $E_-^{cb}\oplus E_-^{u}$ and $E_+^{s}$ have
  complements in
  $\overline{X}$ denoted by $E_-^{ca}\oplus E_-^s$ and
  $E_+^{cu}$. Moreover, $(E_-^{ca}\oplus E_-^s) \cap E_+^{cu}$ is
  infinite-dimensional and has a basis with elements $V_{0,3}(s_1)$,
  $V_{0,4}(s_1)\in X_0$, $V_{k,4}(s_1)\in X_k$, $k\in \mathbb
  Z\setminus \{0\}$.
\end{Lem}
\begin{proof}
  Recall that $E_-^{ca}= \span\{V_{0,3}(s_1),V_{0,4}(s_1)\}$.  Let
  \begin{equation*}
    \begin{aligned}
      E_-^{s}:&=\overline {\span\{V_{k,1}(s_1),V_{k,4}(s_1);\;
        k\in\mathbb Z\setminus\{0\} \}}, \\
      E_+^{cu}:&=\overline{\span\{V_{k,2}(s_1),V_{k,3}(s_1),V_{k,4}(s_1);\;
        k\in \mathbb Z\}},
    \end{aligned}
  \end{equation*}
  where the closures are taken in $\overline X$. It is easy to see
  that these spaces have the desired properties.
\end{proof}

Let $Q$ be the projection in $\overline X$ onto $\Ran \iota(\cdot,\cdot;\lambda_0,0) = E_+^s + (E_-^u\oplus E_-^{cb})$ such that
\begin{equation*}
   \ker Q = E_+^{cu} \cap (E_-^{ca} \oplus E_-^s).
\end{equation*}
Note that $\Ran Q$ and $\ker Q$ are closed subspaces of $\overline X$,
and $Q$ is therefore  continuous.

Equation \eqref{E:iota-eq} is equivalent to the pair of equations
\begin{equation}\label{E:iota-split}
   \begin{aligned}
      Q\iota(U_0^s,U_0^{u} + U_0^{cb};\lambda,\tilde \rho) &= 0,\\
      (I-Q)\iota(U_0^s,U_0^{u} + U_0^{cb};\lambda,\tilde \rho) &= 0.
   \end{aligned}
\end{equation}
\begin{Lem}\label{L:Lyapunov1}
  For $(\lambda,\tilde \rho)$ in a neighbourhood of $(\lambda_0,0)\in
  \mathbb R\times \tilde{\mathcal R}$, the first equation of
  \eqref{E:iota-split} has a unique (up to constant multiples) nonzero
  solution $(U_0^s,U_0^{u} + U_0^{cb})$ which depends smoothly on
  $\lambda$ and $\tilde \rho$ in this neighbourhood.  We write
  $U_0^s(\lambda,\tilde \rho)$, $U_0^{u}(\lambda,\tilde \rho)$ and
  $U_0^{cb}(\lambda,\tilde\rho)$. Furthermore,
  $U_0^s(\lambda_0,0)=U_*(s_1)=U_0^{cb}(\lambda_0,0)$ and
  $U_0^{u}(\lambda_0,0)=0$.
\end{Lem}
\begin{proof}
  For $(\lambda,\tilde \rho)$ fixed, $Q\iota$ is a linear mapping from
  $E_+^s\times (E_-^{u}\oplus E_-^{cb})$ to $\Ran Q$. It is clear that
  \begin{equation*}
     \ker Q\iota(\cdot,\cdot;\lambda_0,0) =  \span\{(U_*(s_1),U_*(s_1))\}.
  \end{equation*}
  By Lemma~\ref{L:iota-smooth} and
  since $\Ran Q$ is closed, it follows that $Q \iota$ is a smooth mapping in its
  arguments. Let $D$ be an affine hyperplane of $E_+^s\times
  (E_-^{u}\oplus E_-^{cb})$ such that $D\cap
  \span\{(U_*(s_1),U_*(s_1))\} = \{(U_*(s_1),U_*(s_1))\}$. The
  implicit function theorem then implies that for $(\lambda,\tilde \rho)$
  close to $(\lambda_0,0)$ the first equation of~\eqref{E:iota-split}
  has a unique solution $(U_0^s,U_0^{u} + U_0^{cb}) =  (U_0^s(\lambda,\tilde \rho),U_0^{u}(\lambda,\tilde
  \rho)+U_0^{cb}(\lambda,\tilde \rho))\in D$ in a neighbourhood of
  $(U_*(s_1),U_*(s_1))$. Moreover, $U_0^s$, $U_0^{u}$ and $U_0^{cb}$
  are smooth in their arguments.
\end{proof}

For $(\lambda, \tilde \rho)$ in the neighbourhood obtained in
Lemma~\ref{L:Lyapunov1}, we let
\begin{equation*}
  \begin{aligned}
    F(\lambda,\tilde \rho) :&= (I-Q) \iota(U_0^s(\lambda,\tilde
    \rho),U_0^{u}(\lambda,\tilde \rho)+U_0^{cb}(\lambda,\tilde\rho)) \\
    &=  (I-Q)
    \int_{-\infty}^{s_1} \Phi^{cs}_-(s_1,\tau;\lambda_0,0)
    e^{2\tau}(\lambda-\lambda_0-\tilde\rho(\tau))B_0 \left[
      U^{cb}_0(\tau;\lambda,\tilde\rho)+\Phi_-^u(\tau,s_1;\lambda,\tilde
      \rho)U_0^u(\lambda,\tilde\rho)\right]\, d\tau \\
    &+ (I-Q) \int_{s_1}^\infty
    \Phi_+^{cu}(s_1,\tau;\lambda_0,0)r^\prime(\tau)^2
    (\lambda-\lambda_0)B_0 \Phi_+^{s}(\tau,s_1;\lambda,0)
    U_0^s(\lambda,\tilde\rho)\, d\tau \\
   \end{aligned}
\end{equation*}
where $U^{cb}_0(s;\lambda,\tilde\rho)$ corresponds to
$U^{cb}_-(s,\lambda,\tilde\rho,U^{cb}_0(\lambda,\tilde\rho))$ so that
$U^{cb}_0(s;\lambda_0,0)=U_*(s)$. We see that solving
\eqref{E:iota-split} is equivalent to solving $F(\lambda,\tilde
\rho)=0$.

On $\Ran (I-Q)\subset \overline X$, the $\overline X$-norm is the same
as the $\mathcal X$-norm, and so we solve $F(\lambda,\tilde \rho)=0$
in $\mathcal X$.
For $k\in \mathbb Z\setminus\{0\}$ let $F_k(\lambda,\tilde
\rho):=\langle W_{k,4}(s_1),F(\lambda,\tilde \rho)\rangle$, and for
$k=0$ and $j=3, 4$ we let $F_{0,j}(\lambda,\tilde \rho):= \langle
W_{0,j}(s_1),F(\lambda,\tilde \rho)\rangle$.  Define $W_{k,4}(s)$ so
that it satisfies the adjoint system~\eqref{E:adjoint}
(and hence~\eqref{E:Fourier_unperturbed_adjoint}) and passes through
$W_{k,4}(s_1)$ for $s=s_1$. As
$W_{k,4}(s_1)\in\Ran(\Phi_{k,+}^{cu}(s_1,s_1)^*) \cap
\Ran(\Phi_{k,-}^{s}(s_1,s_1)^*)=\Ran(\widehat\Phi_{k,+}^{cs}(s_1,s_1)) \cap
\Ran(\widehat\Phi_{k,-}^{u}(s_1,s_1))$, we get for
$k\ne 0$
\begin{equation*}
  \begin{aligned}
    F_k(\lambda,\tilde \rho) &= \int_{-\infty}^{s_1} \left\langle W_{k,4}(s_1)\,,\,
      \Phi^{cs}_-(s_1,\tau;\lambda_0,0)\,
      e^{2\tau}(\lambda-\lambda_0-\tilde\rho(\tau)) B_0 \left[
        U^{cb}_0(\tau;\lambda,\tilde\rho)+\Phi_-^u(\tau,s_1;\lambda,\tilde
        \rho)U_0^u(\lambda,\tilde\rho)\right]\right\rangle\, d\tau \\
     & \qquad + \int_{s_1}^\infty \langle
    W_{k,4}(s_1), \Phi_+^{cu}(s_1,\tau;\lambda_0,0)r^\prime(\tau)^2
    (\lambda-\lambda_0)B_0 \Phi_+^{s}(\tau,s_1;\lambda,0)
    U_0^s(\lambda,\tilde\rho) \rangle\, d\tau
    \\
    &=
    \int_{-\infty}^{s_1} \left\langle W_{k,4}(\tau),
     e^{2\tau}(\lambda-\lambda_0-\tilde\rho(\tau))B_0 \left[
        U^{cb}_0(\tau;\lambda,\tilde\rho)+\Phi_-^u(\tau,s_1;\lambda,\tilde
        \rho)U_0^u(\lambda,\tilde\rho)\right]\right\rangle\, d\tau \\
        &\qquad + \int_{s_1}^\infty \langle W_{k,4}(\tau),r^\prime(\tau)^2
    (\lambda-\lambda_0)B_0 \Phi_+^{s}(\tau,s_1;\lambda,0)
    U_0^s(\lambda,\tilde\rho) \rangle\, d\tau
   \end{aligned}
\end{equation*}
and for $k=0$ and $j=3, 4$ we have similarly
\begin{eqnarray}\label{E:F_0j-eqn}
  F_{0,j}(\lambda,\tilde \rho) &=&
  \int_{-\infty}^{s_1} \left\langle W_{0,j}(\tau), e^{2\tau}(\lambda-\lambda_0-\tilde\rho(\tau))B_0 \left[ U^{cb}_0(\tau;\lambda,\tilde\rho)+\Phi_-^u(\tau,s_1;\lambda,\tilde \rho)U_0^u(\lambda,\tilde\rho)\right]\right\rangle\, d\tau  \\ \nonumber
  && + \int_{s_1}^\infty \langle W_{0,j}(\tau), r^\prime(\tau)^2(\lambda-\lambda_0)B_0 \Phi_+^{s}(\tau,s_1;\lambda,0) U_0^s(\lambda,\tilde\rho) \rangle\, d\tau
\end{eqnarray}

\begin{Lem}
  For $(\lambda,\tilde \rho)$ in a neighbourhood of $(\lambda_0,0)\in
  \mathbb R\times \tilde{\mathcal R}$, the equation
  \eqref{E:iota-split} has a nontrivial solution $(U_0^s(\lambda,\tilde
  \rho),U_0^u(\lambda,\tilde \rho)+U_0^{cb}(\lambda,\tilde
  \rho),\lambda,\tilde \rho)\in E_+^s\times (E_-^u\oplus
  E_-^{cb})\times \mathbb R\times \tilde{\mathcal R}$ if and only if
  $F_k(\lambda,\tilde \rho)=0$ for $k\in \mathbb Z\setminus\{0\}$ and
  $F_{0,j}(\lambda,\tilde \rho)=0$ for $j=3, 4$.
\end{Lem}
\begin{proof}
  Suppose that $(U_0^s(\lambda,\tilde \rho),U_0^u(\lambda,\tilde
  \rho)+U_0^{cb}(\lambda,\tilde \rho),\lambda, \tilde\rho)\in
  E_+^s\times (E_-^u\oplus E_-^{cb})\times \mathbb R\times
  \tilde{\mathcal R}$ solves \eqref{E:iota-split}. It is then clear
  from the definition of $F_k$ that $F_k(\lambda,\tilde \rho)=0$ for
  $k\in \mathbb Z\setminus\{0\}$ and $F_{0,j}(\lambda,\tilde \rho)=0$
  for $j=3, 4$. Conversely, let $F_k(\lambda,\tilde \rho)=0$ for $k\in
  \mathbb Z\setminus\{0\}$ and $F_{0,j}(\lambda,\tilde \rho)=0$ for
  $j=3, 4$.
  {By Lemma \ref{L:Lyapunov1}, the first equation of \eqref{E:iota-split} is satisfied,
  so it remains to check the second equation of \eqref{E:iota-split}.
  Recall that the basis vectors in $\Ran(I-Q)$ are $V_{k,j}(s_1)$, where $j=4$ for $k\ne 0$ and $j=3, 4$ for $k=0$.
  The coefficients of $(I-Q)\iota(U_0^s,U_0^u+U_0^{cb};\lambda,\tilde \rho)$ with respect to this basis are then
  $F_{0,j}$ ($j=3,4$) and $F_{k}$, $k\in \mathbb Z\setminus \{0\}$. Since all these coefficients vanish, the
  conclusion follows.}
\end{proof}
\begin{Lem}\label{L:lambdaprim}
  {The equation $F_{0,3}(\lambda,\tilde \rho)=0$
  defines $\lambda$ as a smooth function of $\tilde \rho$ in a
  neighbourhood of $\tilde \rho=0$ such that $\lambda(0) =  \lambda_0$.
  {Furthermore,}
  \begin{equation*}
    \lambda'(0)\tilde \rho = -\frac{\int_{-\infty}^{s_1} u_*(\tau)^2 \widehat {\tilde \rho}_0(\tau)
      e^{2\tau}\, d\tau}{\int_{-\infty}^\infty \overline u_*(\tau)^2 r(\tau) r'(\tau)\, d\tau}
  \end{equation*}
  where $\widehat {\tilde\rho}_0$ is the Fourier coefficient of
  $\tilde\rho$ corresponding to $k=0$. }
\end{Lem}
\begin{proof}
{By Lemma~\ref{L:iota-smooth} it follows that $F_{0,3}$ is a smooth
  function of $\lambda$ and $\tilde \rho$ in a neighbourhood of
  $(\lambda_0,0)$.  Note that
   \begin{equation*}
        \frac{\partial F_{0,3}}{\partial \lambda}(\lambda_0,0)
        =\int_{-\infty}^{\infty} \langle W_{0,{{3}}}(\tau), r'(\tau)^2 B_0
        U_*(\tau)\rangle\, d\tau
        = \int_{-\infty}^{\infty} \overline{w_{0,{3}}(\tau)} u_*(s)
        r'(\tau)^2\, ds, 
   \end{equation*}
where we have used that
$U_0^{s}(\lambda_0,0)=U_*(s)=U_0^{cb}(\lambda_0,0)$ and
$U_0^{u}(\lambda_0,0)=0$. {We}
recall that $W_{0,3}(s_1)=U_*^{\perp}(s_1)$ and that $W_{0,3}(s)$ satisfies the adjoint system
\eqref{E:adjoint} for $s\in \mathbb R$. It can be verified that
\begin{equation*}
   W_{0,3}(s) = \frac{r(s)}{r'(s)}U_*^\perp(s)
\end{equation*}
for $s\in \mathbb R$, where $U^\perp = (-u_4,u_3,-u_2,u_1)$, and $u_j$
are the components of $U$, $j=1,\dots,4$. {Hence
   \begin{equation}\label{E:lambdaderivative}
        \frac{\partial F_{0,3}}{\partial \lambda}(\lambda_0,0)
        = \int_{-\infty}^{\infty}  u_*(s)^2
        r(\tau) r'(\tau)\, ds>0.
   \end{equation}}
The last inequality follows since the integral is positive (using that $u_*$ is an eigenfunction).}

{
Since $\partial F_{0,3}/\partial \lambda(\lambda_0,0)\ne 0$, we can  solve the equation $F_{0,3}(\lambda,\tilde\rho)=0$ by
the implicit function theorem for $\lambda$ as a function of $\tilde\rho$, and this solution is a smooth function
$\lambda(\tilde \rho)$, defined in a neighbourhood of $\tilde \rho=0$,
such that $\lambda(0)=\lambda_0$, and $\lambda'(0)$ is given by
\begin{equation*}
  \begin{aligned}
    \lambda'(0) \tilde \rho &= -\frac{\partial F_{0,3}}{\partial
      \tilde \rho}(\lambda_0,0) \tilde \rho\bigg/\frac{\partial
      F_{0,3}}{\partial \lambda}(\lambda_0,0)
    = -\frac{\int_{-\infty}^{s_1} u_*(\tau)^2
      \widehat {\tilde \rho}_0(\tau) e^{2\tau}\,
      d\tau}{\int_{-\infty}^\infty u_*(\tau)^2
      r(\tau) r'(\tau)\, d\tau}
  \end{aligned}
\end{equation*}
as claimed.}
\end{proof}

Since we have solved $F_{0,3}=0$ for $\lambda$ in terms of
$\tilde \rho$, the remaining equation corresponding to $k=0$ is
$F_{0,4}(\lambda(\tilde\rho),\tilde\rho)=0$. We define
\begin{equation*}
   \begin{aligned}
     G_0(\tilde \rho) :&=
     \int_{-\infty}^{s_1} \left\langle W_{0,4}(\tau),
       e^{2\tau}(\lambda(\tilde\rho)-\lambda_0-\tilde\rho(\tau))B_0
       \left[ U^{cb}_0(\tau;\lambda(\tilde\rho),\tilde\rho) +
         \Phi_-^u(\tau,s_1;\lambda(\tilde\rho),\tilde
         \rho)U_0^u(\lambda(\tilde\rho),\tilde\rho)\right]\right\rangle\,
     d\tau \\
     &\qquad + \int_{s_1}^\infty \langle W_{0,4}(\tau),
     r^\prime(\tau)^2(\lambda(\tilde\rho)-\lambda_0)B_0
     \Phi_+^{s}(\tau,s_1;\lambda(\tilde\rho),0)
     U_0^s(\lambda(\tilde\rho),\tilde\rho) \rangle\, d\tau,
   \end{aligned}
\end{equation*}
and for $k\ne 0$,
\begin{equation*}
  \begin{aligned}
    G_k(\tilde \rho) :&=  F_k(\lambda(\tilde \rho),\tilde \rho) \\
    &= \int_{-\infty}^{s_1} \left\langle W_{k,4}(\tau),
      e^{2\tau}(\lambda(\tilde\rho)-\lambda_0-\tilde\rho(\tau))B_0
      \left[ U^{cb}_0(\tau;\lambda(\tilde\rho),\tilde\rho) +
        \Phi_-^u(\tau,s_1;\lambda(\tilde\rho),\tilde
        \rho)U_0^u(\lambda(\tilde\rho),\tilde\rho)\right]\right\rangle\,
    d\tau. \\
    &\qquad +
    \int_{s_1}^\infty \langle W_{k,4}(\tau),r^\prime(\tau)^2
    (\lambda(\tilde\rho)-\lambda_0)B_0
    \Phi_+^{s}(\tau,s_1;\lambda(\tilde\rho),0)
    U_0^s(\lambda(\tilde\rho),\tilde\rho) \rangle\, d\tau.
  \end{aligned}
\end{equation*}

\begin{Lem}\label{L:G-mapping}
   The mapping $\mathcal G: \tilde {\mathcal R}\to l_1^2$ defined by
   \begin{equation*}
       \mathcal G(\tilde \rho) = \{G_k(\tilde \rho)\}_{k\in \mathbb Z}
   \end{equation*}
is smooth.
\end{Lem}
\begin{proof}
  We first verify that the range of $\mathcal G$ belongs to $l_1^2$. To do
  this, we split the expression for $G_k(\tilde \rho)$ ($k\ne 0$) into
  three terms, which we deal with separately:
  \begin{equation}\label{e:G}
    \begin{aligned}
      G_k(\tilde \rho) = & (\lambda(\tilde \rho)-\lambda_0)
      \int_{-\infty}^{s_1} \left\langle W_{k,4}(\tau), e^{
          2\tau} U^{cu}(\tau) \right\rangle\, d\tau
      -\int_{-\infty}^{s_1}\left\langle W_{k,4}(\tau), e^{
          2\tau}\tilde \rho(\tau) U^{cu}(\tau) \right\rangle \, d\tau \\
      & + (\lambda(\tilde\rho)-\lambda_0) \int_{s_1}^\infty
      \left\langle W_{k,4}(\tau), r'(\tau)^2U^s(\tau)\right\rangle \, d\tau,
    \end{aligned}
  \end{equation}
where we used the notation
\begin{equation*}
  \begin{aligned}
    U^{cu}(\tau) & := [U^{cb}_0(\tau;\lambda(\tilde\rho),\tilde\rho) + \Phi_-^u(\tau,s_1;\lambda(\tilde\rho),\tilde\rho) U_0^u(\lambda(\tilde\rho),\tilde\rho)], \\
    U^s(\tau) & := \Phi_+^s(\tau,s_1;\lambda(\tilde \rho),\tilde
    \rho) U_0^s(\lambda(\tilde \rho),\tilde \rho).
  \end{aligned}
\end{equation*}
Then ${B_0} U^{cu}(\tau)\in \{0\}\times \{0\}\times \{0\}\times H^2(S^1)$
and ${B_0} U^{s}(\tau)\in \{0\}\times \{0\}\times \{0\}\times H^1(S^1)$.
Furthermore, by 
{its construction}, we have that
$W_{k,4}(s_1)\in \Ran(\widehat\Phi_{k,-}^{u}(s_1,s_1)) \cap
\Ran(\widehat\Phi_{k,+}^{cs}(s_1,s_1))$ for all $k\in
\mathbb{Z}\backslash\{0\}$.  Thus Lemma~\ref{L:adjoint_upper_bound_J+}
implies that for any $\epsilon>0$ there exists a constant $K$ such that
for every $k\in\mathbb Z\backslash\{0\}$ and $s\ge s_1$, {
\begin{equation*}
   \|W_{k,4}(s)\|_{(\mathcal{X}^s)'}\leq K e^{\epsilon(s-s_1)}
\end{equation*}
as the norms on $\mathcal X'$ and $(\mathcal X^{s_1})'$ are equivalent and
$\|W_{k,4}(s_1)\|_{\mathcal{X}'}=1$.}

Now observe that $B_0 U^{s}(\tau)$ vanishes for all components except the
last one, so we only need an estimate on the last component of
$W_{k,4}(s)$, which we denote by $w_{k,4}$. Then
the estimate above gives that there exists a constant $K$ independent
of $k$ such that
\begin{equation}
  \label{E:bound_w_J+}
  |w_{k,4}(s)|\leq K e^{\epsilon(s-s_1)}, \mbox{ for all }
  k\in\mathbb Z\backslash\{0\} \mbox{ and } s\ge s_1,
\end{equation}
as $\mathcal{X}'\equiv H^{-1}\times L^{2}\times H^{-1}\times L^{2}$.
Similarly, from $W_{k,4}(s_1)\in
\Ran(\widehat\Phi_{k,-}^{u}(s_1,s_1))$,
Lemma~\ref{L:adjoint_upper_bound_J-}, $X'=H^{-2}\times H^{-1}\times
H^{-1}\times L^{2}$ and $\|W_{k,4}(s_1)\|_{X'}\leq
\|W_{k,4}(s_1)\|_{\mathcal{X}'}=1$, it follows that there is some
constant $K$ such that
\begin{equation}
  \label{E:bound_w_J-}
|w_{k,4}(s)|\le K e^{|k|(s-s_1)} \mbox{ for all } k\in\mathbb
Z\backslash\{0\} \mbox{ and } s\le s_1.
\end{equation}

First we look at the last integral in \eqref{e:G}.
Let $\widehat{u_k^s}(\tau)$ be the first component of the $k$-th
Fourier coefficient of $\Phi_+^s(\tau,s_1;\lambda(\tilde \rho),\tilde
\rho) U_0^s(\lambda(\tilde \rho),\tilde \rho)$, then the definition of $B_0$ in \eqref{E:def B_0} gives
\begin{equation*}
  \int_{s_1}^\infty \left\langle W_{k,4}(\tau),r'(\tau)^2
    B_0 U^s(\tau)\right\rangle \, d\tau = \int_{s_1}^\infty
  r'(\tau)^2\overline{w_{k,4}(\tau)}\, \widehat{u_k^s}(\tau)\, d\tau.
\end{equation*}
From Lemma~\ref{L:far field dichotomy}, it follows that, for any
$\epsilon>0$ and $\tau\ge s_1$,
\[
\|\Phi_+^s(\tau,s_1;\lambda(\tilde \rho),\tilde \rho)
U_0^s(\lambda(\tilde \rho),\tilde \rho)\|_{\mathcal{X}^\tau} \leq K
e^{-(\lambda(\tilde\rho)^{1/4}-\epsilon)(\tau-s_1)} \|U_0^s(\lambda(\tilde
\rho),\tilde \rho)\|_{\mathcal{X}^{s_1}}.
\]
This implies for the Fourier coefficients $\widehat{u_k^s}$ that
\[
\sum_{k\in\mathbb{Z}\backslash\{0\}}
\left(1+\frac{k^2}{\tau^2}\right)\left|\widehat{u_k^s}(\tau)\right|^2
\leq K^2 e^{-2(\lambda(\tilde\rho)^{1/4}-\epsilon)(\tau-s_1)} \|U_0^s(\lambda(\tilde
\rho),\tilde \rho)\|^2_{\mathcal{X}^{s_1}}.
\]
Combining this with~\eqref{E:bound_w_J+}, we see that for any
$\tilde\rho$, {we have
\begin{equation*}
  \begin{aligned}
    \sum_{k\in \mathbb Z\setminus\{0\}} (1+k^2)&\left( \int_{s_1}^\infty
      r'(\tau)^2\overline{w_{k,4}(\tau)} \widehat u_k^s(\tau)\, d\tau\right)^2
    \le  C \sum_{k\in \mathbb Z\setminus\{0\}}(1+k^2) \left(\int_{s_1}^\infty K
      e^{\epsilon(\tau-s_1)}
      |\widehat u_k^s(\tau)|\, d\tau \right)^2
    \\
    &\le C \sum_{k\in \mathbb Z\setminus\{0\}}\left(\int_{s_1}^\infty
      e^{-\epsilon(\tau-s_1)}d\tau\right) \,
    \left(\int_{s_1}^\infty e^{3\epsilon(\tau-s_1)} (1+k^2)
      |\widehat u_k^s(\tau)|^2\, d\tau \right)\\
    &\le \frac{C }{\epsilon} \int_{s_1}^\infty
    e^{3\epsilon(\tau-s_1)}\sum_{k\in \mathbb Z\setminus\{0\}}
    2\tau^2\left(1+\frac{k^2}{\tau^2}\right) |\widehat
    u_k^s(\tau)|^2\,d\tau\\
    &\le C \int_{s_1}^\infty
    \tau^2\,
    e^{-2(\lambda(\tilde\rho)^{1/4}-4\epsilon)(\tau-s_1)} \|U_0^s(\lambda(\tilde
\rho),\tilde \rho)\|^2_{\mathcal{X}^{s_1}}\,d\tau \le C,
  \end{aligned}
\end{equation*}
where $C=C(\epsilon)$ denotes the different constants occuring.}

Next, we look at the first integral in \eqref{e:G}. Let $\widehat{u_k^{cu}}$ be the
first component of the $k$-th Fourier coefficient of
$U^{cb}_0(\tau;\lambda(\tilde\rho),\tilde\rho) +
\Phi_-^u(\tau,s_1;\lambda(\tilde\rho),\tilde\rho)
U_0^u(\lambda(\tilde\rho),\tilde\rho)$.  The definition of $B_0$ gives
that the first integral can be written as
\[
\int_{-\infty}^{s_1} \left\langle W_{k,4}(\tau), e^{2\tau} B_0
  U^{cu}(\tau) \right\rangle\, d\tau = \int_{-\infty}^{s_1}
e^{2\tau}\, \overline{w_{k,4}(\tau)}\, \widehat{u_k^{cu}}(\tau)\,
d\tau.
\]
As $\Phi_-^u$ leads to
solutions with an $X$-norm that is exponentially decaying at $-\infty$
and $U^{cb}_0$ is bounded in the $X$-norm, there exists a constant $K$
such that the Fourier coefficients $\widehat{u_k^{cu}}(\tau)$ satisfy
\[
\sum_{k\in\mathbb{Z}\backslash\{0\}}
\left(1+{k^2}\right)\left|\widehat{u_k^{cu}}(\tau)\right|^2\le
\sum_{k\in\mathbb{Z}\backslash\{0\}}
\left(1+{k^2}\right)^2\left|\widehat{u_k^{cu}}(\tau)\right|^2\leq K^2
\]
for all $\tau\leq s_1$.   Together with the
fact that~\eqref{E:bound_w_J-} implies that $|w_{k,4}(\tau)|\leq
 K$ for all $\tau\le s_1$, this gives
\begin{equation*}
  \begin{aligned}
    \sum_{k\in \mathbb Z\setminus\{0\}} (1+k^2)&\left(
      \int_{-\infty}^{s_1} e^{2\tau} \, \overline{w_{k,4}(\tau)}\,
      \widehat{u_k^{cu}}(\tau) \, d\tau\right)^2 \le K^2 \sum_{k\in
      \mathbb Z\setminus\{0\}}(1+k^2)\left( \int_{-\infty}^{s_1} e^{2s_1}
      e^{2(\tau-s_1)} \,
      |\widehat{u_k^{cu}}(\tau)| \, d\tau\right)^2 \\
    &\le K^2 e^{2s_1} \sum_{k\in \mathbb Z\setminus\{0\}}
    \left(\int_{-\infty}^{s_1} e^{2(\tau-s_1)}d\tau\right) \,
    \left(\int_{-\infty}^{s_1} (1+k^2)e^{2(\tau-s_1)}
      |\widehat{u_k^{cu}}(\tau)|^2\, d\tau \right)\\
    &\le \frac{K^2 e^{2s_1}}{2} \,\int_{-\infty}^{s_1}
    e^{2(\tau-s_1)} \sum_{k\in \mathbb
      Z\backslash\{0\}}(1+k^2)|\widehat{u_k^{cu}}(\tau)|^2\,
    d\tau  \leq \frac{K^4 e^{2s_1}}{4}.
  \end{aligned}
\end{equation*}

Finally, let $\nu(\tau)$ be the first component of $\tilde \rho(\tau)
U^{cu}(\tau)$, so that $\nu(\tau) = \tilde \rho(\tau) u^{cu}(\tau)$. As
$u^{cu}\in H^2(S^1)$, its $H^2$ norm is uniformly bounded on
$(-\infty,s_1]$ and $\tilde \rho \in L^2(J_-;H^{1/2}(S^1),e^{2s}ds)$,
Lemma~\ref{L:product} implies that $\nu\in
L^2(J_-;H^{1/2}(S^1),e^{2s}ds)$.

Denote the Fourier coefficients of $\nu$ by $\widehat{\nu}_k$. Then
the second integral of \eqref{e:G} can be written as
\[
\int_{-\infty}^{s_1}\left\langle W_{k,4}(\tau), e^{ 2\tau}\tilde
  \rho(\tau) U^{cu}(\tau) \right\rangle \, d\tau =\int_{-\infty}^{s_1} e^{ 2\tau}\, \overline{w_{k,4}(\tau)} \,
\widehat{\nu}_k(\tau)\,d\tau
\]
and the estimate~\eqref{E:bound_w_J-} on the decay of $w_{k,4}$
implies that there exists a constant $C$ such that
\begin{equation*}
  \begin{aligned}
    \left( \int_{-\infty}^{s_1} e^{2\tau}\, \overline{w_{k,4}(\tau)}
      \, \widehat{\nu_k}(\tau)\,d\tau\right)^2 &\le
    \left(\int_{-\infty}^{s_1} K^2 e^{2(|k|+1)(\tau-s_1)}\right)
    \left(\int_{-\infty}^{s_1} e^{2\tau} \widehat{\nu_k}(\tau)^2\,d\tau\right)\\
    &\le \frac{K^2}{2(|k|+1)}\, \int_{-\infty}^{s_1}e^{2\tau}\widehat
    \nu_k(\tau)^2\, d\tau \le \frac{C}{\sqrt{1+k^2}}\,
    \int_{-\infty}^{s_1}e^{2\tau}\widehat \nu_k(\tau)^2\, d\tau,
  \end{aligned}
\end{equation*}
and so
\begin{equation*}
  \begin{aligned}
    \sum_{k\in \mathbb Z\setminus\{0\}} (1+k^2)&\left(
      \int_{-\infty}^{s_1} e^{2\tau}\, \overline{w_{k,4}(\tau)} \,
      \widehat{\nu_k}(\tau)\,d\tau\right)^2 \le C\sum_{k\in \mathbb Z\setminus\{0\}}
    (1+k^2)^{1/2}\int_{-\infty}^{s_1}e^{2\tau}\widehat \nu_k(\tau)^2\, d\tau
    \\&= C\|\nu\|_{L^2(J_-;H^{1/2}(S^1),e^{2s}ds)}^2 <\infty.
  \end{aligned}
\end{equation*}
Hence the second term also belongs to $l_1^2$, and so the proof of the
claim that the range of $\mathcal G$ is contained in $l_1^2$ is complete.

Smoothness follows since the integrands are smooth in $\tilde \rho$
and since the derivatives of arbitrary order of the evolution
operators $\Phi_+^s(\tau,s_1;\lambda(\tilde \rho),\tilde \rho)$ and
$\Phi_-^u(\tau,s_1;\lambda(\tilde \rho),\tilde \rho)$ belong to the
same exponentially weighted space as the evolution operators
themselves (see Theorem~\ref{T:dich-nonsmooth} and Lemma~\ref{L:far
  field dichotomy}).
\end{proof}

Finally we consider $\mathcal{G}'(0)$.  Since $U_*$ is radially
symmetric (and hence belongs to $X_0$) we have for $k\ne 0$ that
\begin{equation*}
  \begin{aligned}
    G_k'(0) \tilde \rho &= \int_{-\infty}^{\infty} \left\langle
    W_{k,4}(\tau),r'(\tau)^2 (\lambda'(0) \tilde \rho-
    \tilde\rho(\tau)) B_0
    U_*(\tau)\right\rangle\, d\tau    \\
    &= -\int_{-\infty}^{s_1} \overline{w_{k,4}(\tau)} \widehat {\tilde
      \rho}_k(\tau) u_*(\tau) e^{2\tau}\, d\tau,
  \end{aligned}
\end{equation*}
where $\widehat {\tilde\rho}_k$ is the $k$-th Fourier coefficient of
  $\tilde\rho$.
For  $k=0$ we have
\begin{equation*}
G_0'(0)\tilde \rho = -\int_{-\infty}^{s_1} \overline{w_{0,4}(s)} \,\widehat {\tilde\rho}_0(s) u_*(s) e^{2s}\, ds
- \frac{\int_{-\infty}^{\infty}\overline{w_{0,4}(\tau)} u_*(\tau) r'(\tau)^2 \, d\tau}{\int_{-\infty}^\infty u_*(\tau)^2 r(\tau)r'(\tau)\, d\tau}
    \int_{-\infty}^{s_1} \widehat {\tilde \rho}_0(s) u_*(s)^2 e^{2s}\, ds.
\end{equation*}
To rewrite the preceding expressions, we define
\begin{equation}\label{eta}
  \eta_k(s) := \overline{w_{k,4}(s)} u_*(s) \chi_{(-\infty,s_1]}(s)
\end{equation}
for $k\in \mathbb Z\setminus\{0\}$, and set
\begin{equation*}
  \eta_0(s) := \left[\overline{w_{0,4}(s)} + u_*(s)\,
    \frac{\int_{-\infty}^{\infty}\overline{w_{0,4}(\tau)} u_*(\tau)
      r'(\tau)^2\, d\tau}{\int_{-\infty}^\infty u_*(\tau)^2 r(\tau) r'(\tau)\, d\tau}\right] \, u_*(s)
  \chi_{(-\infty,s_1]}(s).
\end{equation*}
Then we may write
\[
\mathcal{G}'(0)\tilde \rho = \left\{ - \int_{-\infty}^{s_1} e^{2\tau}
  \eta_k(\tau) \widehat{\tilde \rho}_k(\tau) \, d\tau
\right\}_{k\in\mathbb{Z}}.
\]
For any $k\in\mathbb{Z}$, we have $\eta_k
e^{ik\cdot}\in\tilde{\mathcal{R}}$: indeed, \eqref{E:bound_w_J-}
shows that $|w_{k,4}(\tau)|\leq Ke^{|k|(\tau-s_1)}$ for any $\tau\leq
s_1$ and $k\in\mathbb{Z}\backslash\{0\}$, while $|w_{0,3}|$ and
$|w_{0,4}|$ are bounded on $J_-$, so that there exists a constant~$C$ such
that
\begin{equation}\label{E:upperbound eta_k}
  \int_{-\infty}^{s_1} |\eta_k(\tau)|^2\, e^{2\tau}d\tau \le
  \sup_{\tau\in (-\infty,s_1)} |u_*(\tau)|^2 \, C \int_{-\infty}^{s_1}
  e^{ (2|k|+2)(\tau - s_1)}\, d\tau \le \sup_{\tau\in (-\infty,s_1)}
  u_*(\tau)^2 \frac{C}{2|k|+2}.
\end{equation}

From the definition of $\mathcal{G}'(0)\tilde \rho$, it follows
immediately that $\mathcal G'(0)\tilde \rho= 0$ if and only if
\begin{equation}\label{keg}
\int_{-\infty}^{s_1} e^{2\tau}
  \eta_k(\tau) \widehat{\tilde \rho}_k(\tau) \, d\tau =0
\end{equation}
for all $k\in\mathbb{Z}$. Thus, if we define
\begin{equation*}
  \mathcal M := \overline{\span \{ \eta_k e^{ik\varphi};\; k\in \mathbb Z\} },
\end{equation*}
where the closure is taken in $\tilde {\mathcal R}$, then  it can be
seen that $\mathcal M$ is the orthogonal complement in $\tilde
{\mathcal R}$ of $\ker \mathcal G'(0)$, and so $\tilde {\mathcal R} = \ker \mathcal G'(0) \oplus \mathcal M$.

\begin{Lem}\label{L:decompose}
  $\mathcal G'(0):\mathcal M\to l^2_1$ is a linear homeomorphism. Furthermore, the spaces $\ker \mathcal G'(0)$ and $\mathcal M$ are both infinite-dimensional.
\end{Lem}

\begin{proof}
{It is clear that $\mathcal G'(0):\widetilde R \to l_1^2$ is bounded since by Lemma \ref{L:G-mapping}, $\mathcal G$ is smooth in a neighbourhood of $0$. }

We need to investigate the subspace $\cal M$.
Let $\eta\in \mathcal M$ be arbitrary, then
\begin{equation}\label{E:eta-expression}
  \eta(s,\varphi) = \sum_{k\in \mathbb Z}
  a_k \eta_k (s) e^{ik\varphi}.
\end{equation}
The upper bound estimate~\eqref{E:upperbound eta_k} implies that
\begin{equation}\label{E:eta-upper}
  \|\eta\|_{\tilde {\mathcal R}}^2 = \sum_{k\in \mathbb Z}
      (1 + k^2)^{1/2} |a_k|^2 \int_{-\infty}^{s_1} |\eta_k(\tau)|^2 \,e^{2\tau}
      d\tau \le C' \sum_{k\in \mathbb Z}
  |a_k|^2.
\end{equation}

Next we derive a lower bound for~$\|\eta\|_{\tilde {\mathcal
    R}}^2$. Since $u_*(s_1) \ne 0$, there exist $\widehat\epsilon$ and $\delta>0$
such that $u_*(s)^2 >{\widehat\epsilon}^2$ for every $s\in (s_1 -
\delta,s_1)$. Lemma~\ref{L:adjoint_lower_bound} shows that, for
$k$ large, $W_{k,4}$ is close to solutions of the system at infinity,
both in the $X$ and $\mathcal X$ norms.  This allows us to get a
lower bound on $|w_{k,4}(s)|$ for $k$ large. Let $\epsilon>0$ and $K$
as in Lemma~\ref{L:adjoint_lower_bound}. As $W_{k,4}(s_1)\in
\Ran(\widehat\Phi_k^{u}(s_1,s_1))$, it follows that
$\widehat\Phi_k^{u}(s,s_1)W_{k,4}(s_1)=W_{k,4}(s)$, and hence
\begin{equation*}
   \|W_{k,4}(s)-e^{|k|(s-s_1)}(P_k^{s})^*W_{k,4}(s_1)\|_{\mathcal{X}'}\leq
   \epsilon e^{|k|(s-s_1)}\|W_{k,4}(s_1)\|_{\mathcal{X}'}=\epsilon
   e^{|k|(s-s_1)}
\end{equation*}
for $|k|>K$.
 Thus we get for the fourth component~$w_{k,4}(s)$
  \begin{equation*}
    \begin{aligned}
      |w_{k,4}(s)| &\geq e^{|k|(s-s_1)}|((P_k^{s})^*W_{k,4}(s_1))_4| -
      |w_k(s)-e^{|k|(s-s_1)}((P_k^{s})^*W_{k,4}(s_1))_4| \\
      & \geq e^{|k|(s-s_1)}|((P_k^{s})^*W_{k,4}(s_1))_4|
      -\|W_{k,4}(s)-e^{|k|(s-s_1)}(P_k^{s})^*W_{k,4}(s_1)\|_{\mathcal{X}'}\\
      & \geq e^{|k|(s-s_1)}\left(|((P_k^{s})^*W_{k,4}(s_1))_4| -
        \epsilon \|W_{k,4}(s_1)\|_{\mathcal{X}'}\right)
      \geq e^{|k|(s-s_1)}\left(|w_{k,4}(s_1)| -
        2\epsilon \|W_{k,4}(s_1)\|_{\mathcal{X}'}\right)\\
      & \geq e^{|k|(s-s_1)}\left(|w_{k,4}(s_1)| - 2\epsilon \right).
    \end{aligned}
  \end{equation*}
  To get a lower bound on $w_{k,4}(s_1)$, first note that
  Lemma~\ref{L:adjoint_lower_bound} implies that
  $\|W_{k,4}(s_1)-(P_k^{s})^*W_{k,4}(s_1)\|_{\mathcal{X}'}\leq
  \epsilon$ and that
  $\Ran(P_k^{s})^*=\span\{(-|k|,1,0,0)^T,(0,0,-|k|,1)^T\}$. Thus there
  exist $\alpha_k,\,\beta_k\in\mathbb{C}$ and $W_k \in
  \mathcal{X}'_k$ with $\|W_k\|_{\mathcal{X}'} \leq 1$ such that
\[
W_{k,4}(s_1)  = \alpha_k(-|k|,1,0,0)^Te^{ik\cdot} +
\beta_k(0,0,-|k|,1)^Te^{ik\cdot} + \epsilon W_k.
\]

 As $V_{k,1}(s_1)\in
E_+^s$, its first component has to be a multiple of
$K_{|k|}(\lambda_0^{1/4}r(s_1))$ and hence
\[
V_{k,1}(s_1) = C_ke^{ik\cdot}\,
\begin{pmatrix}
K_{|k|}(\lambda_0^{1/4}r_1)\\
r'(s_1)\lambda_0^{1/4}K_{|k|}'(\lambda_0^{1/4}r_1)\\
\lambda_0^{1/2}K_{|k|}(\lambda_0^{1/4}r_1)\\
r'(s_1)\lambda_0^{3/4}K_{|k|}'(\lambda_0^{1/4}r_1)
\end{pmatrix} C_ke^{ik\cdot}\,\begin{pmatrix}\underline{v}_k\\\lambda_0^{1/2}\,\underline{v}_k
\end{pmatrix}
\mbox{ with } \underline v_k=\begin{pmatrix}
    K_{|k|}(\lambda_0^{1/4}r_1)\\
    r'(s_1)\lambda_0^{1/4}K_{|k|}'(\lambda_0^{1/4}r_1)
\end{pmatrix}
\]
where $r_1=r(s_1)$ and $C_k$ is such that $\|V_{k,1}(s_1)\|_{\mathcal
  X}=1$,
i.e.,
\[
C_k^2 = \frac{1}{(1+\lambda_0)
  \left[(k^2+1)\left(K_{|k|}(\lambda_0^{1/4}r_1)\right)^2 +
    \sqrt{\lambda_0}\left(r'(s_1)K_{|k|}'(\lambda_0^{1/4}r_1)\right)^2 \right]}.
\]
{Since
$\langle W_{k,4}(s_1), V_{k,1}(s_1)\rangle =0$, we have}
\begin{equation}\label{E:w_k estimate}
0 = C_k\,\left(-|k|K_{|k|}(\lambda_0^{1/4}r_1)+
\lambda_0^{1/4}r'(s_1)K_{|k|}'(\lambda_0^{1/4}r_1)\right)
\,(\alpha_k+\beta_k\sqrt{\lambda_0}) +
\epsilon \langle W_k,V_{k,1}(s_1)\rangle.
\end{equation}
From (9.6.23) in~\cite{AS72}, we see that $K_{|k|}(z)>0$ for any $z>0$
and (9.6.26) implies $K'_{|k|}(z) = -
K_{|k|-1}(z)-\frac{{|k|}}{z}\,K_{|k|}(z)$ for any $z>0$, hence
$K'_{|k|}(z)<0$ for any $z>0$. So we see that
$-|k|K_{|k|}(\lambda_0^{1/4}r_1)<0$ and
$\lambda_0^{1/4}r'(s_1)K_{|k|}'(\lambda_0^{1/4}r_1)<0$. {A short
analysis gives that
\begin{equation*}
 -2/\sqrt{1+\lambda_0}\le  C_k\,(-|k|K_{|k|}(\lambda_0^{1/4}r_1) + \lambda_0^{1/4}r'(s_1)K_{|k|}'(\lambda_0^{1/4}r_1))
\le -1/\sqrt{2(1+\lambda_0)}.
\end{equation*}}%
Furthermore, $|\langle W_k,V_{k,1}(s_1)\rangle|\leq
\|W_k\|_{\mathcal{X}'} \|V_{k,1}(s_1)\|_{\mathcal{X}}\leq 1$, and we
can conclude from~\eqref{E:w_k estimate} that
$\alpha_k=-\beta_k\sqrt{\lambda_0} + \mathcal{O}(\epsilon)$. Finally,
$\|W_{k,4}(s_1)\|_{\mathcal{X}'}=1$ gives
$\frac{1+2k^2}{1+k^2}\,(\alpha_k^2+\beta_k^2)=1-\mathcal{O}(\epsilon^2)$,
and hence
$(1+\lambda_0)\beta_k^2=\frac{1+k^2}{1+2k^2}-\mathcal{O}(\epsilon)$. Thus
there exists a $C>0$ such that $|w_{k,4}(s_1)|>C$ for all $|k|>K$. This
implies that there exists a $\tilde C>0$ such that $|w_{k,4}(s)|>\tilde
C\,e^{|k|(s-s_1)}$ for every $s\leq s_1$ and $|k|>K$.

Combining the lower bounds on $u_*(s)$ and $w_{k,4}(s)$, we find that
there exists a $\delta>0$ such that for $|k|>K$
\begin{equation}\label{E:eta-fourier}
  \begin{aligned}
    \int_{-\infty}^{s_1} \eta_k(\tau)^2 e^{2\tau}\, d\tau &\ge
    {\widehat\epsilon}^2 \tilde C^2 \int_{s_1 - \delta}^{s_1}
    e^{(2|k|+2) (\tau - s_1)}\, d\tau = \frac{\widehat\epsilon^2\tilde
      C^2}{2|k|+2}\left(1 - e^{-(2 |k|+2) \delta}
    \right)\\
    &\ge \frac{\widehat\epsilon^2\tilde C^2}{2|k|+2}\left(1 - e^{-2
        \delta}\right) \ge \frac{C}{(1 + k^2)^{1/2}},
  \end{aligned}
\end{equation}
for some positive $k$-independent constant $C$.  Since for all
$k\in\mathbb{Z}$, $\displaystyle\int_{-\infty}^{s_1} \eta_k(\tau)^2
e^{2\tau}\, d\tau>0$, {the constant $C$ above can be modified so that
\begin{equation*}
   \int_{-\infty}^{s_1} \eta_k(\tau)^2 e^{2\tau}\, d\tau\geq \frac{C}
{{(1 + k^2)^{1/2}}}\end{equation*}
also for $|k|\leq K$. Hence it follows that
\begin{equation*}
  \|\eta\|_{\tilde {\mathcal R}}^2 = \sum_{k\in \mathbb Z}
  (1 + k^2)^{1/2} |a_k|^2 \int_{-\infty}^{s_1} \eta_k(\tau)^2 \,
  e^{2\tau}d\tau
  \ge C\sum_{k\in \mathbb Z}|a_k|^2.
\end{equation*}}

The upper and lower bounds on $\|\eta\|_{\tilde {\mathcal R}}^2$ show
that $\eta \in \mathcal M$ if and only if $\eta$ is given by
\eqref{E:eta-expression} and $\{a_k\}_{k\in \mathbb Z}\in l^2$.

As we have seen that the mapping $\mathcal G'(0)$ is bounded above, it
is sufficient to show that it is bounded below to conclude that
$\mathcal G'(0)$ is a linear homeomorphism from $\mathcal M$ to
$l^2_1$. From its definition, it follows that
\[
\mathcal G'(0)\eta=\left\{a_k \int_{-\infty}^{s_1}
|\eta_k(\tau)|^2 e^{2\tau}\, d\tau \right\}_{k\in\mathbb{Z}},
\]
{and so by \eqref{E:eta-fourier} and \eqref{E:eta-upper}}, we see that
\begin{equation*}
    \|\mathcal G'(0)\eta\|_{l^2_1}^2 = \sum_{k\in \mathbb Z} (1 +
    k^2) |a_k|^2 \left( \int_{-\infty}^{s_1}
      |\eta_k(s)|^2e^{2s} \, ds
    \right)^2
    \ge \tilde C\sum_{k\in \mathbb Z} |a_k|^2 \ge C' \|\eta\|_{\tilde
      {\mathcal R}}^2.
\end{equation*}

It remains to show that the spaces $\ker \mathcal G'(0)$ and $\mathcal M$ have infinite dimension. For $\mathcal M$, this follows directly from its definition. Next, consider the characterization of $\ker \mathcal G'(0)$ given in (\ref{keg}). We proved above that the functions $w_{k,4}(s)$ that appear in the definition (\ref{eta}) of $\eta_k$ satisfy $|w_{k,4}(s_1)|\geq C$ uniformly in $|k|\geq K$, which implies that the space $\ker \mathcal G'(0)$ has infinite dimension as claimed.
\end{proof}

We are now ready to complete the proof of Theorem~\ref{T:main}.
\begin{proof}[Proof of Theorem~\ref{T:main}]
  By Lemma~\ref{L:iota-eq}, if $(\lambda,\tilde \rho)$ is sufficiently
  close to $(\lambda_0,0)$ then $\lambda$ is an embedded eigenvalue
  for $\Delta^2 + \tilde \theta + \tilde \rho$ if and only if
  \eqref{E:iota-eq} holds.  We have also seen that \eqref{E:iota-eq}
  is equivalent to $F(\lambda,\tilde \rho) = 0$, which allowed us to
  solve for $\lambda$ as a function of $\tilde \rho$ and finally
  obtain the equation $\mathcal G(\tilde \rho) =0$, where $\mathcal
  G:\tilde{\mathcal R} \to l^2_1$. By Lemma~\ref{L:decompose}, $\tilde
  {\mathcal R} = \ker \mathcal G'(0) \oplus \mathcal M$, and $\mathcal
  G'(0):\mathcal M\to l^2_1$ is a linear homeomorphism. Hence for
  $\tilde \rho\in \tilde {\mathcal R}$ we may write $\tilde \rho = \xi
  + \eta$, where $\xi\in \ker \mathcal G'(0)$ and $\eta\in \mathcal
  M$. By the implicit function theorem, we can solve for $\eta$ in
  terms of $\xi$, and this equation defines a smooth manifold in a
  neighbourhood of $0$ with infinite dimension and codimension, since
  the spaces $\ker \mathcal G'(0)$ and $\mathcal M$ are infinite-dimensional
  by Lemma~\ref{L:decompose}.
\end{proof}
\end{section}

\begin{section}{Conclusions and open problems}

In this paper, we considered the planar bilaplacian with a smooth, radially symmetric and compactly supported potential $\theta$ and described the set of perturbations of the potential in the space $\mathcal R = L^2([0,r_1];H^{1/2}(S^1),r\,dr)$ for which an embedded eigenvalue persists. We expect that the space $\mathcal R$ can be replaced by the Sobolev space $H^{1/2}(B_{r_1}(0))$ of $H^{1/2}$-functions of two variables that have support in the ball $B_{r_1}(0)$.

One restriction of our work is that we consider only potentials with
compact support: The reason is that we were forced to work with
different function spaces of solutions for $r$ small and for $r$
large. For $r$ small, we have some freedom in choosing the space, as
any space of the form $X=H^{1+\alpha}\times H^\alpha\times H^1\times
L^2$, with $0<\alpha\leq 1$, ensures that an exponential dichotomy
exists. For $r$ large, due to the structure of the equations, there is
no such freedom, the regularity on the first two components has to be
same as the regularity of the last two.  So it is unclear whether
there exists an exponential dichotomy when the support of $\rho$ is
not compact. It would be interesting to see whether our hypothesis
that $\rho$ has compact support could be replaced by an appropriate
decay condition on~$\rho$.

For the original potential $\theta$, we see no obstacles in removing the condition that $\theta$ has compact support. It should be possible to replace this condition by the long range condition $|\theta'(r)|\le C (1+ r)^{-1-\beta}$ for some $\beta>0$. It should also be
possible to remove the condition that $\theta$ is radially symmetric, although considerably more work will be needed without this condition.

We believe that the methods put forward in this paper can be used to study
other operators. In particular, the exponential-dichotomy results established
in \cite{aS98} are for systems of reaction-diffusion equations, so we believe that the only
obstacle for extending our results to selfadjoint systems are the presence of nonsmooth
potentials. For other operators, it might not be possible to modify the
function spaces involved to prove the existence of exponential dichotomies. These are
difficult problems that have to be studied in future work.

\end{section}



\end{document}